\newif\ifarxiv
\arxivtrue

\ifarxiv
\documentclass{article}
\usepackage[square,numbers]{natbib}
\usepackage{geometry}
\usepackage{sectsty}
\sectionfont{\fontsize{12}{12}\selectfont}
\subsectionfont{\fontsize{11}{11}\selectfont}
\subsubsectionfont{\fontsize{11}{11}\selectfont}
\else
\documentclass[3p]{elsarticle}
\journal{Journal of Computational Physics}
\fi

\usepackage[utf8]{inputenc}

\usepackage[font=small]{caption}
\usepackage[font=footnotesize]{subcaption}

\usepackage{amsmath}
\usepackage{url}
\usepackage{bm}

\usepackage{amssymb}
\usepackage{amsthm}
\usepackage{mathtools}
\usepackage{amsfonts}
\usepackage{microtype}
\usepackage{mathpazo}
\usepackage{yfonts}
\usepackage[normalem]{ulem}

\usepackage{graphicx}
\usepackage{afterpage}
\usepackage{setspace}
\usepackage{xcolor}
\usepackage{color}
\usepackage{floatrow}
\usepackage{isomath}
\usepackage[export]{adjustbox}
\usepackage{cancel}

\newfloatcommand{capbtabbox}{table}[][\FBwidth]
\usepackage{blindtext}

\usepackage{xifthen}
\usepackage{ifdraft}
\usepackage{paralist}                 
\usepackage[modulo]{lineno}

\makeatletter
\newcommand{\hypertargetraised}[1]{\Hy@raisedlink{\hypertarget{#1}{}}}
\makeatother

\usepackage{hyperref}
\usepackage[capitalize]{cleveref}
\crefname{equation}{Eq.}{Eqs.}
\Crefname{equation}{Equation}{Equations}
\crefrangelabelformat{equation}{(#3#1#4--#5#2#6)}
\crefmultiformat{equation}{Eqs. (#2#1#3}{, #2#1#3)}{#2#1#3}{#2#1#3}
\Crefmultiformat{equation}{Equations (#2#1#3}{, #2#1#3)}{#2#1#3}{#2#1#3}

\newcommand\pr[1]{\cref{#1}}

\newif\ifUseTikz
\UseTikztrue 
\ifUseTikz
\usepackage{tikz,pgfplots,import,tkz-euclide}
\pgfplotsset{compat=newest}

\usetikzlibrary{arrows.meta}
\usetikzlibrary{backgrounds}
\usetikzlibrary{external}
\usetikzlibrary{pgfplots.groupplots}
\usetikzlibrary{plotmarks}
\usetikzlibrary{shapes}

\pgfplotsset{plot coordinates/math parser=false}
\newlength\figureheight
\newlength\figurewidth

\pgfkeys{/pgf/images/include external/.code={\includegraphics[]{#1}}}

\makeatletter
\tikzset{
    every picture/.style={
        execute at begin picture={
            \let\ref\@refstar
        }
    }
}
\makeatother

\definecolor{plt-blue}{rgb}{0.0078,0.2980,0.7961}
\definecolor{plt-orange}{rgb}{1.0000,0.6431,0.2627}
\definecolor{plt-purple}{rgb}{1.0000,0.2863,0.5255}
\definecolor{plt-violet}{rgb}{0.6118,0.1765,1.0000}
\else 
\fi

\usepackage{tabularx}     
\usepackage{booktabs}     
\usepackage{multirow}
\usepackage{collcell}     
\usepackage{dcolumn}      

\usepackage{amsthm}

\usetikzlibrary{calc}                                       
\usetikzlibrary{decorations}                  
\newcounter{icoord}            
\pgfkeys{/tikz/.cd,                                                                                                                      
    curved pipe/.cd,                                                                                                                     
    radius/.store in=\CurvedPipeRadius,                                                                                                  
    radius=10pt,                                                                                                                         
    step/.store in=\CurvedPipeStep,
    step=1pt,                         
    shade/.style={left color=blue!22,right color=blue!37,middle color=blue!5}
}                                         
                                                                                                                                         
\pgfdeclaredecoration{curved pipe}{initial}
{%
  \state{initial}[width=\CurvedPipeStep,next state=cont] {%
    \pgfmoveto{\pgfpoint{\CurvedPipeStep}{\CurvedPipeRadius}}%
    \pgfpathlineto{\pgfpoint{0.3\pgflinewidth}{\CurvedPipeRadius}}%
    \setcounter{icoord}{0}%
    \pgfcoordinate{lastup-\number\value{icoord}}{\pgfpoint{1pt}{\CurvedPipeRadius}}%
    \pgfcoordinate{lastdown-\number\value{icoord}}{\pgfpoint{1pt}{-1*\CurvedPipeRadius}}%
  }                                                
  \state{cont}[width=\CurvedPipeStep]{%
     \stepcounter{icoord}%
    \pgfcoordinate{lastup-\number\value{icoord}}{\pgfpoint{\CurvedPipeStep}{\CurvedPipeRadius}}
     \pgfcoordinate{lastdown-\number\value{icoord}}{\pgfpoint{\CurvedPipeStep}{-1*\CurvedPipeRadius}}
     \pgfcoordinate{tmpup}{\pgfpoint{\CurvedPipeStep+0.3pt}{\CurvedPipeRadius}}
     \pgfcoordinate{tmpdown}{\pgfpoint{\CurvedPipeStep+0.3pt}{-1*\CurvedPipeRadius}}
     \pgfmathsetmacro\myshadingangle{0}%
     \path[curved pipe/shade,shading angle=\myshadingangle]                                                                              
        (lastup-\the\numexpr\value{icoord}-1) 
     -- (tmpup) to[out=180,in=180] (tmpdown) -- (lastdown-\the\numexpr\value{icoord}-1)
     to[out=180,in=180] cycle;%
     \pgfmoveto{\pgfpointanchor{lastup-\the\numexpr\value{icoord}-1}{center}}%
     \pgfpathlineto{\pgfpointanchor{lastup-\number\value{icoord}}{center}}%
     \pgfmoveto{\pgfpointanchor{lastdown-\the\numexpr\value{icoord}-1}{center}}%
     \pgfpathlineto{\pgfpointanchor{lastdown-\number\value{icoord}}{center}}%
  }                                                                                                                                      
  \state{final}[width=\CurvedPipeStep]                           
  { 
    \pgfmoveto{\pgfpointdecoratedpathlast}
    \fill (tmpup) to[out=0,in=0] (tmpdown) to[out=-180,in=-180] cycle;           
  }                                                                                                                                      
}

\definecolor{c1}{rgb}{.10588,.6197,.4666}
\definecolor{c2}{rgb}{.851,.3725,0}
\definecolor{c3}{rgb}{.4588,.4392,.702}

\usetikzlibrary{external}
\usepackage{expl3}
\ExplSyntaxOn
\newcommand\latinabbrev[1]{
  \peek_meaning:NTF . {
    #1\@}%
  { \peek_catcode:NTF a {
      #1.\@ }%
    {#1.\@}}}
\ExplSyntaxOff

\def\eg{\latinabbrev{e.g}}
\def\etal{\latinabbrev{et al}}

\def\ie{\latinabbrev{i.e}}

\newcommand\vct[1]{{\ensuremath{\bm{#1}}}}
\newcommand\convop[1]{{\ensuremath{\mathcal{#1}}}}
\newcommand\Real{{\ensuremath{\mathbb{R}}}}
\newcommand\bigO[1]{{\ensuremath{\mathcal{O}\!\left(#1\right)}}}   

\newcommand\radius{{\ensuremath{\varepsilon}}}
\newcommand\params{{\ensuremath{s}}}
\newcommand\paramt{{\ensuremath{\theta}}}
\newcommand\normal{{\vct{n}}}

\newcommand\LaplaceSLR{{\ensuremath{\convop{S}^{L}}}}  
\newcommand\LaplaceDLR{{\ensuremath{\convop{D}^{L}}}}

\newcommand\LaplaceSL{{\ensuremath{S^{L}}}}          
\newcommand\LaplaceDL{{\ensuremath{D^{L}}}}
\newcommand\LaplaceCF{{\ensuremath{K^{L}}}}
\newcommand\StokesSLR{{\ensuremath{\convop{S}}}}       
\newcommand\StokesDLR{{\ensuremath{\convop{D}}}}

\newcommand\StokesSL{{\ensuremath{S}}}              
\newcommand\StokesDL{{\ensuremath{D}}}
\newcommand\StokesCF{{\ensuremath{K}}}
\newcommand\GreensFn{\ensuremath{G}}         

\newcommand\RigidBodyVel{{\ensuremath{\vct{V}}}}
\newcommand\RigidBodySpace{{\ensuremath{\mathfrak{V}}}}
\newcommand\CompletionFlowDen{{\ensuremath{\vct{\nu}}}}

\newcommand{\Time}{{\ensuremath{\text{T}}}}                      
\newcommand{\Tsolve}{{\ensuremath{\Time_{solve}}}}
\newcommand{\Tsetup}{\ensuremath{\Time_{setup}}}
\newcommand{\Teval}{\ensuremath{\Time_{eval}}}

\newcommand\Ns[1][ ]{{\ensuremath{N_{s}^{#1}}}}        
\newcommand\Nt[1][ ]{{\ensuremath{N_{\theta}^{#1}}}}    
\newcommand\Nobj{{\ensuremath{B}}}             
\newcommand\Nelem{\ensuremath{K}}             

\newcommand{\Nunknown}{\ensuremath{N}}                           
\newcommand{\Nnodes}{\ensuremath{N_{nodes}}}                       
\newcommand{\gmrestol}{\ensuremath{\epsilon_{\textsc{GMRES}}}}   
\newcommand{\gmresiter}{\ensuremath{N_{\textsc{GMRES}}}}         
\newcommand{\quadtol}{\ensuremath{\epsilon_{\textrm{quad}}}}     

\newcommand{\SlenderElem}{\ensuremath{\Gamma}}      
\newcommand{\Near}{\ensuremath{\mathcal{N}}}

\newcommand\SurfJacobian{{\ensuremath{J}}}

\newcommand\ModalIntegOp{{\ensuremath{h}}}

\newcommand{\bi}{\begin{itemize}}
\newcommand{\ei}{\end{itemize}}
\newcommand{\ben}{\begin{enumerate}}
\newcommand{\een}{\end{enumerate}}
\newcommand{\be}{\begin{equation}}
\newcommand{\ee}{\end{equation}}
\newcommand\tbox[1]{{\textrm{#1}}}
\newtheorem{thm}{Theorem}

\newtheorem{pro}[thm]{Proposition}

\newtheorem{rmk}[thm]{Remark}

\newcommand{\eps}{\radius}         
\newcommand{\kk}{^{(k)}}         
\newcommand{\oa}{^{(\alpha)}}         
\newcommand\Nsk{\ensuremath{N_{\params}\kk}}            
\newcommand\Ntk{\ensuremath{N_{\paramt}\kk}}
\newcommand{\sk}{\params\kk}
\newcommand{\tk}{\paramt\kk}
\newcommand{\xx}{\vct{x}}
\newcommand{\yy}{\vct{y}}
\newcommand{\uu}{\vct{u}}

\newcommand{\su}{\yy}              
\newcommand{\bp}{\vct{\phi}}       
\newcommand{\shat}{\hat{\vct{s}}}     

\newcommand\NystCorrec{{\ensuremath{E}}}

\ifarxiv
\date{}
\else
\date{\today}
\fi

\begin{document}

\ifarxiv
\title{\Large\bf\sffamily Efficient Convergent Boundary Integral Methods for Slender Bodies}
\author{Dhairya Malhotra\footnote{Center for Computational Mathematics, Flatiron Institute, New York, NY 10010.
  Email: \texttt{dmalhotra@flatironinstitute.org}},\:
  Alex Barnett\footnote{Center for Computational Mathematics, Flatiron Institute, New York, NY 10010.
  Email: \texttt{abarnett@flatironinstitute.org}}}
\maketitle
\else
\begin{frontmatter} 
\title{\bf\sffamily Efficient Convergent Boundary Integral Methods for Slender Bodies}
\author{Dhairya Malhotra\footnote{Center for Computational Mathematics, Flatiron Institute, New York, NY 10010.
  Email: \texttt{dmalhotra@flatironinstitute.org}},
  Alex Barnett\footnote{Center for Computational Mathematics, Flatiron Institute, New York, NY 10010.
  Email: \texttt{abarnett@flatironinstitute.org}}}
\fi

\begin{abstract}
  %
  The interaction of
  fibers in a viscous (Stokes) fluid plays
  a crucial role in industrial and biological processes,
  such as sedimentation, rheology, transport, cell division, and locomotion.
%
%
Numerical simulations generally rely on 
slender body theory (SBT),
an asymptotic, nonconvergent approximation whose error blows up
as fibers approach each other.
%
Yet convergent
boundary integral equation (BIE) methods which completely resolve
the fiber surface have so far been impractical due to the prohibitive cost
of layer-potential quadratures in such high aspect-ratio 3D geometries.
We present a 
high-order Nystr\"om quadrature scheme
with aspect-ratio independent cost,
making such BIEs practical.
It combines centerline panels (each with a small number of poloidal Fourier modes),
toroidal Green's functions,
generalized Chebyshev quadratures, HPC parallel implementation, and FMM acceleration.
%
We also present new BIE formulations for slender bodies
that lead to well conditioned linear systems upon discretization.
%
We test Laplace and Stokes Dirichlet problems, and Stokes mobility problems,
for slender rigid closed fibers with (possibly varying) circular cross-section,
at separations down to $1/20$ of the slender radius,
reporting convergence typically to at least 10 digits.
%
We use this to quantify the breakdown of numerical SBT for
close-to-touching rigid fibers.
We also apply the methods to time-step the sedimentation
of 512 loops with up to $1.65$ million unknowns at around 7 digits of accuracy.
\end{abstract}


\ifarxiv
\else
\end{frontmatter} 
\fi

\section{Introduction}

Understanding viscous hydrodynamics in the presence of high aspect ratio bodies, such as
filaments, rods, and rings,
is key to many areas of science and industry,
including the rheology of fiber and polymer suspensions \cite{AML19rev}
and permeability of fiber structures \cite{tomadakis05}.
Numerical simulation is crucial, since theoretical approximations may
only apply in the dilute limit, or not at all \cite{guazzellirev}.
  For instance, the gravity induced sedimentation of rigid fibers
  (with applications to paper and pulp manufacture)
  is shaped by spatio-temporal
  correlations in density, velocity and orientation,
  and instabilities such as streamers and flocculation, that are an active research area
  \cite{saint05,Gustavsson09,guazzellibook,naz17}.
  Modern experimental methods allow a detailed comparison with numerics \cite{AML19rev}.
  Simulation also sheds light on many areas of
  biological fluids, including
locomotion by driven flagellae \cite{Fauci88},
the collective effects of such swimmers \cite{laugabook},
and transport by arrays of driven cilia \cite{Ding14}.
Molecular motors moving along a bed of flexible filaments
can lead to a swirling instability in oocytes \cite{stein21}.
In cells, the hydrodynamics of large assemblies of microtubules
or actin filaments
controls cell motion, transport, and division \cite{naz17}.

The inner step in modeling such phenomena is
solving a quasistatic elliptic boundary value problem (BVP) for the incompressible
Stokes equations in the exterior of the collection of slender bodies in question;
see \pr{e:stokes-eq0,e:stokes-eq1,e:stokes-eq2} below.
From this, the linear relationship between
velocities (and angular velocities) and forces (and torques) may be extracted.
In the case of rigid bodies, this is encoded by a matrix,
leading to two distinct tasks: \cite{KK91book,Pozrikidis1992}
1) the \textit{resistance problem} where velocities and angular velocities are
prescribed (as in porosity, where they vanish), and 2) its inverse, the
\textit{mobility problem} where body forces and torques are prescribed
(as in sedimentation).
In the former case statistics of the velocity field are
of interest \cite{tomadakis05,afk16}.
In the latter case the resulting velocities enable the body dynamics to be
evolved in a time-stepping scheme for a system of ordinary differential equations
(see, e.g., \cite{Tornberg2006,junwang}).
The case of flexible fibers is in some ways easier, since no BVP solution is needed
(simply an evaluation of prescribed Greens function sources);
however, more specialized time-stepping schemes are needed
to handle bending and inextensibility \cite{naz17,AML19rev}. 
Stochastic averaging over many expensive long-time simulations is often needed
to extract meaningful bulk quantities (viscosity, transport rates, etc),
making computational efficiency a pressing concern.
In this work we consider the rigid (and non-Brownian)
case where the challenge of an accurate BVP solution is foremost.
 
\subsection{Prior work}

We now overview numerical methods for the Stokes BVP and associated
mobility problem for many slender bodies;
also see \cite[Sec.~3]{AML19rev} in the flexible case.
Firstly, direct discretization
(e.g., finite elements), while in principle convergent as the mesh size
tends to zero, is inefficient and thus rarely used.
This is due to the high cost of meshing (or remeshing every time step)
the exterior domain and its inability to handle unbounded domains.
More common volumetric approaches include Lattice-Boltzmann methods \cite{Wu10}
in which fictitious gas particles interact within a grid of cells
covering the domain, and immersed boundary methods  where
fibers are overlayed on a finite difference grid \cite{limpeskin}.
However, both of these are accurate only when the grid spacing
is somewhat less than the smallest features or separations \cite{Wu10},
forcing the number of grid points to be huge if accuracy is needed.
For a comparison including finite elements for aggregates of spheres see \cite{schlauch13}.

The most popular approach is nonlocal \textit{slender body theory} (SBT)
\cite{keller76,johnson80,gotz00}.
This expresses the fluid velocity via a 1D integral over a given force density living on
the union of the centerlines of the bodies.
The integral kernel is the stokeslet Green's function,
plus a ``doublet'' correction due to Johnson \cite{johnson80}.
The self-interaction of a body requires a special form with local and nonlocal terms.
SBT was originally derived using matched asymptotics as $\eps\to 0$,
where $\eps$ is a slenderness parameter, with
leading neglected term $\bigO{\eps^2 \log \eps^{-1}}$ \cite{johnson80}.
Ellipsoid-like rounding of open fiber ends is assumed;
for closed fibers (as we consider in this work)
there is no such complication \cite{ueda}.
For flexible fibers the force is given locally by the geometry,
thus evolving the dynamics involves simply applying the 1D SBT integral operator
(modulo a tension solve) \cite{Tornberg2004}.
In contrast, for rigid fibers one must \textit{solve} a 1D integral equation for the
unknown force density \cite{Tornberg2006} (the so-called \textit{slender body
  inverse problem} \cite{mori21inv}).
There are other subtleties.
Since the classical self-interaction operator is in fact
logarithmically divergent \cite{gotz00,ueda} \cite[App.~B]{Tornberg2004},
{\em regularization} of the kernel is often needed for numerical stability.
Only recently has SBT been derived \cite{koens18} from boundary integral equations (discussed below).
Very recent rigorous error bounds, by Ohm and coworkers,
include ${\cal O}(\eps^{3/2})$ in the flexible
case (which required defining a new {\em slender body BVP} \cite{moriCPAM}),
and numerical tests give best-fit errors around ${\cal O}(\eps^{1.7})$
\cite[Sec.~5]{mitchell2022}.

On the computational side, SBT has been scaled to around $10^3$ fibers
on one workstation
\cite{Gustavsson09,naz17,maxian21} through the use of fast multipole (FMM) acceleration \cite{Greengard1987} or particle-mesh Ewald methods \cite{saint05,Tornberg2006,maxian21}.
Building on straight-fiber recurrences of Gustavsson--Tornberg \cite{Tornberg2006},
efficient quadratures for velocity evaluation near a fiber
have been developed \cite{afKlinteberg2020},
and product quadratures for fiber self-interaction \cite{tornbergprod}.

Yet, despite their wide use, classical and regularized
SBT do not give convergent numerical methods:
for any finite $\eps>0$,
the error in solving the desired Stokes BVP (or in the flexible case,
slender body BVP), does not vanish
as the centerline discretization becomes finer.
Especially problematic is the growth towards $\bigO{1}$ errors
as fibers approach $\bigO{\eps}$ or smaller separation;
even at a separation of $10\eps$, \cite[Fig.~8]{mitchell2022} shows a
maximum-norm centerline velocity error of about $0.003$, as $\eps\to 0$,
in the flexible loop case.
In \pr{ss:results-sbt} we explore the analogous breakdown
for a resistance problem in the rigid loop case. 
The cause is failure of the inner SBT asymptotic expansion when there
is another fiber within distance $\bigO{\eps}$.
This has led to 
various ad hoc methods to better handle the close-touching case,
such as added lubrication forces \cite{butler02,saint05} or blending
to the centerline \cite{maxian21}.
However, the upshot is that in any given fixed-slenderness simulation,
there is no way to vary a convergence parameter in order to assess the size of the
errors. The errors induced by SBT in most practical settings are unknown.

We should note that SBT is one in a family of
popular lower-order approximate models used in slender
hydrodynamic interactions, such as local SBT
\cite{butler02,saint05}, pairwise Rotne--Prage--Yamakawa tensors
in Stokesian dynamics \cite{brady1988stokesian},
regularized stokeslets (e.g. \cite{bouzarth11,Ding14}),
bead models (e.g. \cite{delmottebead}), and short-range approximations.
Such methods are reviewed in \cite{Tornberg2006,AML19rev}.
A common difficulty is that the regularization scale, or
bead diameter, must be of the same order as the radius $\eps$ for accuracy,
making such models arbitrarily expensive as $\eps\to 0$.

Finally, via potential theory with the free-space Green's functions
(stokeslet and/or stresslet) for the Stokes system,
a Stokes BVP can be reformulated as a boundary integral equation (BIE)
involving an unknown 
density on the surfaces of the bodies
\cite{Ladyzhenskaya,KK91book,Pozrikidis1992,HW}.
There are various choices in formulation for the
resistance problem (Stokes Dirichlet BVP) and for the
mobility problem for rigid bodies; we discuss these below.
The discretized BIE results in a dense linear system,
but with many less unknowns ($N$) than for volume discretization
to the same accuracy.
An iterative solution for the density often converges rapidly if a well-conditioned
formulation is chosen. Here, each matrix-vector multiplication may be accelerated
by a Stokes FMM \cite{tornberg2008fast,stkfmm} or Ewald method \cite{bagge22},
giving optimal $\bigO{N}$ or $\bigO{N\log N}$ complexity.
With the density found, solution evaluation can also benefit from such acceleration.
BIE methods have had much success for bodies of moderate aspect ratio,
including mobility simulations with nearly $10^5$ spheres \cite{yan20}.
Accurate and efficient {\em quadrature schemes} to discretize the
weakly-singular
integral operator 
is still an active research challenge,
especially when surfaces approach each other or evaluation is needed near to a surface.
The root cause is the $1/r$ singularity on-surface, and typically $1/r^2$ singularity
for surface-to-volume interactions.
%
While we cannot review all such high-order quadratures,
for 3D Stokes the major types are:
Galerkin on triangle patches \cite{Pozrikidis1992,sauterschwab,delia14},
and Nystr\"om methods including partitions of unity \cite{bruno01,ying06,Malhotra-jcp-2019},
triangular patches \cite{bremer3d},
local weight corrections on regular grids \cite{tlupova18,zeta3d},
quadrature by expansion (QBX) \cite{lineqbx,afk16},    
line extrapolation \cite{hedgehog},
and spectral schemes exploiting spherical harmonics
for spheres \cite{coronasphere,yan20}
or (by grid rotations) smooth deformations of spheres \cite{veera11}.

At high aspect ratio ($\eps\to 0$), conventional BIE quadratures suffer:
the requirement that the surface patches or grids retain a local $\bigO{1}$
aspect ratio forces the $N$ per body needed for fixed accuracy to grow like $\eps^{-1}$.
This is growth is illustrated by the
2nd-order slender body tests of Keaveny--Shelley \cite{keaveny11},
the 16th-order torus tests of Bremer--Gimbutas
\cite[Tbl.~3]{bremer3d}, and the rounded rods of aspect ratio 10 of
Bagge--Tornberg \cite[Sec.~6.3]{baggerods} via QBX.
A related difficulty is the need to replace any ``completion flow''---%
a term traditionally applied using a stokeslet and rotlet at a single interior point
\cite{Pozrikidis1992}---%
with centerline sources
whose number also grows like $\eps^{-1}$ \cite{keaveny11}
\cite[Sec.~7.2.1]{baggerods}.
These issues motivated the present work.
The method of fundamental solutions has also been
used for Stokes flows \cite[Ch.~7]{Pozrikidis1992}, but also suffers
at least as much as BIE at high aspect ratio.

\subsection{Contributions}

We present a convergent Nystr\"om discretization scheme
(that we dub CSBQ for ``convergent slender-body quadrature''),
that is high-order accurate for the Stokes BIE over a wide range
of aspect ratios $\eps^{-1}$, and has constant cost in the limit $\eps\to 0$.
It exploits the idea that the density, although highly anisotropic
with respect to the surface metric,
is smooth as a function of centerline ($s$) and angular ($\theta$) coordinates.
The centerline is discretized into panels, on which
the density is interpolated from a tensor product grid in $s$ and $\theta$
involving a small
number (usually $\sim10$) of uniform angular nodes.
To handle the weakly singular on-surface kernels,
precomputed custom quadratures are used both for
inner toroidal Green's function evaluations and for outer centerline integrals.
This minimizes the number of auxiliary nodes needed, boosting efficiency.
We show that, even at aspect ratio only 10, our set-up and solution times are 
$10\times$ faster than a recent optimized more general torus
solver (\pr{ss:results-biest}), and this ratio grows linearly with aspect ratio.
By interpolation onto new target-dependent panels, we provide
accurate velocity evaluations with a cost that is only logarithmic in distance.
By choosing parameters adaptively, CSBQ also
accurately handles close-to-touching fibers: we show 11-digit accuracy for
a separation of $\eps/20$.

Through our HPC implementation of these tools \cite{CSBQ-zenodo},
we envision that
accuracy may be controlled (via simple convergence studies)
in rigid-fiber Stokes simulations,
\textit{at a CPU cost comparable to that of a numerical SBT simulation}.
Towards this goal,
we show that typically 6-digit accuracy is achieved with a quadrature set-up throughput
of 20,000 unknowns per second per core, which is
comparable to a few Stokes FMM calls at this accuracy.
In settings where SBT is inaccurate, or has unknown accuracy,
such as the common case of
moderate aspect ratio ($10$--$10^2$)
and/or close-to-touching fibers where lubrication forces play a role,
CSBQ for the first time enables a reliable reference solution at an acceptable cost.

{
A key contribution is a layer
representation (an $\radius$-dependent admixture of single- and double-layer)
which leads to a $\bigO{1}$ condition number as $\radius\to 0$,
for both Laplace and Stokes BVPs (\pr{s:bie}).
As our numerical tests show, achieving reasonable
condition numbers and iteration counts would otherwise be impossible.
For the Stokes mobility problem, we use this in
a new projected combined field integral equation
\pr{e:new-mobility-bie}, and prove that the latter
is uniquely solvable (\pr{t:BIE}).
} 

Our secondary contributions are:
\bi
\item We use the presented tools to
  explore the breakdown of the SBT solution for resistance (drag)
  as a function of the slenderness $\eps$ and separation of two rigid tori;
  see \pr{ss:results-sbt}.
  For this we solve the SBT inverse problem with high order accuracy.
  This complements a recent study \cite[Sec.~5.3]{mitchell2022} in the flexible-fiber case.
\item
  In mobility problem tests (sedimentation of up to 512 rigid loops),
we combine 
our BVP solver with a high-order timestepper to accurately evolve the body dynamics,
approaching the time of first numerical collision.
We study parallel strong scaling for this problem.
  \ei

We expect that CSBQ could also much accelerate the
solution of the slender body BVP for the flexible fiber case,
for which a high-order accurate, but non-accelerated, BIE quadrature was
devised recently \cite[App.~B]{mitchell2022}.

\begin{rmk}
  We develop and test the presented quadratures in the scalar Laplace
  as well as the Stokes case.
  The former can be directly applied to the electrostatics of thin wires and loops.
  The techniques would also easily adapt to the 
  Helmholtz (and possibly Maxwell) equations with high aspect
  ratio geometries when the radius is subwavelength.
  These are crucial to the modeling of
  electromagnetic scattering or radiation from
  arbitrarily curved (or close-to-touching) wires or antennae,
  for which 1D integral equations are only exact in the
  isolated, straight, axisymmetric antenna case \cite{haslamwire}.
  \end{rmk}

\subsection{Limitations}

CSBQ as presented is limited to
the common case of bodies with circular cross-section, although their
radius may smoothly vary along the centerline.
This allows all $\theta$-integrals to be evaluated rapidly via
toroidal Green's functions, for which we precompute custom quadratures.
We present and test only on closed fibers (loops);
note that the open-fiber case has endpoint singularities that add an extra complication,
unless a parabolic radius scaling is assumed as in SBT
\cite{johnson80,gotz00,Tornberg2004}.
In order to tackle singular geometries in classical Laplace and Stokes BVPs,
we restrict to rigid body problems (resistance and mobility).
Thus our work is not (directly) applicable to the
recently-formulated slender-body BVP needed for flexible-fiber hydrodynamics
\cite{moriCPAM}.

Finally, we do not address the explicit prevention of collisions:
we solve the mathematical Stokes BVPs and the resulting rigid-body
dynamics until the bodies become too close to resolve (see \pr{r:collisions}).

\subsection{Organization of the paper}
The following \pr{s:problem-setup} defines the Stokes BVPs under study
and their layer potentials.
\pr{s:numerical-algo} presents the core quadrature techniques.
\pr{s:bie} the new well-conditioned scaling for slender-body
combined-field representations, for both Laplace and Stokes BVPs,
and motivates them via extracting operator eigenvalues on the straight periodic fiber.
Several different types of numerical tests are reported in \pr{s:results},
including complicated fibers, close-touching geometries,
a study of the breakdown of SBT in a close-to-touching mobility setting,
comparison with an existing solver,
and the sedimentation in time of many rigid tori.
We summarize and discuss open problems in \pr{s:concl}.
Three appendices organize the details of the proofs, generalized
Chebyshev quadrature construction, and the numerical SBT inverse problem.
In \pr{t:notation}, we list some frequently used symbols for easy reference.

\begin{table}[!htbp]
  \centering
  \caption{\label{t:notation} Index of frequently used symbols.}
  \resizebox{\textwidth}{!}{
  \begin{tabular}{l l}
    \toprule 
    Symbol                                                    & Description \\
    \midrule 
    $\yy(\params, \paramt)$                                   & surface parameterization along centerline \\
                                                              & and angular direction respectively \\
    $\vct{x}_c(\params)$                                        & slender-body centerline coordinates \\
    \radius(\params)                                          & cross-sectional radius \\
    \normal(\params, \paramt)                                 & outward unit surface normal \\
    \midrule 
    $I_k$                                                     & centerline panel interval in $\params$\\
    $\SlenderElem_k$                                          & slender body surface element \\
    $\Near_{\SlenderElem_k}$                                      & near region of $\SlenderElem_k$ \\
    \Nsk, \Ntk                                                & element discretization orders in $\params$ and $\paramt$ \\
    \Nobj                                                     & number of rigid bodies \\
    \Nelem                                                    & number of slender body elements \\
    \Nunknown                                                 & total number of unknowns \\
    \bottomrule 
  \end{tabular}
  \begin{tabular}{l l}
    \toprule 
    Symbol                                                    & Description \\
    \midrule 
    $\LaplaceSLR[\cdot]$, $\StokesSLR[\cdot]$                 & Laplace, Stokes single-layer potentials\\
    $\LaplaceDLR[\cdot]$, $\StokesDLR[\cdot]$                 & Laplace, Stokes double-layer potentials\\
    $\LaplaceSL$, $\StokesSL$, $\LaplaceDL$, $\StokesDL$      & Laplace, Stokes single- and double-layer\\
                                                              & boundary integral operators \\
    \midrule 
    $m_0$                                                       & number of generalized Chebyshev nodes\\
    \quadtol                                                  & accuracy of layer potential quadrature \\
    \gmrestol                                                 & tolerance for GMRES solve \\
    \gmresiter                                                & average number of iterations for each \\
                                                              & solve of the boundary integral equation \\
    \Tsetup                                                   & setup time for quadrature \\
    \Teval                                                    & evaluation time for quadrature \\
    \Tsolve                                                   & total solve time (excluding setup time) \\
    \bottomrule 
  \end{tabular}
  }
\end{table}


\section{Problem setup and notation\label{s:problem-setup}}
Consider $\Nobj$ distinct slender bodies $\Omega = \bigcup_{b=1}^{\Nobj} \Omega_b$, where $\Omega_b\subset\Real^3$.
Each body $\Omega_b$ is described by 
its centerline, restricted in this work to be a smooth closed curve (loop),
plus a function $\radius$
giving the radius of the circular cross section at each point on the centerline.
%
The bodies are suspended in a Stokesian (linear viscous) fluid
where the constant dynamic viscosity has been nondimensionalized to 1.
The fluid velocity $\vct{u}$ and pressure $p$ in the exterior of $\Omega$ are governed by the Stokes equations,
\begin{align}
    -\Delta \vct{u} + \nabla p &= \vct{0} & \text{in } \Real^3 \setminus \overline{\Omega} \label{e:stokes-eq0} , \\
    \nabla \cdot \vct{u} &= 0 & \text{in } \Real^3 \setminus \overline{\Omega} \label{e:stokes-eq1}.
\end{align}
Here \pr{e:stokes-eq0} and \pr{e:stokes-eq1} denote the momentum balance and
incompressibility constraint respectively.
In addition, we also assume that the fluid velocity at infinity decays to zero,
\begin{align}
  \vct{u}(\vct{x}) &\rightarrow \vct{0}
    & \text{as } |\vct{x}| \rightarrow \infty \label{e:stokes-eq2}.
\end{align}
In this work we consider the Stokes problem with Dirichlet boundary
conditions (as arises in the resistance setting where body motions are
known),
and the Stokes mobility problem (corresponding to solving the
unknown rigid motions of bodies with specified forces and torques).
The boundary conditions in each case are discussed below.

\paragraph{Stokes Dirichlet problem} The fluid phase $\Real^3 \setminus \overline{\Omega}$ satisfies the Stokes equations
\pr{e:stokes-eq0,e:stokes-eq1} with decay boundary conditions at infinity
\pr{e:stokes-eq2}, and specified fluid velocity on the geometry boundary $\partial\Omega$,
\begin{align}
  \vct{u} &= \vct{u}_0 & \text{on } \partial\Omega  \label{e:stokes-dirichlet3} .
\end{align}
The velocity field $\vct{u}$ in the fluid phase $\Real^3 \setminus \overline{\Omega}$ has a unique solution
\cite[p.~60--64]{Ladyzhenskaya}
which is to be determined.
A common application is to the resistance problem \cite[\S4.9]{Pozrikidis1992}:
one solves this BVP with $\vct{u}_0$ equal to minus
the value on $\partial\Omega$ of a given background flow.
Then $\vct{u}$ is interpreted as the change from the background
flow due to the presence of the bodies $\Omega$ with non-slip
boundary conditions.

%
%

\paragraph{Stokes mobility problem}
The bodies $\Omega_b$, $b=1,\dots,\Nobj$, are embedded in a viscous fluid satisfying the Stokes equations
\pr{e:stokes-eq0,e:stokes-eq1} with decay at infinity \pr{e:stokes-eq2}.
The $b$\textsuperscript{th} body has a net force $\vct{F}_b$ and
a net torque $\vct{T}_b$ acting about a fiduciary point $\vct{x}^c_b\in\mathbb{R}^3$.
The bodies undergo rigid body motions with velocity $\RigidBodyVel$ of the form
\begin{align}
  \RigidBodyVel(\vct{x}) &= \vct{v}_b + \vct{\omega}_b \times (\vct{x} - \vct{x}^c_b), & \text{for all } \vct{x} \in \Omega_b , \; b=1,\dots,\Nobj,
  \label{rbm}
\end{align}
where $\vct{v}_b$ is the unknown translational velocity and $\vct{\omega}_b$ is the unknown angular velocity of $\Omega_b$ about the point $\vct{x}^c_b$.
A slip velocity boundary condition $\vct{u}_s$ between the rigid bodies and the fluid is prescribed, so that
\begin{align}
  \vct{u} &= \RigidBodyVel + \vct{u}_s  & \text{on } \partial \Omega.
  \label{e:mobility-bc}    
\end{align}
We are given $\vct{u}_s$, and $\vct{x}^c_b$, $\vct{F}_b$, and $\vct{T}_b$
for each $b=1,\dots,\Nobj$.
The rigid body motion $\RigidBodyVel$ (\ie, $\vct{v}_b$ and $\vct{\omega}_b$ for each $b=1,\dots,\Nobj$), and the flow $\vct{u}$, are not known and must be determined.
Such a solution exists and is unique, as shown
via potential theory and Fredholm theory \cite[\S4.9]{Pozrikidis1992}
\cite[Secs.~2.5.1 and 4.6]{manasthesis} \cite{corona3dmob}.

\paragraph{Boundary integral equations}
Solutions to homogeneous linear constant coefficient elliptic PDEs
can be represented by layer potentials, \ie, convolution of a Green's function
for the PDE with a boundary density function.
For the Stokes equations, \pr{e:stokes-eq0,e:stokes-eq1,e:stokes-eq2},
the fluid velocity $\vct{u}$ can be represented as
\cite{Ladyzhenskaya,Pozrikidis1992,HW}
\begin{align}
\vct{u}(\vct{x}) = \StokesSLR[\vct{\sigma}](\vct{x})
&\coloneqq \int_{\partial\Omega} S(\vct{x}-\vct{y}) \vct{\sigma}(\vct{y}) dS_\vct{y} &
\text{for all } \vct{x} \in \Real^3 \setminus \overline{\Omega},
\label{e:stokes-sl}
\end{align}
where $S(\vct{r}) :=
\frac{1}{8\pi}\left( \frac{\vct{I}}{|\vct{r}|} + \frac{\vct{r} \vct{r}^{T}}{|\vct{r}|^3} \right)$
is the (tensor-valued) single-layer Stokes velocity kernel,
$\vct{\sigma}$ is an unknown vector-valued boundary density function,
and $dS$ is the surface area element.
By construction, $\vct{u}$ satisfies the PDE everywhere in $\Real^3 \setminus \overline{\Omega}$.
We will also need the double-layer velocity potential
\begin{align}
\vct{u}(\vct{x}) = \StokesDLR[\vct{\sigma}](\vct{x})
&\coloneqq \int_{\partial\Omega} D(\vct{x}-\vct{y};\vct{n}_\vct{y}) \vct{\sigma}(\vct{y}) dS_\vct{y} &
\text{for all } \vct{x} \in \Real^3 \setminus \overline{\Omega},
\label{e:stokes-dl}
\end{align}
where
$D(\vct{r};\vct{n}):= -\frac{3}{4\pi} (\vct{r}\cdot\vct{n}) \frac{\vct{r} \vct{r}^{T}}{|\vct{r}|^5}$
is the (tensor-valued) double-layer velocity kernel,
and $\vct{n}_\vct{y}$ is the outward-pointing unit surface normal
at $\vct{y}\in\partial\Omega$.

Applying the boundary conditions to a representation as above
(which may involve jump relations)
results in an integral equation for $\vct{\sigma}$,
which may be discretized to give a
linear system whose solution vector is
$\vct{\sigma}$ at the set of surface nodes.
To effectively use such integral equation methods,
we need layer potential representations that lead to
second-kind Fredholm integral equations, hence to
well-conditioned
linear systems in the discretized unknown boundary density.
We also need efficient quadrature schemes to evaluate the singular integrals
required for computing layer potentials.
We discuss the algorithms for evaluating layer potentials for slender bodies
in \pr{s:numerical-algo} and
give boundary integral equation (BIE) formulations
for the Stokes Dirichlet and mobility problems in \pr{s:bie}.

\section{Numerical algorithms\label{s:numerical-algo}}
In \pr{ss:disc}, we describe the discretization of the slender body surface
and then discuss the construction of boundary integral operators in \pr{ss:quad}.
In \pr{ss:parallel}, we give details of the parallel algorithms.

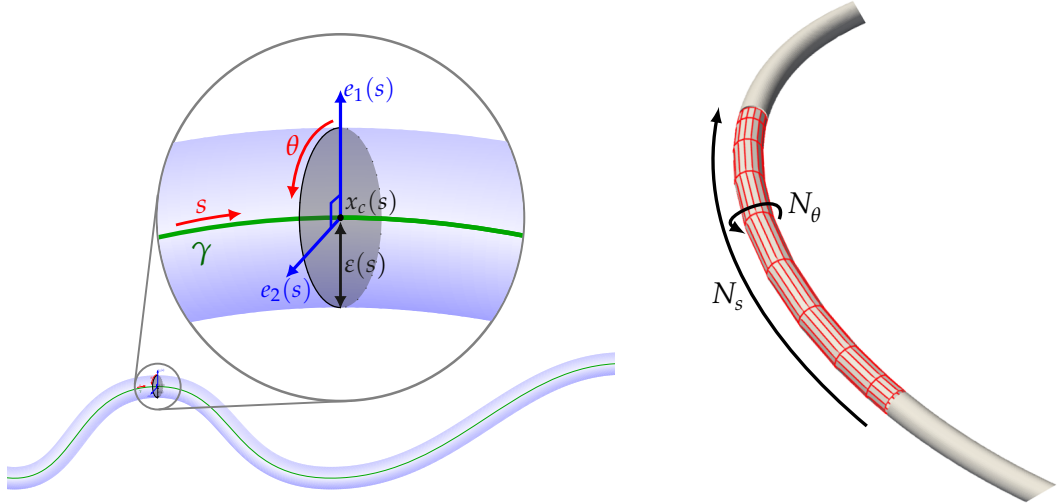
\begin{figure}[htb!]
  \centering
  $\vcenter{\hbox{\resizebox{0.49\textwidth}{!}{\begin{tikzpicture}
  \draw[white] (-1.28,-0.8) rectangle (12.0,9.7); 

  \begin{scope}[local bounding box=scope1]
    \clip (-1.28,-0.8) rectangle (12.0,2.5);

    \draw[decorate,decoration=curved pipe,curved pipe/radius=7.0,looseness=0.7]
     plot [smooth, tension=0.9] coordinates { (-3,1) (-1,-0.5) (2,1.5) (6,-0.5) (12,2) (15,0)}; 

    \draw[decorate,looseness=0.7, color=black!35!green]
     plot [smooth, tension=0.9] coordinates { (-3,1) (-1,-0.5) (2,1.5) (6,-0.5) (12,2) (15,0)}; 
  \end{scope}

  
  \begin{scope}[shift={(-2,-1.5)}, transform canvas={scale=8}, shift={(0.75,0.65)}] 
    \draw[color=gray, line width=0.2pt] (2,1.5) circle (0.502);
    \clip (2,1.5) circle (0.5);

    \draw[decorate,decoration=curved pipe,curved pipe/radius=7.0,looseness=0.7]
     plot [smooth, tension=0.9] coordinates { (-1,-0.5) (2,1.5) (6,-0.5)}; 

    \draw[line width=0.35pt, decorate,looseness=0.7, color=black!35!green]
     plot [smooth, tension=0.9] coordinates { (-1,-0.5) (2,1.5) (6,-0.5)}; 
    \node[color=black!60!green, scale=0.26] at (1.62,1.40) {$\gamma$};

    \draw[color=red, line width=0.22pt, -{Latex[length=1.0pt,width=1.0pt]}] (1.55,1.48)  to [out=14,in=188, looseness=1] (1.73,1.512); 
    \node[color=red, scale=0.22] at (1.62,1.53) {$\params$};

    \draw[color=red, line width=0.22pt, -{Latex[length=1.0pt,width=1.0pt]}] (1.98,1.765)  to [out=202,in=83, looseness=0.9] (1.866,1.55); 
    \node[color=red, scale=0.22] at (1.87,1.7) {$\paramt$};

    \draw[color=black, line width=0.22pt, {Latex[length=1.0pt,width=1.0pt]}-{Latex[length=1.0pt,width=1.0pt]}] (2.0,1.49)  to (2.0,1.25); 
    \node[color=black, scale=0.20] at (2.07,1.37) {$\radius(\params)$};

    \node[color=black, scale=0.19] at (2.085, 1.54) {$x_c(\params)$}; 

    \draw[fill=black, fill opacity=0.4, line width=0.0pt, darkgray, densely dotted] (2.0,1.5) ellipse (3.2pt and 7.02pt); 

    \begin{scope}
      \clip (1.0,1.0) rectangle (2.0,2.0);
      \draw[line width=0.1pt] (2.0,1.5) ellipse (3.2pt and 7.0pt);
    \end{scope}


    \draw[color=blue, line width=0.25pt, -{Latex[length=1.0pt,width=1.0pt]}] (2, 1.5) to +(0.0, 0.35); 
    \node[color=blue, scale=0.19] at (2.08, 1.85) {$e_1(\params)$};

    \draw[color=blue, line width=0.25pt, -{Latex[length=1.0pt,width=1.0pt]}] (2, 1.5) to +(-0.15, -0.165); 
    \draw[color=blue, line width=0.20pt] (1.974, 1.47) to (1.974, 1.54) to (2, 1.563); 
    \node[color=blue, scale=0.19] at (1.853, 1.3) {$e_2(\params)$};

    \draw[color=black, line width=0.25pt] (2, 1.5) circle (0.125pt); 

  \end{scope}
  \begin{scope}[shift={(-2,-1.5)}, transform canvas={scale=8}, shift={(0.75,0.65)}] 
    \clip (2,1.5) circle (0.5);
    \clip (2.0,1.0) rectangle (2.5,2.0);
    \draw[line width=0.35pt, decorate,looseness=0.7, color=black!35!green]
     plot [smooth, tension=0.9] coordinates { (-1,-0.5) (2,1.5) (6,-0.5)}; 
    \draw[color=black, line width=0.25pt] (2, 1.5) circle (0.125pt); 
  \end{scope}
  \draw[color=gray, line width=1.2pt] (2,1.5) circle (0.502);
  \draw[color=gray, line width=1.2pt] (1.497,1.55) -- (2.0,5.7);
  \draw[color=gray, line width=1.2pt] (2.02,1.0) -- (6.25,1.187);

    \node[color=black!60!green, scale=0.26] at (1.62,1.40) {$\gamma$};

    \draw[color=red, line width=0.22pt, -{Latex[length=1.0pt,width=1.0pt]}] (1.55,1.48)  to [out=14,in=188, looseness=1] (1.73,1.512); 
    \node[color=red, scale=0.22] at (1.62,1.53) {$\params$};

    \draw[color=red, line width=0.22pt, -{Latex[length=1.0pt,width=1.0pt]}] (1.98,1.765)  to [out=202,in=83, looseness=0.9] (1.866,1.55); 
    \node[color=red, scale=0.22] at (1.87,1.7) {$\paramt$};

    \draw[color=black, line width=0.22pt, {Latex[length=1.0pt,width=1.0pt]}-{Latex[length=1.0pt,width=1.0pt]}] (2.0,1.49)  to (2.0,1.25); 
    \node[color=black, scale=0.20] at (2.07,1.37) {$\radius(\params)$};

    \node[color=black, scale=0.19] at (2.085, 1.54) {$x_c(\params)$}; 

    \draw[fill=black, fill opacity=0.4, line width=0.0pt, darkgray, densely dotted] (2.0,1.5) ellipse (3.2pt and 7.02pt); 

    \begin{scope}
      \clip (1.0,1.0) rectangle (2.0,2.0);
      \draw[line width=0.1pt] (2.0,1.5) ellipse (3.2pt and 7.0pt);
    \end{scope}


    \draw[color=blue, line width=0.25pt, -{Latex[length=1.0pt,width=1.0pt]}] (2, 1.5) to +(0.0, 0.35); 
    \node[color=blue, scale=0.19] at (2.08, 1.85) {$e_1(\params)$};

    \draw[color=blue, line width=0.25pt, -{Latex[length=1.0pt,width=1.0pt]}] (2, 1.5) to +(-0.15, -0.165); 
    \draw[color=blue, line width=0.20pt] (1.974, 1.47) to (1.974, 1.54) to (2, 1.563); 
    \node[color=blue, scale=0.19] at (1.853, 1.3) {$e_2(\params)$};
    
    \draw[color=black, line width=0.25pt] (2, 1.5) circle (0.125pt); 
\end{tikzpicture}}}}$
  %
  \hspace{0.05\textwidth}
  $\vcenter{\hbox{\resizebox{0.29\textwidth}{!}{\begin{tikzpicture}
    \node[anchor=north west,inner sep=0] at (0,0) {\includegraphics[angle=0,width=5.5cm]{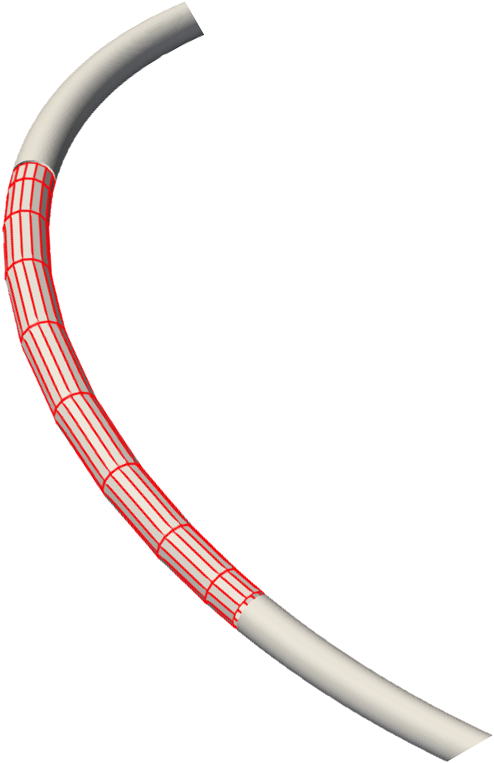}};
    \node at (1.25,-3.5) {\Large $\Nt$};
    \draw[ultra thick, -latex] (0.82,-3.7)  to [out=70,in=150, looseness=2.5] (0.25,-3.95);

    \node at (-0.05,-5) {\Large $\Ns$};
    \draw[ultra thick, -latex] (2.3, -7.2)  to [out=140,in=260, looseness=0.95] (-0.2, -1.8);
  \end{tikzpicture}}}}$
  \caption{\label{f:elem-disc}
    Left: A slender body is described by a centerline $\gamma$ parameterized
    as $\vct{x}_c(\params)$,
    and cross-sectional radius $\radius(\params)$.
    The vectors $d\vct{x}_c/ds$, $\vct{e}_1(\params)$ and $\vct{e}_2(\params)$
    are orthogonal, with $\vct{e}_1(\params)$
    setting the $\paramt$ angular origin at each $\params$.
  ~Right: The slender body surface
  $\Gamma$ ($=\partial\Omega_b$ for the $b^\text{th}$ body) is partitioned into
  surface elements $\Gamma_k$, $k=1,\dots,\Nelem$, each discretized
  by the tensor product of $\Nsk$ Gauss--Legendre ``panel'' nodes in $\params$
  with $\Ntk$ equispaced nodes in $\paramt$.
  }
\end{figure}

\subsection{Surface discretization \label{ss:disc}}

The geometry of the $b$\textsuperscript{th} slender body $\Omega_b$ is
shown in \pr{f:elem-disc} (left);
for notational simplicity we drop the dependence on $b$.
\pr{t:notation} provides a summary of commonly used notations.
The smooth closed centerline curve $\gamma$ is
parameterized by $\params \in [0,1]$, so that
$\gamma = \vct{x}_c([0,1))$.
  The circular cross-sectional radius is given by the smooth
  periodic function $\radius(\params)$.
  In addition, we require an orientation vector $\vct{e}_1(\params)$ of unit length at each point on the centerline,
orthogonal to the centerline. This orientation vector may be given or determined
automatically as described below in \pr{sss:orient}. A second orientation unit vector
$\vct{e}_2$ is computed as,
$
  \vct{e}_2 = \left. \frac{d\vct{x}_c}{d\params} \middle/ \left| \frac{d\vct{x}_c}{d\params}  \right| \right. \times \vct{e}_1.
$
  Points $\su$ on the slender body surface $\Gamma = \partial \Omega_b$
  are then given by,
\begin{equation}
    \su(\params,\paramt) = \vct{x_c}(\params) + \radius(\params)
    [\vct{e}_1(\params)\cos\paramt + \vct{e}_2(\params)\sin\paramt],
    \qquad \params\in[0,1),\; \paramt\in[0,2\pi).
    \label{e:surface-coord}
\end{equation}

We discretize a slender body surface by
constructing a piecewise polynomial approximation
of $\vct{x}_c$, $\radius$, and $\vct{e}_1$ along the centerline to a given accuracy tolerance.
This results in partitioning of $[0,1]$ into $\Nelem$ intervals
(or ``panels'') $\{I_1,\dots,I_\Nelem\}$,
with corresponding partition of $\Gamma$ into surface elements $\{\SlenderElem_1, \dots, \SlenderElem_\Nelem\}$.
A surface
element $\SlenderElem_k$ is approximated to
polynomial order (number of nodes)
$\Nsk$ in $\params$, and Fourier
discretization order (angular nodes) $\Ntk$ in $\paramt$.
Its nodes are thus
$\su\kk_{ij} = \su(\sk_i, \tk_j)$, where
$\sk_i$ for $i=1,\dots,\Nsk$
are Gauss--Legendre (GL) nodes for the interval $I_k$,
while $\tk_j$ for $j=1,\dots,\Ntk$
are uniform nodes in $[0,2\pi)$; see \pr{f:elem-disc} (right).
The total number of surface nodes for this body is thus
$\Nnodes := \sum_{k=1}^{\Nelem} \Nsk \Ntk$.

Smooth functions $f(\params,\paramt)$ on the slender body surface $\partial\Omega$
are represented to high-order accuracy by their point values at these nodes,
$f\kk_{ij} = f(\sk_i, \tk_j)$.
For instance, approximating a vector-valued density on the single body
requires $3\Nnodes$ unknowns.
A function can then be evaluated
at any another parameter pair $(\params,\paramt)$
by using barycentric Lagrange interpolation
from the nodes of the containing panel in $\params$ \cite{Berrut2004}
and trigonometric interpolation in $\paramt$.
By differentiating the interpolant, we can also compute surface gradients
of these functions.

\subsubsection{Constructing an orientation vector automatically\label{sss:orient}}
The orientation vector $\vct{e}_1$ can be any smoothly varying vector
defined at each point on the centerline and orthonormal to it.
For a given slender body $\Omega_b$, we set $\vct{e}_1$ at $\params=0$ to be any random vector
orthonormalized with respect to the centerline tangent vector ${d\vct{x}_c}/{ds}$.
Then, using this as the initial value, we solve the following ordinary differential equation (ODE) in $s\in(0,1)$,
\[
\frac{d\vct{e}_1}{ds} = \left. \vct{e}_1 \times \left( \frac{d^2\vct{x}_c}{ds^2} \times
\frac{d\vct{x}_c}{ds} \right) \middle/ \left| \frac{d\vct{x}_c}{ds}  \right|^2 \right. .
\]
{Since ${d\vct{e}_1}/{ds}$ is orthogonal to $\vct{e}_1$, the magnitude $|\vct{e}_1|$ remains constant for all $s$.
From the above expression for ${d\vct{e}_1}/{ds}$, we can show that $\frac{d}{ds}(\vct{e}_1 \cdot {d\vct{x}_c}/{ds}) = 0$ whenever $\vct{e}_1$ and ${d\vct{x}_c}/{ds}$ are orthogonal.
Since we choose $\vct{e}_1$ orthogonal to the centerline at $s=0$, it remains orthogonal to the centerline for all $s$.
}
We solve for $\vct{e}_1$ using spectral discretization of the ODE
on the Gauss--Legendre nodes of each panel in $\params$.
(Note that the accuracy of the scheme is not crucial, only that
the approximate solution be smooth.)
The derivatives \,${d\vct{x}_c}/{ds}$\, and \,${d^2\vct{x}_c}/{ds^2}$\, are computed through numerical differentiation
of the polynomial representation of $\vct{x}_c$ on the piecewise centerline panels.
Since the numerical solution and its piecewise polynomial representation along the centerline are not exact,
each time we evaluate $\vct{e}_1(\params)$, we re-orthonormalize it with respect to the centerline.

\subsection{Layer potential operators \label{ss:quad}}
We now describe the numerical computation of boundary integrals of the form in \pr{e:stokes-sl,e:stokes-dl}.
Our method is kernel independent, and therefore can be applied to most elliptic PDE kernels.
In this work we demonstrate it for Laplace single- and double-layer potentials
($\LaplaceSLR$ and $\LaplaceDLR$), and Stokes single- and double-layer
velocity potentials
($\StokesSLR$ and $\StokesDLR$).
The target point $\vct{x}$ can be either on-surface (such as when solving the boundary integral equation) or off-surface.
For on-surface targets, the boundary integral is singular and we need special quadrature rules.
For off-surface targets, when the targets are far away from the boundary, the integral is smooth
and we can use standard quadratures; however, for targets close to the boundary the integrand
is sharply peaked and using standard quadrature rules is not feasible.
In the following sections, we describe how our algorithm handles each of these cases.


\subsubsection{Nystr\"{o}m discretization with local corrections \label{sss:nystrom}}
We consider a single slender body surface $\Gamma=\partial\Omega_b$,
discretized into $\Nelem$
surface elements $\{\SlenderElem_1, \dots, \SlenderElem_\Nelem\}$,
as described in \pr{ss:disc}.
Also consider a surface density function $\sigma$ discretized similarly with the function values
$\sigma\kk_{ij}$ given at the surface discretization nodes $\su\kk_{ij}$.
Let $\GreensFn$ be a kernel function whose convolution with $\sigma$ on
$\Gamma$ gives the desired potential.
The potential at a target point $\vct{x}$ is given by summing the potential from each slender body element $\SlenderElem_k \in \partial\Omega$,
\begin{align}
u(\vct{x})
=
\int_{\Gamma} \GreensFn(\vct{x}, \vct{y}) \sigma(\vct{y}) dS_\vct{y}
= \sum\limits_{k=1}^{\Nelem} \int_{\SlenderElem_k} \GreensFn(\vct{x}, \vct{y}) \sigma(\vct{y}) dS_\vct{y}.
\end{align}
The kernel $\GreensFn(\vct{x},\vct{y})$     
can be any linear combination of single-layer $S(\vct{x}-\vct{y})$
and double-layer $D(\vct{x}-\vct{y}; \vct{n}_\vct{y})$ kernels.
The potential from each element $\SlenderElem_k$ is given by integrating along
its interval $I_k$ in $\params$ and over $\paramt \in [0,2\pi)$.
When the target $\xx$ is sufficiently far from $\SlenderElem_k$, the integrand is smooth
and the integral can be computed numerically using the existing Gauss--Legendre quadrature rule in $\params$ and periodic trapezoidal rule in $\paramt$,
\begin{align}
u(\vct{x})
=& \sum\limits_{k=1}^{\Nelem} \int_{I_k} \int_{0}^{2\pi} \GreensFn(\xx, \su(\params,\paramt)) \, \sigma(\params,\paramt) \, \SurfJacobian(\params,\paramt) \, d\paramt d\params
\label{e:elem-integral-s-t} \\
\approx& \sum\limits_{k=1}^{\Nelem} \sum\limits_{i=1}^{\Nsk} \sum\limits_{j=1}^{\Ntk}
\, \GreensFn(\xx,\su\kk_{ij}) w\kk_{ij} \, \sigma\kk_{ij}
\label{e:nystrom-no-correction}
\end{align}
where $\displaystyle \SurfJacobian = \left| \frac{\partial \su}{\partial\params} \times \frac{\partial \su}{\partial\paramt} \right|$ is the Jacobian
of \pr{e:surface-coord}, and $w\kk_{ij}$ are the quadrature weights.
These weights are defined as $w\kk_{ij} = \SurfJacobian\kk_{ij} w^{(GL,k)}_i \cdot 2 \pi/\Ntk$, where
$\SurfJacobian\kk_{ij} \coloneqq \SurfJacobian(\sk_i, \tk_j)$ are the values at the surface discretization nodes,
$w^{(GL,k)}_i$ are the weights for the $\Nsk$-th order Gauss--Legendre quadrature rule
on $I_k$, and
$2\pi/\Ntk$ are the equal periodic trapezoidal quadrature rule weights.
Evaluating \pr{e:nystrom-no-correction} requires \bigO{\Nnodes} work per target point, where $\Nnodes$ is the number of surface discretization points.
When the number of targets is also \bigO{\Nnodes}, this requires \bigO{\Nnodes^2} total work.
This can be accelerated through the use of a fast multipole method (FMM) \cite{Rokhlin1985,Greengard1987} and computed in $\bigO{\Nnodes}$ time.
We use the PVFMM library \cite{Malhotra2015,Malhotra2016}, which is an optimized distributed-memory implementation of the kernel independent FMM \cite{Ying2004}
and supports several elliptic PDE kernels.

\begin{figure}[htb!]
  \resizebox{0.7\textwidth}{!}{\begin{tikzpicture}
    \node[anchor=south west,inner sep=0] at (0,0) 
    {\includegraphics[width=10cm]{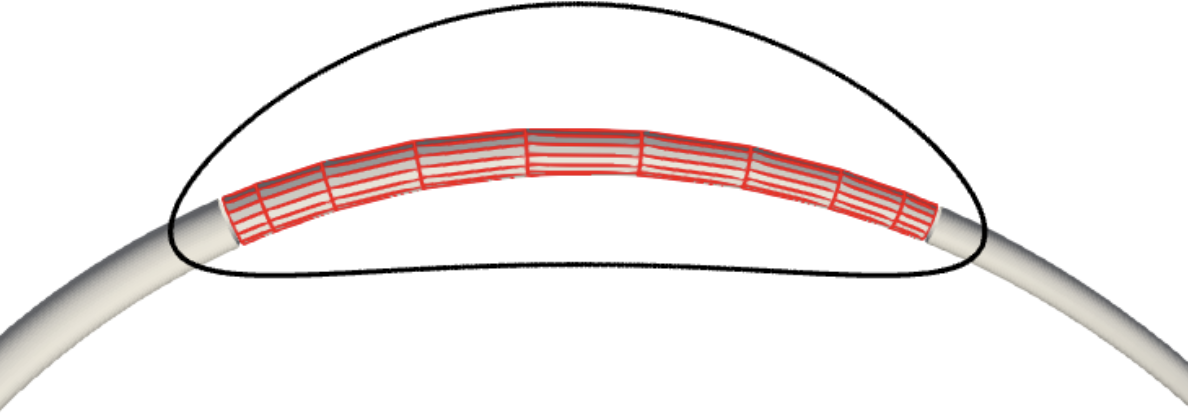}};
    \node at (4.7, 3.1) {$\Near_{\SlenderElem_k}$};

    \node at (6.0, 2.45) {\textcolor{red}{$\SlenderElem_k$}};
    \draw [fill=orange, fill opacity=0.35] (2.23,1.72) circle (0.79cm);
    \draw [fill=orange, fill opacity=0.35] (2.76,1.91) circle (0.96cm);
    \draw [fill=orange, fill opacity=0.35] (3.50,2.10) circle (1.10cm);
  \end{tikzpicture}}
  \caption{\label{f:bernstein-ellipse}
  The slender body surface is discretized into elements $\SlenderElem_k$.
  For a given quadrature accuracy tolerance $\quadtol$, there is a region
  $\Near_{\SlenderElem_k}$ around each element such that for target points outside this region,
  standard quadrature rules (tensor product of Gauss--Legendre and periodic trapezoidal) can be used.
  For target points inside $\Near_{\SlenderElem_k}$ we use special quadrature rules.
  We approximate the region $\Near_{\SlenderElem_k}$ by a set of overlapping spheres centered at the surface discretization nodes.
  This simplifies the task of identifying the target points within $\Near_{\SlenderElem_k}$.
  }
\end{figure}
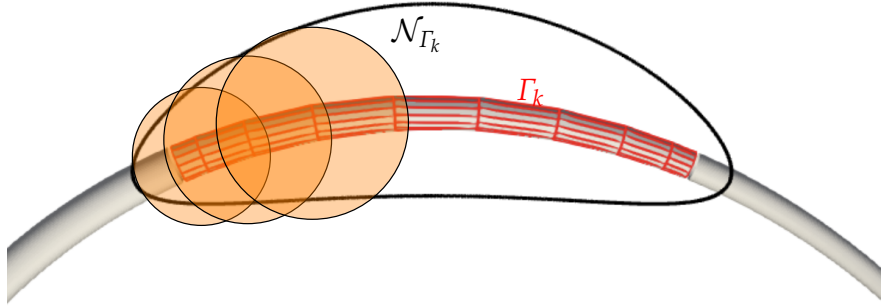

Error estimates for Gauss--Legendre quadrature are given in \cite[Thm.~19.3]{Trefethen2019},
and for periodic trapezoidal quadrature in \cite{Trefethen2014}.
These estimates depend on how far the integrand can be extended analytically into the complex
parameter plane.
In our application, this is related to the distance of the target point from the surface.
Given an accuracy tolerance $\quadtol$,
the $k$\textsuperscript{th} element quadrature in \pr{e:nystrom-no-correction}
is accurate for target points outside a region $\Near_{\SlenderElem_k}$
around $\SlenderElem_k$, as shown in \pr{f:bernstein-ellipse}.
Such regions have been determined accurately
for quadratures on curves in \cite{afKlinteberg2020} and
for quadratures on surfaces in \cite{afKlinteberg2022}.
In our setting, when the radius $\radius$ is small, $\Near_{\SlenderElem_k}$
takes the form of a prolate ellipsoid
(a body of revolution generated by a Bernstein
ellipse) distorted by the map $\xx_c(I_k)$ giving the centerline panel
\cite[Sec.~3.2]{afKlinteberg2020}.
For efficiency, we approximate this $\Near_{\SlenderElem_k}$ 
by a union of overlapping spheres centered at each surface discretization node.
The radius of the spheres is chosen based on
the accuracy tolerance $\quadtol$, the dimensions of the element,
and the orders of the Gauss--Legendre and periodic trapezoidal quadrature rules.
%
By sorting the target points and these spheres on a space-filling curve
we can efficiently determine the target points which lie within
each neighborhood $\Near_{\SlenderElem_k}$
for all elements $\SlenderElem_k$ in a collection of slender bodies.
We describe this algorithm in more detail in \pr{ss:parallel}.
For all such pairs of elements and target points $\vct{x} \in \Near_{\SlenderElem_k}$,
we add a correction to \pr{e:nystrom-no-correction},
\begin{align}
u(\vct{x}) \approx& \sum\limits_{k=1}^{\Nelem} \sum\limits_{i=1}^{\Nsk} \sum\limits_{j=1}^{\Ntk} \, \GreensFn(\xx,\su\kk_{ij}) w\kk_{ij} \, \sigma\kk_{ij}
+
%
%
\sum\limits_{\bigl\{k : \, \vct{x}\in\Near_{\SlenderElem_k} \bigr\}}
\sum\limits_{i=1}^{\Nsk} \sum\limits_{j=1}^{\Ntk} \NystCorrec_{ij}\kk(\xx) \, \sigma\kk_{ij}
. \label{e:nystrom-correction}
\end{align}
where $\NystCorrec_{ij}\kk(\vct{x})$ are
target-specific corrections to the discretized linear operator.
{
The target $\vct{x}$ can be either
an off-surface point anywhere in $\Near_{\SlenderElem_k}$,
or it can be an on-surface point, in which case it is one of the surface discretization nodes of $\SlenderElem_k$.
We discuss the construction of $\NystCorrec_{ij}\kk(\xx)$ for each of these two cases in \pr{sss:near-sing-quad} and
\pr{sss:sing-quad} respectively.
}
These corrections have the effect of subtracting the incorrect first term above
and adding the correct potential.
The correct potential is computed using special target-specific quadrature rules described below.
These quadratures can be relatively expensive to evaluate on-the-fly each iteration.
Therefore, for every surface element $\SlenderElem_k$ we precompute
the needed correction matrix elements $\NystCorrec\kk_{ij}(\su^{(k')}_{i'j'})$
for all surface nodes $\su^{(k')}_{i'j'}$ in its near field.
The result may be kept in memory as a sparse $\Nnodes\times\Nnodes$ matrix,
whose product with the current
density vector is added to the result of the FMM.
This makes subsequent applications of the
discretized boundary integral operator faster.
Since the corrections are local, we expect only $\bigO{1}$ storage and evaluation cost per target point.
We refer to the computation of these corrections as the quadrature setup step,
and the subsequent application of the integral operator as the evaluation step.

\begin{figure}[htb!]
  \resizebox{.7\textwidth}{!}{\begin{tikzpicture}
    \path [draw=none,fill=white!0,even odd rule] (4,0) circle (1.5) (4,0) circle (0.75);
    \path [draw=none,fill=blue!30,even odd rule] (4,0) circle (1.5) (4,0) circle (0.75);
    \draw[color=red, ultra thick] (0,0) ellipse (4cm and 1cm);
    \draw[fill=blue!30,draw=none] (180:0.75)+(4,0) arc (180:0:0.75) -- (0:1.5)+(4,0) arc (0:180:1.5) -- cycle;

    \draw[color=blue, ultra thick] (4.7,-0.8) circle (0.7pt);
    \node at (5, -0.8) {\color{blue} \large $\vct{x}_k$};
    \node at (-2.5, 0.38) {\color{red} \Large $\gamma_{circ}$};
    \node at (0, -0.5) {\Large $\paramt$};
    \node (c) at (-1, -0.75) {};
    \node (d) at ( 1, -0.75) {};
    \draw[ultra thick, ->] (c)  to [out=-5,in=185, looseness=1] (d);
    \draw[thick, -latex] (4, 0) to (4.75, 1.3);
    \draw[thick, -latex] (4, 0) to (4.7, 0.27);
    \node at (4.4,0.0) {\large $\alpha$};
    \node at (4.75,0.75) {\large $2\alpha$};
  \end{tikzpicture}}
  \caption{\label{f:modal-greens-fn}
  A circular source loop $\gamma_{circ}$ (in red) of unit radius,
  parameterized by $\paramt$ with points on the loop given by $\su_{circ}(\paramt)$.
  An annular sheet (in blue) with radius between $\alpha$ and $2\alpha$, and a target point $\vct{x}_k$ on this annulus.
  We use the methods described in \pr{sss:modal-greens-fn} and \pr{s:cheb-quad} to build a generalized Chebyshev quadrature rule in $\paramt$
  that can integrate the potential at any target point on the annulus.
  }
\end{figure}
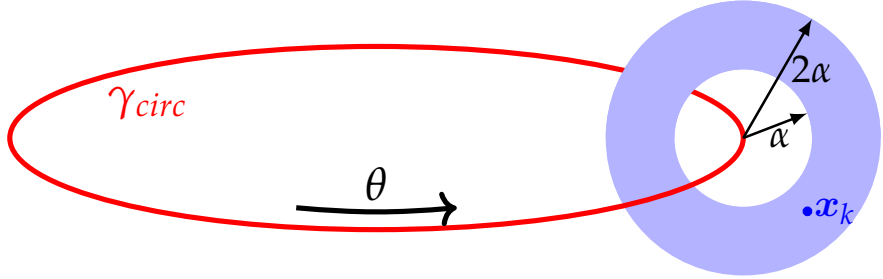

\subsubsection{Fast evaluation of angular integrals \label{sss:modal-greens-fn}}

{
In order to precompute the correction weights $E\kk_{ij}(\xx)$ appearing
in \eqref{e:nystrom-correction}, one must approximate the layer potential
\eqref{e:elem-integral-s-t} for targets $\xx$ close to the $k$th panel.
For this we use iterated integration, with the angular $\theta$-integral
performed first.
In this section we explain how, given $\xx$ and a single parameter
$\params$, we evaluate efficiently this angular integral
\begin{align}
I(\xx, \params) \coloneqq &
  \int_{0}^{2\pi} \GreensFn(\xx, \su(\params,\paramt)) \, \sigma(\params,\paramt) ~ \SurfJacobian(\params,\paramt) ~ d\paramt
  .
  \label{Idef}
\end{align}
For targets $\xx$ further than one radius from the source circle
$\su(\params,\paramt)$, $\paramt\in[0,2\pi)$,
simple use of the existing periodic trapezoidal rule nodes $\paramt_j$
is accurate, and no precomputation is necessary.
Otherwise, for closer targets
we precompute a special set of target-specific quadrature
weights $h_j(\xx,s)$, $j=1,\dots,\Nt$,
so that \eqref{Idef} is accurately approximated by the sum
\begin{align}
I(\xx, \params) 
  \approx &
  \sum\limits_{j=1}^{\Nt} \ModalIntegOp_j({\vct{x}}, \params) \, \sigma(\params, \paramt_j),
  \label{hwei}
\end{align}
involving the density $\sigma(s,\paramt_j)$ sampled only at the
periodic trapezoid nodes.
In the following subsections these weights $h_j(\xx,s)$,
for particular choices of nodes $\params$,
will be used to approximate the outer $\params$-integral over each
panel $I_k$, in order to calculate $E\kk_{ij}(\xx)$.
}

{
We now turn to the rapid calculation
of $h_j(\xx,s)$
for targets $\xx$ arbitrarily close to the source circle.
Here, a conventional adaptive quadrature in $\paramt$ would
be too slow, because the application to
singular on-surface $\params$-quadratures involves many such nearby targets.
Thus, for efficiency, we introduce another level of precomputation:
we precompute custom quadrature schemes in $\paramt$ that handle the
full range of near target locations using a small number of nodes.
To explain this, consider for now a standardized
unit radius source circle, denoted by $\gamma_{circ}$, parameterized as
$\su_{circ}(\paramt)$ such that the target point lies nearest the point
$\theta=0$; see \pr{f:modal-greens-fn}.
}
We want to compute the potential from this source at target points $\vct{x}$ close to it,
by evaluating integrals of the form $\int_0^{2\pi} \GreensFn(\xx,\su_{circ}(\paramt))f(\paramt)\, d\paramt$, where $f$ is smooth
(the product of the density and the Jacobian).
When $f(\paramt)=e^{i l \paramt}$, these potentials are called modal or toroidal
Green's functions \cite{youngmart}.

We construct a custom quadrature rule accurate for
all targets in an annular region with radius between $\alpha$ and $2\alpha$,
centered at $\su_{circ}(0)$ and orthogonal to the circular source at $\paramt=0$;
see \pr{f:modal-greens-fn}.
We sample a uniform distribution (\ie, a polar grid uniform in angle and radius)
of target points $\vct{x}_k$ for $k = 1, \dots, k_0$ on the annulus,
and define modal integrand functions $\varphi_{kl}(\paramt) \coloneqq
\GreensFn(\xx_k,\su_{circ}(\paramt)) e^{i l \paramt}$ for
$l \in \{0, 1,\dots, l_0\}$.
{
The maximum mode number $l_0$ is set slightly larger than $\Nt/2$
to account for the product with the Jacobian.
}
We choose the number of points $k_0$ to be sufficiently large
so that the integrand corresponding to a target point anywhere on the annulus
can be represented as a linear combination of $\varphi_{kl}$
up to the desired accuracy tolerance.
Typically, $k_0 \approx 600$ is sufficient for 14-digits accuracy.
Next, we build a panel based (composite)
Gauss--Legendre quadrature rule in the interval $[-\pi,\pi]$
with the panels refined dyadically around $\paramt=0$ until
the smallest panel size is commensurate with the dimensions of the annulus.
The order of the Gauss--Legendre rule is chosen to be sufficiently high to
integrate the product of any two integrands;
{this is required for constructing generalized Chebyshev quadrature rules in the next step.}
We follow the algorithm in \pr{s:cheb-quad}
to build a generalized Chebyshev quadrature rule
with nodes $\{\paramt\oa_1, \dots, \paramt\oa_{m_0}\}$
and weights $\{\omega\oa_1, \dots, \omega\oa_{m_0}\}$
that can integrate all integrands $\varphi_{kl}$.
This quadrature rule is thus accurate for evaluating the
potential from $\gamma_{circ}$
for all source functions $f(\paramt) = e^{i l \paramt}$ for $|l| \leq l_0$
and any target point $\vct{x}$ on the annulus.
We construct such a quadrature rule for each value of $\alpha \in \{2^{-1}, \dots, 2^{-30}\}$.
The number of quadrature nodes $m_0$ is found to be independent of the distance $\alpha$,
and is significantly smaller than what would be required for
an adaptive panel based quadrature rule.
For instance, for $l_0 = 8$ and quadrature accuracy tolerance $\quadtol = 10^{-10}$,
we have $m_0 \approx 38$ (for all distances $\alpha$).

In our application, the kernel function $\GreensFn$ is scale invariant; therefore,
these precomputed generalized Chebyshev quadratures can be applied to source loops of any radius
when the annulus region is also appropriately scaled.
Due to rotational symmetry,
by rotating in $\paramt$ by an appropriate angle, the quadrature can be applied to targets anywhere in the volume of revolution of the annulus.
To compute the inner integral in \pr{e:elem-integral-s-t}, for a given $\params$ and target $\vct{x}$,
we determine $\paramt_0$ corresponding to the closest point $\su(\params, \paramt_0)$ to $\vct{x}$.
We pick the quadrature rule with nodes $\{\paramt\oa_1, \dots, \paramt\oa_{m_0}\}$
and weights $\{\omega\oa_1, \dots, \omega\oa_{m_0}\}$\, such that
$\alpha \leq |\su(\params,\paramt_0)-\vct{x}| / \radius(\params) < 2 \alpha$.
Then, the integral \eqref{Idef} can be approximated as
\begin{align*}
I(\xx,\params)
\approx&
  \sum\limits_{m=1}^{m_0} \GreensFn(\xx,\su(\params,\paramt_0 + \paramt\oa_m)) ~ \sigma(\params,\paramt_0 + \paramt\oa_m) ~ \SurfJacobian(\params,\paramt_0 + \paramt\oa_m) ~ \omega\oa_m .
\end{align*}

{
To convert the above formula to the form \eqref{hwei},
it only remains to insert a trigonometric polynomial
interpolant for $\sigma(\params,\cdot)$ in terms of
the samples $\sigma(\params,\paramt_j)$.
Firstly, using the Fourier series ${ \sigma(\params,\paramt) \approx \sum_{l=1-\Nt/2}^{\Nt/2} \hat{\sigma}_l(\params) ~ e^{i l \paramt} }$,
the above becomes
\begin{align}
 I(\xx, \params) \approx &
 \sum_{l=1-\Nt/2}^{\Nt/2} \hat{\sigma}_l(\params) ~ \hat{\ModalIntegOp}_l(\xx, \params),
\label{hath}
\end{align}
} where
\begin{align}
  \hat{\ModalIntegOp}_l({\vct{x}}, \params) \coloneqq
  \sum\limits_{m=1}^{m_0} \GreensFn(\xx,\su(\params,\paramt_0 + \paramt\oa_m)) ~ e^{il (\paramt_0 + \paramt\oa_m)} ~ \SurfJacobian(\params,\paramt_0 + \paramt\oa_m) ~ \omega\oa_m,
  \qquad l = 1-\Nt/2, \dots, \Nt/2.
  \label{e:modal-greens-fn-coeff}
\end{align}
{
Secondly, inserting the trapezoid approximation
to the Euler--Fourier formula,
${\hat{\sigma}_l(\params) \approx \Nt^{-1} \sum_{j=1}^{\Nt} e^{-2\pi i l j/\Nt}\sigma(\params,\paramt_j)}$,
into \eqref{hath} and swapping the order of summation
recovers \eqref{hwei}
with the desired weights
}
\begin{align}
  \ModalIntegOp_j(\xx,s) \;=\;
  \frac{1}{\Nt} \sum\limits_{l=1-\Nt/2}^{\Nt/2} \hat{\ModalIntegOp}_l(\xx,s) \, e^{-2 \pi i j l / \Nt},
  \qquad j = 1, \dots, \Nt,
  \label{e:modal-greens-fn-wts}
\end{align}
an inverse discrete Fourier transform (DFT) of
{the vector of weights in} \eqref{e:modal-greens-fn-coeff}.

{
In practice, we exploit some further accelerations.
}
We only need to evaluate $\hat{\ModalIntegOp}_l$ for $l \geq 0$,
since $\hat{\ModalIntegOp}_{-l} = \hat{\ModalIntegOp}_{l}^*$.
The Fourier modes evaluated at the quadrature nodes ($e^{il\paramt\oa_m}$) are precomputed
and stored along with the quadrature rule to avoid expensive evaluation
of complex exponentials each time.
Since $\Nt$ is small ($\leq 100$), we build the operator matrix for the inverse DFT and apply it as a matrix-vector product (or matrix-matrix product when batched).

Once the coefficients $h_j(\xx, \params)$ have been constructed via \pr{e:modal-greens-fn-coeff,e:modal-greens-fn-wts}, the integral \eqref{Idef} can then be approximated to high order for any smooth density $\sigma$ using \pr{hwei}.
The needed values $\sigma(\params,\paramt_j)$ are interpolated in $\params$ from the known surface values
$\sigma\kk_{ij}$ in the manner explained in the next section.

\begin{figure}[htb!]
  \resizebox{0.85\textwidth}{!}{\begin{tikzpicture}
    \node[anchor=south west,inner sep=0] at (0,0) {\includegraphics[width=13.85cm]{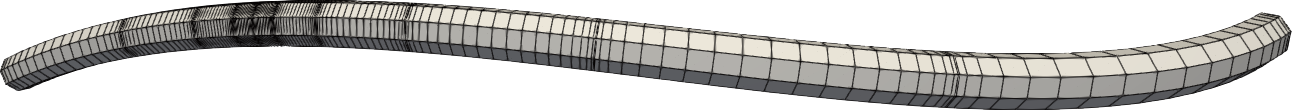}};

    \draw[color=red, ultra thick] (2.73,0.5) circle (1pt);
    \node at (2.95, 0.45) {\color{red} \large $\vct{x}$};

    \draw[ultra thick, ->] (6.7,0.1) to (5.5,0.3);
    \node [rotate=-4.5] at (9.5, -0.17) {\small adaptively refined GL panel quadrature};
  \end{tikzpicture}}
  \caption{\label{f:near-sing-quad-s}
  A slender body element $\SlenderElem_k$ with a target point $\vct{x}$ close to it.
  To evaluate the potential at $\vct{x}$, for the outer integral
  we use a panel based Gauss--Legendre quadrature
  rule in $\params$. Panels are refined
  adaptively near the target $\vct{x}$, giving a total of $m_a$ nodes.
  }
\end{figure}
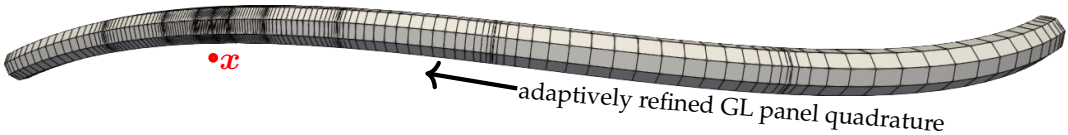

\subsubsection{Near-singular quadrature \label{sss:near-sing-quad}}
When evaluating the potential from a slender body element $\SlenderElem_{k}$
at a target point $\vct{x}$ that is near it (\ie ~$\vct{x}\in\Near_{\SlenderElem_k}$) but off-surface,
we use an adaptive panel based Gauss--Legendre quadrature rule
to compute the integral in $\params$ along the length of the element,
as shown in \pr{f:near-sing-quad-s}.
{
We use error estimates for Gauss--Legendre quadrature (\cite[Thm.~19.3]{Trefethen2019})
to determine the optimal quadrature order for integrands that are analytic within a Bernstein $\rho$-ellipse,
with a parameter $\rho=2.5$ chosen empirically to
minimize the total number of resulting nodes.
Then, we build the adaptive quadrature rule by subdividing the panels so that the target evaluation point $\vct{x}$
lies outside of each sub-panel near field ${\cal N}$ as
in \pr{sss:nystrom}.
}
%
The integral in $\paramt$ is computed as described in \pr{sss:modal-greens-fn}.
The potential at $\vct{x}$ due to $\SlenderElem_k$ can be approximated as,
\begin{align}
\int_{I_k} \int_{0}^{2\pi} \GreensFn(\xx,\su(\params,\paramt)) ~ \sigma(\params,\paramt) ~ \SurfJacobian(\params,\paramt) ~ d\paramt d\params
~\approx~
\sum_{m=1}^{m_a} w^{pan}_m \sum\limits_{j=1}^{\Ntk}
\ModalIntegOp_j(\vct{x}, \params^{pan}_m) \,
\sigma(\params^{pan}_m, \paramt_j)
,
\end{align}
where $\{s^{pan}_1, \dots, s^{pan}_{m_a}\}$ and $\{w^{pan}_1, \dots, w^{pan}_{m_a}\}$
are the complete sets of nodes and weights of the adaptive panel quadrature rule in
$\params$,
while $\ModalIntegOp_j$ are the components of \pr{e:modal-greens-fn-wts}.
We use barycentric Lagrange interpolation in $s$ to interpolate from
the given $\Nsk$ nodes $\params\kk_i$ to the panel quadrature nodes $\params_m^{pan}$.
Let $P$ be the resulting interpolation matrix, so that
$\sigma(\params^{pan}_m, \paramt_j) = \sum_{i=1}^{\Nsk} P_{mi} \sigma\kk_{ij}$
for each $j=1,\dots,\Ntk$.
Then the off-surface
Nystr\"{o}m corrections $\NystCorrec_{ij}\kk(\xx)$ in \pr{e:nystrom-correction}
are given by
\begin{align}
\NystCorrec_{ij}\kk(\xx)
~=~ 
\sum_{m=1}^{m_a}  P_{mi} ~ w^{pan}_m ~ \ModalIntegOp_j(\vct{x}, \params^{pan}_m)
~-~ \GreensFn(\xx, \su\kk_{ij}) ~ w\kk_{ij}
,
\qquad i=1,\dots,\Nsk, \quad j=1,\dots,\Ntk,
\end{align}
where $w\kk_{ij}$ are the surface element far-field quadrature weights that were used in \pr{e:nystrom-no-correction}.

\begin{figure}[htb!]
  \resizebox{0.85\textwidth}{!}{\begin{tikzpicture}
    \node[anchor=south west,inner sep=0] at (0,0) {\includegraphics[width=13.85cm]{figs/s-quad/adap-quad.png}};

    \draw[color=red, ultra thick] (2.7,0.9) circle (1pt);
    \node at (2.6, 0.5) {\color{red} \large $\vct{x}$};

    \draw[ultra thick, ->] (4.3,1.45) to (3.1,1.5);
    \node [rotate=-6] at (5.55, 1.3) {log singularity};

    \draw[ultra thick, ->] (10.5,0.85) to (12.1,0.9);
    \node [rotate=-4.5] at (9.6, 0.97) {$|\params-\params_0|^{-\alpha}$};

    \draw[ultra thick, ->] (3.3,-0.0) to (2.9,0.65);
    \node [rotate=0] at (3.5, -0.25) {special quadrature};
    \node [rotate=0] at (3.5, -0.70) {for $p(\params) \log |\params-\params_0| + q(\params)$};

    \draw[ultra thick, ->] (6.9,0.1) to (5.5,0.3);
    \node [rotate=-3] at (10, -0.2) {dyadically refined GL panel quadrature};
  \end{tikzpicture}}

  \vspace{2em}
  \resizebox{0.85\textwidth}{!}{\begin{tikzpicture}
    \node[anchor=south west,inner sep=0] at (0,0) {\includegraphics[width=13.85cm]{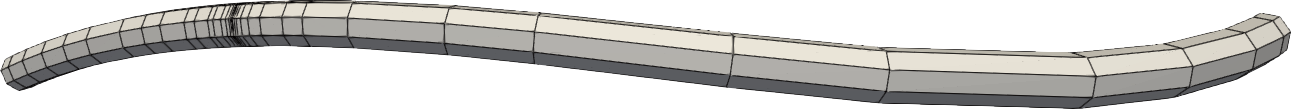}};
    \draw[color=red, ultra thick] (2.49,0.89) circle (1pt);
    \node at (2.5, 0.5) {\color{red} \large $\vct{x}$};

    \draw[ultra thick, ->] (2.5,1.45) to (1.3,1.3);
    \node [rotate=-4] at (6.6, 1.25) {special quadrature rule replaces several GL panels};
    \draw[ultra thick, ->] (10.9,0.95) to (12.5,1.05);
  \end{tikzpicture}}
  \caption{\label{f:singular-quad-s}
    Illustration of our special quadrature rule in $\params$ for
    an on-surface target point $\vct{x}=\su(\params_0,\paramt_0)$.
    The $\params$-integral is singular with the asymptotic forms shown.
  Top: This integral could be approximated using dyadically refined Gauss--Legendre panel quadrature,
  plus a special singular quadrature rule for the panels on either side of $\vct{x}$.
  Bottom: For efficiency we replace such a panel quadrature rule by a {\em single}
  generalized Chebyshev quadrature rule with far fewer ($m_0$) nodes.
  }
\end{figure}

\subsubsection{Singular quadrature \label{sss:sing-quad}}
When evaluating the potential from a slender body element $\SlenderElem_k$ at a target point $\vct{x}$ that is on-surface
{(\ie\, one of the surface discretization nodes of $\SlenderElem_k$)},
for the Laplace and Stokes kernels, the integrand in $\params$ has a logarithmic singularity at the target parameter $\params_0$.
{This follows from the well-known logarithmic singularity of the
modal Green's functions (see \pr{sss:modal-greens-fn}) with respect to
distance from the source circle;
e.g., see \cite[Sec.~5]{youngmart}
for the Laplace case
(the $l$th mode involving the
associated Legendre function ${\cal Q}_{l-1/2}$),
and \cite[Sec.~3.1]{veera09} for the mode $l=0$ for Stokes.
}
We could compute this integral using a panel based quadrature rule as shown in \pr{f:singular-quad-s} (top).
In such an approach, the panels are refined dyadically around the target until the smallest panel length is approximately $\radius(\params)$.
The Gauss--Legendre quadrature rule is used for panels not touching $\vct{x}$, while
a special singular quadrature rule is used for the panels on either side of $\vct{x}$.
We constructed the latter singular rule using the method described in \pr{s:cheb-quad} using integrands of the form
$p(\params) \log|\params-\params_0| + q(\params)$, where $p(\params)$ and $q(\params)$ are polynomials in $\params$.

However, using such
dyadically refined panel quadratures can become expensive when $\radius(\params)$ is small,
simply because the panels must refine down to the $\radius$ scale
before the logarithmic singularity becomes dominant (note that at larger scales
the singularity tends to the pointwise Green's function, thus is
$1/r$ or $1/r^2$).
This has been realized in the electromagnetic setting \cite{haslamwire}.
Instead, we precompute generalized Chebyshev quadrature rules that integrate
the {\em entire} length of $\SlenderElem_k$ at once.
For this we generate integrand functions using several ($\sim20$) straight cylindrical slender body elements
with different aspect ratios (the ratio of the length of the slender body element
to its radius $\radius$) {sampled uniformly} in the range $\beta$ to $2\beta$.
We did not try using elements with different curvature since the number of integrands was already very large; \eg\, for the Stokes single-layer kernel we have $6 \times 20 \Ns \Nt$ scalar integrand functions ($6\times$ is because the kernel is a symmetric $3 \times 3$ tensor).
We first use the dyadically refined panel quadratures described above to discretize these integrands,
then use the method in \pr{s:cheb-quad} to replace this by
a generalized Chebyshev quadrature with far fewer ($m_0$) quadrature nodes.
{Each quadrature rule constructed in this way works for all aspect ratios in the interval $(\beta, 2\beta)$.}
We precompute and store separate quadrature rules for
each of the $\Ns$ positions of the target within the panel
(\ie~each Gauss--Legendre node in \params,
since the on-surface targets are just the surface discretization nodes)
and for different slender element aspect ratios for $\beta \in \{2^{-7}, 2^{-6}, \dots, 2^{30}\}$.
For 10-digit accuracy and $\Nt=10$, typically $m_0$ is in the range 32 to 47.

To compute on-surface Nystr\"{o}m corrections $\NystCorrec_{ij}\kk(\xx)$ in \pr{e:nystrom-correction},
we select the appropriate quadrature rule for the aspect ratio of $\SlenderElem_k$.
Let $\{s^{sing}_1, \dots, s^{sing}_{m_0}\}$ and $\{w^{sing}_1, \dots, w^{sing}_{m_0}\}$
be the nodes and weights for the singular quadrature in $\params$.
As in the previous subsection,
we build a Lagrange interpolation matrix $P$
to interpolate from the $\params$ discretization nodes to the
special quadrature nodes such that 
$\sigma(\params^{sing}_m, \paramt_j) = \sum_{i=1}^{\Nsk} P_{mi} \sigma\kk_{ij}$.
The Nystr\"{o}m corrections are given by
\begin{align}
\NystCorrec_{ij}\kk(\xx)
~=~ 
\sum_{m=1}^{m_0}  P_{mi} ~ w^{sing}_m ~ \ModalIntegOp_j(\vct{x}, \params^{sing}_m)
~-~ \GreensFn(\xx,\su\kk_{ij}) ~ w\kk_{ij}
,
\qquad i=1,\dots,\Nsk, \quad j=1,\dots,\Ntk,
\label{e:sing-nyst-correc}
\end{align}
where $\ModalIntegOp_j$ are again given by \pr{e:modal-greens-fn-wts}
and $w\kk_{ij}$ are the far-field quadrature weights that were used in \pr{e:nystrom-no-correction}.

\subsection{Parallel algorithms \label{ss:parallel}}
We support distributed memory parallelism using MPI.
This requires partitioning the data across MPI processes, while also achieving good load balance across processes.
The data consists of a global array of slender body elements that is partitioned across processes.
We load balance for the quadrature setup stage by assigning a cost estimate to each element using the cost analysis for $\Tsetup$ in \pr{ss:cost},
and repartition this global array so that each local section has similar total cost.
For the layer potential evaluation stage, the far-field computation ({\ie} the N-body sum) is the dominant cost, and this has a different cost estimate than the quadrature setup.
Therefore, for the evaluation stage, all the source points (the discretization nodes $\su\kk_{ij}$) and the target points are repartitioned equally across processes; this is done within the PVFMM library.

The quadrature setup step is essentially a local operation for each slender body element.
However, identifying the set of target points in the near region $\Near_{\SlenderElem_k}$
of a slender body element $\SlenderElem_k$ requires communication since these points may be on different processes.
As described previously in \pr{sss:nystrom}, we approximate $\Near_{\SlenderElem_k}$
by a set of overlapping spheres centered at the surface discretization nodes $\su\kk_{ij}$.
To identify the target points that overlap with these spheres requires an octree-like data structure
to partition the space hierarchically and allow for searching in the neighborhood of a tree node.
However, instead of using a standard tree data structure, we use space-filling curves with Morton ordering (also called a hashed octree \cite{warren93}).
It provides most of the same functionalities (such as $\bigO{\log N}$ searching),
while having several performance advantages due to better memory access patterns.
The spheres and the target points are sorted in Morton order using a distributed sorting algorithm \cite{Sundar2013}.
Each tree node can be identified as a contiguous section of this Morton sorted array and these sections can be identified through binary searching in logarithmic time complexity.
The partitioning of this array across processes also gives a partitioning of the domain across processes.
The radius of a sphere determines its depth in the tree (according to $d=\lfloor \log_2 r \rfloor$ where $r$ is the radius),
and we search among the sphere's neighboring tree nodes at that depth to identify target points that overlap with the sphere.
For spheres that have neighboring nodes outside of the domain of the current process,
we send the spheres to the processes containing those neighboring nodes.
After identifying the target points in each slender body element's near region, these points are sent to the process where the element originated.
Duplicate target points must be removed since a point can overlap with more than one sphere of the same element.
For each target and element pair ($\vct{x}, \SlenderElem_k$) such that $\vct{x} \in \Near_{\SlenderElem_k}$,
the local Nystr\"{o}m corrections $\NystCorrec_{ij}\kk(\xx)$ can now be computed as described in \pr{sss:near-sing-quad,sss:sing-quad}.
We keep track of where each target point originated, since in the evaluation stage, the corrections from all near-elements of each target must be gathered and added to the final target potential.

\subsection{Computational Cost \label{ss:cost}}
The quadrature setup step requires computing the local Nystr\"om corrections for the singular and near-singular interactions.
The algorithm in \pr{ss:parallel} to identify the near interaction targets is dominated by the cost of the sorting algorithm
and requires $\bigO{\Nunknown / p \, \log \Nunknown}$ time,
where $\Nunknown$ is the total number of unknowns and $p$ is the number of processes.
The number of near interactions is geometry dependent; however, we will assume that the total number of singular and near-singular interactions is proportional to the number of surface discretization nodes.
We also assume that the number of generalized Chebyshev quadrature nodes in $\paramt$ and $\params$ scale linearly with the discretization order $\Nt$ and $\Ns$ respectively; this is justified by a heuristic geometric convergence
with respect to each of these parameters.
Then, \pr{e:modal-greens-fn-coeff,e:modal-greens-fn-wts}, requires $\bigO{\Nt^2}$ work to compute all $\ModalIntegOp_j$.
For each target $\vct{x}$, \pr{e:sing-nyst-correc} requires computing $\bigO{\Ns}$ modal Green's functions and $\bigO{\Ns}$ work (for each $i, j$ index) to apply the interpolation operator.
Therefore, building the singular corrections has a cost of $\bigO{\Ns \Nt^2 + \Nt \Ns^2}$ per target point.
If we assume a constant number (on average) of dyadically refined panels in $\params$ for near-singular interactions, then we arrive at the same asymptotic cost estimate for the near-singular interactions as well.
If the discretization orders $\Nt$ and $\Ns$ are assumed to not vary too much across all the elements, then the overall quadrature setup time is,
\begin{align}
  \Tsetup = \bigO{\frac{\Nunknown}{p} \log \Nunknown + \frac{\Nunknown}{p} (\Ns \Nt^2 + \Nt \Ns^2) } .
\end{align}

The quadrature evaluation step is dominated by N-body sum (first term in \pr{e:nystrom-no-correction}).
When the sum is computed directly, the evaluation time is $\Teval = \bigO{\Nunknown^2 / p}$.
This can be accelerated to $\Teval = \bigO{\Nunknown / p}$ (neglecting communication costs) using the PVFMM library.
However, due to extra overheads and lower parallel efficiency of FMM for small problem sizes,
direct evaluation may be faster for problem sizes smaller than $\sim 30K$ points per CPU core on up to $40$ cores.

\section{Boundary integral equation formulations\label{s:bie}}
Although our focus is (vector-valued) Stokes BVPs, it will be beneficial first to study the
simpler (scalar) Laplace BVP.
Thus in \pr{ss:dirichlet-bvp}, for Laplace and Stokes Dirichlet boundary value problems, we present new indirect
boundary integral equation (BIE) formulations which exhibit
small condition numbers uniformly as the cross-sectional radius $\radius\to 0$.
This allows us to solve these BIEs efficiently using iterative solvers such as GMRES.
Then, using the analogous formulation, in \pr{ss:mobility-formulation}, we develop a new BIE formulation for the Stokes mobility problem.
Unlike existing formulations, it remains well-conditioned for slender-body geometries.

\subsection{Dirichlet boundary value problems \label{ss:dirichlet-bvp}} 
We first consider the Laplace Dirichlet exterior BVP
\begin{align}
    \Delta u &= 0 & \text{in } \Real^3 \setminus \overline{\Omega} \label{e:laplace-dirichlet0} , \\
    u(\vct{x}) &\rightarrow 0 & \text{as } |\vct{x}| \rightarrow \infty \label{e:laplace-dirichlet1}, \\
    u &= u_0 & \text{on } \partial\Omega  \label{e:laplace-dirichlet2} ,
\end{align}
where $u_0 \in C(\partial\Omega)$ is given surface voltage data,
and the solution $u$ may be interpreted as an electrostatic potential.
The solution exists and is unique \cite[Ch.~6]{LIE}.
Recall the single-layer kernel $G(\xx,\yy) = (4\pi |\xx-\yy|)^{-1}$
and double-layer kernel $\partial G(\xx,\yy)/\partial \vct{n}_\yy$,
which induce the layer representations $\LaplaceSLR$ and $\LaplaceDLR$
respectively \cite{LIE}.
We will use a ``combined field'' representation,
\begin{align}
    u &= (\LaplaceDLR + \eta\LaplaceSLR)[\sigma] & \text{in } \Real^3 \setminus \overline\Omega, \label{e:laplace-bie-rep}
\end{align}
where $\eta>0$ is a fixed mixing parameter.
Then, $u$ satisfies \pr{e:laplace-dirichlet0,e:laplace-dirichlet1} by construction.
Taking the limit as $\vct{x} \rightarrow \partial\Omega$, invoking the jump relations,
and substituting in \pr{e:laplace-dirichlet2}, we get a second-kind boundary integral equation (BIE) in the unknown $\sigma$, namely
\begin{align}
    (I/2 + \LaplaceDL + \eta\LaplaceSL)\sigma &= u_0  & \text{on } \partial\Omega , \label{e:laplace-bie}
\end{align}
where $\LaplaceSL$ and $\LaplaceDL$ are the
weakly singular (hence compact) boundary integral operators
given by restricting the single- and double-layer representations
to $\partial\Omega$; note that $\LaplaceDL$ is taken in the principal value sense.
We can solve this BIE for $\sigma$, and then use \pr{e:laplace-bie-rep} to evaluate $u$ in $\Real^3 \setminus \overline{\Omega}$.

A word is needed about the representation choice
\pr{e:laplace-bie-rep}, which is not commonly used for Laplace.
A plain single-layer representation $u=\LaplaceSLR[\sigma]$
would lead to a first-kind integral equation, whose discretization
is thus ill-conditioned; yet this is sometimes used in 3D low-order settings.
A plain double-layer (setting $\eta=0$ above)
cannot represent the $\bigO{1/r}$ term associated with non-zero net charge,
and furthermore the resulting integral equation,
while second-kind, is not uniquely solvable:
its operator $I/2 + D$ has a null space of dimension equal to the
number of bodies comprising $\Omega$.
A common remedy to recover unique
solvability is to add a rank-1 operator (per body) to $I/2 + D$
\cite[\S38]{mikhlinbook} \cite[Thm.~6.23]{LIE}.
%
Using the representation \pr{e:laplace-bie-rep} is an alternative remedy;
its proof of unique solvability hinges on uniqueness for the interior
Robin BVP when $\eta>0$
(see \cite[\S4.2]{manasthesis} for the trickier $\Real^2$ version).
It is analogous to the popular combined-field Helmholtz
\cite[\S3.9.4]{sauterschwab} or
``completed'' Stokes representations \cite{hsiao85,hebeker}.
In all known prior work $\eta$ is chosen as an $\bigO{1}$ constant.

Even though \pr{e:laplace-bie} is second-kind, with fixed $\eta\approx1$ we
observed that it becomes increasingly ill-conditioned for slender-body geometries as the fiber radius $\radius \rightarrow 0$.
When solving the discretized BIE using GMRES, for extremely slender geometries, the result is
failure to converge even after hundreds of iterations.
However, we found that this can be remedied by a specific
$\radius$-dependent scaling of $\eta$,
empirically resulting in uniformly small condition number as $\radius \to 0$.

\begin{figure}[htb!]
\resizebox{0.7\textwidth}{!}{\begin{tikzpicture}

\draw[very thick,red, ->] (0, 0,0) -- (4, 0,0);
\node at (2.5,0.22,0) {\params};

\draw[very thick,black, ->] (-1, 0.2,0) -- (-1, 1,0);
\draw[very thick,black, ->] (-1, -0.22,0) -- (-1, -1,0);
\node at (-1,0,0) {2\radius};

\node at (0,0,0) [cylinder, shape aspect=4, draw,minimum height=8cm,minimum width=2cm, cylinder uses custom fill, cylinder end fill=blue!40,cylinder body fill=blue!100, opacity=0.1] {};
\node at (0,0,0) [cylinder, shape aspect=4, draw,minimum height=8cm,minimum width=2cm, ultra thick] {};

\draw[thick, dashed] (4,-1,0) -- (6,-1,0);
\draw[thick, dashed] (4, 1,0) -- (6, 1,0);
\draw[thick, dashed] (-6,-1,0) -- (-3.25,-1,0);
\draw[thick, dashed] (-6, 1,0) -- (-3.25, 1,0);

\draw[ultra thick, ->, red] (-2,-1.15,0)  to [out=240,in=120, looseness=2] (-2,1.15,0);
\node at (-2.88,0.2,0) {\Large \paramt};

\draw[very thick,black, ->|] (-0.5, -1.6,0) -- (-3.25, -1.6,0);
\draw[very thick,black, ->|] ( 1.5, -1.6,0) -- ( 4.00, -1.6,0);
\node at (0.5,-1.65,0) {$L_{\sc{period}} = 2 \pi$};

\end{tikzpicture}}
\caption{\label{f:inf-cylinder}
  The infinite cylinder of radius $\radius$ with a periodic boundary condition with respect to length.
The surface of one period is parameterized by $(\params,\paramt)\in[0,2\pi)^2$.
For different values of $\radius$, in \pr{p:cylinder-eigenvalues} we plot the eigenvalues of the Laplace single-layer, double-layer, and the proposed new admixture boundary integral operators.}
\end{figure}
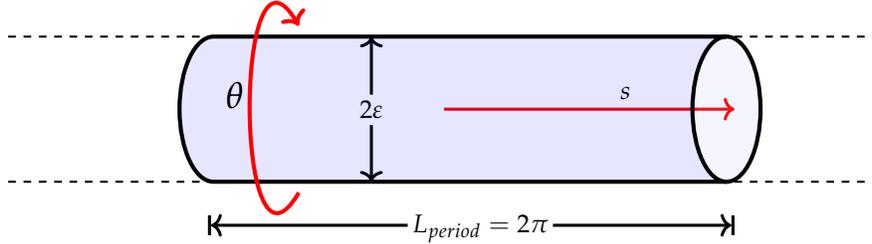

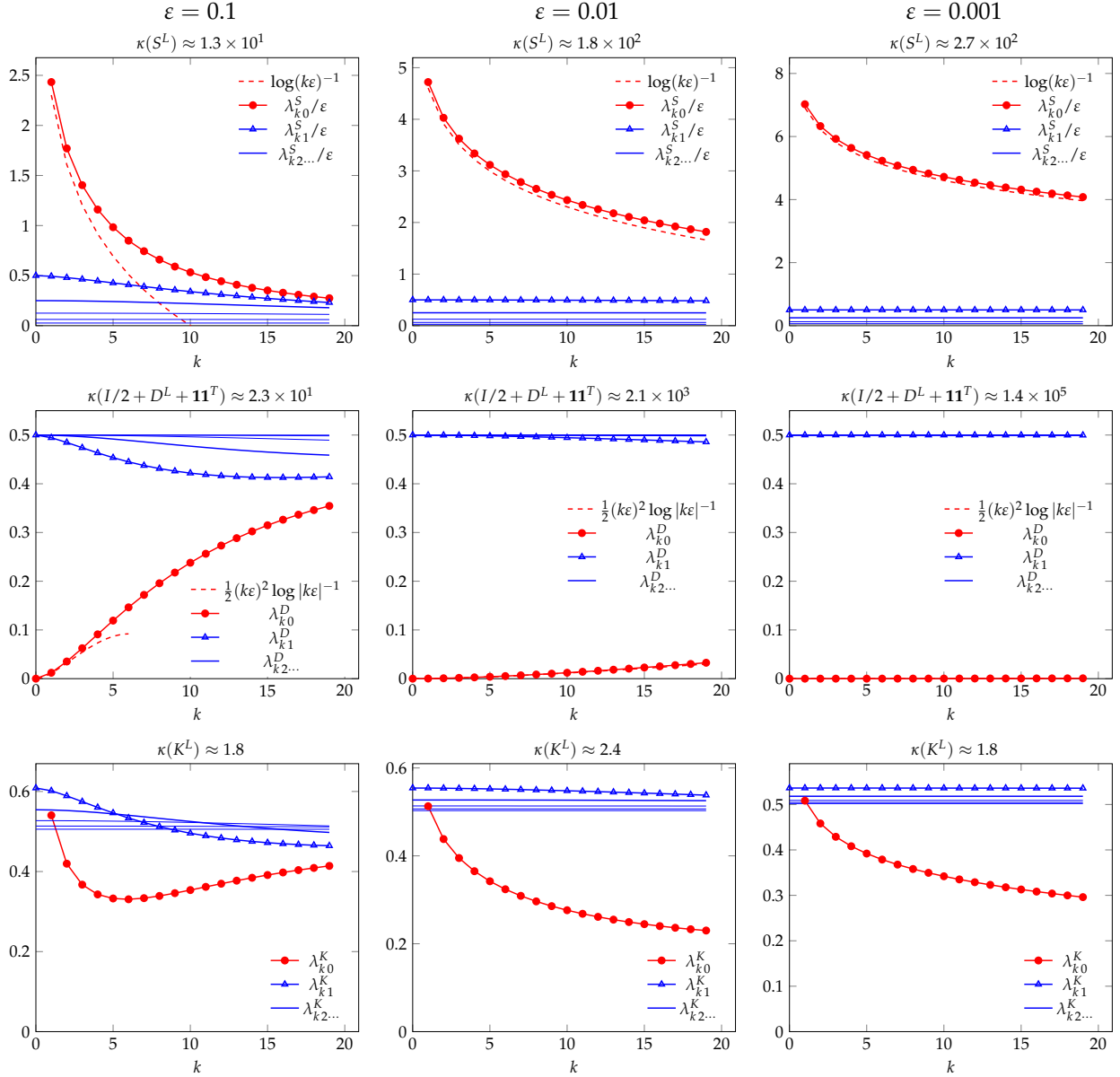
\begin{figure}[htb!]
\resizebox{0.99\textwidth}{!}{\begin{tikzpicture}
  \begin{scope}[shift={(-8,0)}]
    \node at (3.5,6.7) {\Large $\radius = 0.1$};
    \node at (3.5,6) {$\kappa(\LaplaceSL) \approx 1.3\times10^1$};
    \begin{axis}[xmin=0, ymin=0, xlabel={$k$},legend style={draw=none,at={(0.97,0.97)},anchor=north east}]
      \addplot [thick,dashed,color=red] table [x={k},y={logk}] {data/spectrum-M20-eps0.1}; \addlegendentry{$\log(k \radius)^{-1}$};
      \addplot [thick,color=red,mark=*] table [x={k},y={S00}] {data/spectrum-M20-eps0.1}; \addlegendentry{$\lambda^{S}_{k\,0} / \radius$}

      \addplot [thick,color=blue,mark=triangle] table [x={k},y={S01}] {data/spectrum-M20-eps0.1}; \addlegendentry{$\lambda^{S}_{k\,1} / \radius$}
      \addplot [thick,color=blue] table [x={k},y={S02}] {data/spectrum-M20-eps0.1}; \addlegendentry{$\lambda^{S}_{k\,2\cdots} / \radius$}
      \addplot [color=blue] table [x={k},y={S04}] {data/spectrum-M20-eps0.1};
      \addplot [color=blue] table [x={k},y={S08}] {data/spectrum-M20-eps0.1};
      \addplot [color=blue] table [x={k},y={S19}] {data/spectrum-M20-eps0.1};
    \end{axis}
  \end{scope}
  
  \begin{scope}[shift={(-8,-7.5)}]
    \node at (3.5,6) {$\kappa(I/2 + \LaplaceDL + {\bf 11}^{T}) \approx 2.3\times10^1$};
    \begin{axis}[xmin=0, ymin=0, xlabel={$k$},legend style={draw=none,at={(0.97,0.2)},anchor=east}]
      \addplot [thick,dashed,color=red,domain=0.1:6,samples=100] (\x,{0.5*0.1*0.1*\x*\x*(log2(1/(0.1*\x))/1.4427+0.149*0)}); \addlegendentry{$\frac{1}{2}(k \radius)^2 \log |k \radius|^{-1}$};
      \addplot [thick,color=red,mark=*] table [x={k},y={D00}] {data/spectrum-M20-eps0.1}; \addlegendentry{$\lambda^{D}_{k\,0}$};

      \addplot [thick,color=blue,mark=triangle] table [x={k},y={D01}] {data/spectrum-M20-eps0.1}; \addlegendentry{$\lambda^{D}_{k\,1}$};
      \addplot [thick,color=blue] table [x={k},y={D02}] {data/spectrum-M20-eps0.1}; \addlegendentry{$\lambda^{D}_{k\,2\cdots}$};
      \addplot [color=blue] table [x={k},y={D04}] {data/spectrum-M20-eps0.1};
      \addplot [color=blue] table [x={k},y={D08}] {data/spectrum-M20-eps0.1};
      \addplot [color=blue] table [x={k},y={D19}] {data/spectrum-M20-eps0.1};
    \end{axis}
  \end{scope}

  \begin{scope}[shift={(-8,-15.0)}]
    \node at (3.5,6) {$\kappa(\LaplaceCF) \approx 1.8$};
    \begin{axis}[xmin=0, ymin=0, xlabel={$k$},legend style={draw=none,at={(0.97,0.03)},anchor=south east}]
      \addplot [thick,color=red,mark=*] table [x={k},y={SD00}] {data/spectrum-M20-eps0.1}; \addlegendentry{$\lambda^K_{k\,0}$};

      \addplot [thick,color=blue,mark=triangle] table [x={k},y={SD01}] {data/spectrum-M20-eps0.1}; \addlegendentry{$\lambda^K_{k\,1}$};
      \addplot [thick,color=blue] table [x={k},y={SD02}] {data/spectrum-M20-eps0.1}; \addlegendentry{$\lambda^K_{k\,2\cdots}$};
      \addplot [color=blue] table [x={k},y={SD04}] {data/spectrum-M20-eps0.1};
      \addplot [color=blue] table [x={k},y={SD08}] {data/spectrum-M20-eps0.1};
      \addplot [color=blue] table [x={k},y={SD19}] {data/spectrum-M20-eps0.1};
    \end{axis}
  \end{scope}

  \begin{scope}[shift={(0,0)}]
    \node at (3.5,6.7) {\Large $\radius = 0.01$};
    \node at (3.5,6) {$\kappa(\LaplaceSL) \approx 1.8\times10^2$};
    \begin{axis}[xmin=0, ymin=0, xlabel={$k$},legend style={draw=none,at={(0.97,0.97)},anchor=north east}]
      \addplot [thick,dashed,color=red] table [x={k},y={logk}] {data/spectrum-M20-eps0.01}; \addlegendentry{$\log(k \radius)^{-1}$};
      \addplot [thick,color=red,mark=*] table [x={k},y={S00}] {data/spectrum-M20-eps0.01}; \addlegendentry{$\lambda^{S}_{k\,0} / \radius$}

      \addplot [thick,color=blue,mark=triangle] table [x={k},y={S01}] {data/spectrum-M20-eps0.01}; \addlegendentry{$\lambda^{S}_{k\,1} / \radius$}
      \addplot [thick,color=blue] table [x={k},y={S02}] {data/spectrum-M20-eps0.01}; \addlegendentry{$\lambda^{S}_{k\,2\cdots} / \radius$}
      \addplot [color=blue] table [x={k},y={S04}] {data/spectrum-M20-eps0.01};
      \addplot [color=blue] table [x={k},y={S08}] {data/spectrum-M20-eps0.01};
      \addplot [color=blue] table [x={k},y={S19}] {data/spectrum-M20-eps0.01};
    \end{axis}
  \end{scope}
  
  \begin{scope}[shift={(0,-7.5)}]
    \node at (3.5,6) {$\kappa(I/2 + \LaplaceDL + {\bf 11}^{T}) \approx 2.1\times10^3$};
    \begin{axis}[xmin=0, ymin=0, xlabel={$k$},legend style={draw=none,at={(0.97,0.5)},anchor=east}]
      \addplot [thick,dashed,color=red,domain=0.1:19,samples=100] (\x,{0.5*0.01*0.01*\x*\x*(log2(1/(0.01*\x))/1.4427+0.149*0)}); \addlegendentry{$\frac{1}{2}(k \radius)^2 \log |k \radius|^{-1}$};
      \addplot [thick,color=red,mark=*] table [x={k},y={D00}] {data/spectrum-M20-eps0.01}; \addlegendentry{$\lambda^{D}_{k\,0}$};

      \addplot [thick,color=blue,mark=triangle] table [x={k},y={D01}] {data/spectrum-M20-eps0.01}; \addlegendentry{$\lambda^{D}_{k\,1}$};
      \addplot [thick,color=blue] table [x={k},y={D02}] {data/spectrum-M20-eps0.01}; \addlegendentry{$\lambda^{D}_{k\,2\cdots}$};
      \addplot [color=blue] table [x={k},y={D04}] {data/spectrum-M20-eps0.01};
      \addplot [color=blue] table [x={k},y={D08}] {data/spectrum-M20-eps0.01};
      \addplot [color=blue] table [x={k},y={D19}] {data/spectrum-M20-eps0.01};
    \end{axis}
  \end{scope}

  \begin{scope}[shift={(0,-15.0)}]
    \node at (3.5,6) {$\kappa(\LaplaceCF) \approx 2.4$};
    \begin{axis}[xmin=0, ymin=0, xlabel={$k$},legend style={draw=none,at={(0.97,0.03)},anchor=south east}]
      \addplot [thick,color=red,mark=*] table [x={k},y={SD00}] {data/spectrum-M20-eps0.01}; \addlegendentry{$\lambda^K_{k\,0}$};

      \addplot [thick,color=blue,mark=triangle] table [x={k},y={SD01}] {data/spectrum-M20-eps0.01}; \addlegendentry{$\lambda^K_{k\,1}$};
      \addplot [thick,color=blue] table [x={k},y={SD02}] {data/spectrum-M20-eps0.01}; \addlegendentry{$\lambda^K_{k\,2\cdots}$};
      \addplot [color=blue] table [x={k},y={SD04}] {data/spectrum-M20-eps0.01};
      \addplot [color=blue] table [x={k},y={SD08}] {data/spectrum-M20-eps0.01};
      \addplot [color=blue] table [x={k},y={SD19}] {data/spectrum-M20-eps0.01};
    \end{axis}
  \end{scope}

  \begin{scope}[shift={(8,0)}]
    \node at (3.5,6.7) {\Large $\radius = 0.001$};
    \node at (3.5,6) {$\kappa(\LaplaceSL) \approx 2.7\times10^2$};
    \begin{axis}[xmin=0, ymin=0, ymax=8.5, xlabel={$k$},legend style={draw=none,at={(0.97,0.97)},anchor=north east}]
      \addplot [thick,dashed,color=red] table [x={k},y={logk}] {data/spectrum-M20-eps0.001}; \addlegendentry{$\log(k \radius)^{-1}$};
      \addplot [thick,color=red,mark=*] table [x={k},y={S00}] {data/spectrum-M20-eps0.001}; \addlegendentry{$\lambda^{S}_{k\,0} / \radius$};

      \addplot [thick,color=blue,mark=triangle] table [x={k},y={S01}] {data/spectrum-M20-eps0.001}; \addlegendentry{$\lambda^{S}_{k\,1} / \radius$}
      \addplot [thick,color=blue] table [x={k},y={S02}] {data/spectrum-M20-eps0.001}; \addlegendentry{$\lambda^{S}_{k\,2\cdots} / \radius$}
      \addplot [color=blue] table [x={k},y={S04}] {data/spectrum-M20-eps0.001};
      \addplot [color=blue] table [x={k},y={S08}] {data/spectrum-M20-eps0.001};
      \addplot [color=gray] table [x={k},y={S19}] {data/spectrum-M20-eps0.001};
    \end{axis}
  \end{scope}
  
  \begin{scope}[shift={(8,-7.5)}]
    \node at (3.5,6) {$\kappa(I/2 + \LaplaceDL + {\bf 11}^{T}) \approx 1.4\times10^5$};
    \begin{axis}[xmin=0, ymin=0, xlabel={$k$},legend style={draw=none,at={(0.97,0.5)},anchor=east}]
      \addplot [thick,dashed,color=red,domain=0.1:19,samples=100] (\x,{0.5*0.001*0.001*\x*\x*(log2(1/(0.001*\x))/1.4427+0.149*0)}); \addlegendentry{$\frac{1}{2}(k \radius)^2 \log |k \radius|^{-1}$};
      \addplot [thick,color=red,mark=*] table [x={k},y={D00}] {data/spectrum-M20-eps0.001}; \addlegendentry{$\lambda^{D}_{k\,0}$};

      \addplot [thick,color=blue,mark=triangle] table [x={k},y={D01}] {data/spectrum-M20-eps0.001}; \addlegendentry{$\lambda^{D}_{k\,1}$};
      \addplot [thick,color=blue] table [x={k},y={D02}] {data/spectrum-M20-eps0.001}; \addlegendentry{$\lambda^{D}_{k\,2\cdots}$};
      \addplot [color=blue] table [x={k},y={D04}] {data/spectrum-M20-eps0.001};
      \addplot [color=blue] table [x={k},y={D08}] {data/spectrum-M20-eps0.001};
      \addplot [color=blue] table [x={k},y={D19}] {data/spectrum-M20-eps0.001};
    \end{axis}
  \end{scope}

  \begin{scope}[shift={(8,-15.0)}]
    \node at (3.5,6) {$\kappa(\LaplaceCF) \approx 1.8$};
    \begin{axis}[xmin=0, ymin=0, xlabel={$k$},legend style={draw=none,at={(0.97,0.03)},anchor=south east}]
      \addplot [thick,color=red,mark=*] table [x={k},y={SD00}] {data/spectrum-M20-eps0.001}; \addlegendentry{$\lambda^K_{k\,0}$};

      \addplot [thick,color=blue,mark=triangle] table [x={k},y={SD01}] {data/spectrum-M20-eps0.001}; \addlegendentry{$\lambda^K_{k\,1}$};
      \addplot [thick,color=blue] table [x={k},y={SD02}] {data/spectrum-M20-eps0.001}; \addlegendentry{$\lambda^K_{k\,2\cdots}$};
      \addplot [color=blue] table [x={k},y={SD04}] {data/spectrum-M20-eps0.001};
      \addplot [color=blue] table [x={k},y={SD08}] {data/spectrum-M20-eps0.001};
      \addplot [color=blue] table [x={k},y={SD19}] {data/spectrum-M20-eps0.001};
    \end{axis}
  \end{scope}
\end{tikzpicture}}
\caption{\label{p:cylinder-eigenvalues} Plot of the eigenvalues of the Laplace single-layer operator ($\LaplaceSL$) on the first row, the exterior double-layer operator $I/2 + \LaplaceDL$ on the second row,
and the new admixture operator $K^L = I/2 + \LaplaceDL + \LaplaceSL / (2 \radius \log \radius^{-1} )$ on the third row
for infinite periodic cylinder (\pr{f:inf-cylinder}) with period length $L_{period} = 2\pi$ and different radii $\radius=0.1$, $0.01$, and $0.001$.
The eigenvectors are of the form shown in \pr{e:cylinder-eigenvectors}, with index $k$ (on the horizontal axis) corresponding to the Fourier mode in $\params$.
The condition numbers $\kappa$ are for discretizations with 20 Fourier modes in both $\params$ and $\paramt$.
For the new admixture formulation, the condition number of the discretized operator $K^L$ remains uniformly small as $\radius\to0$.
}
\end{figure}

To understand these issues, we
study the spectrum of the boundary integral operators
in the subspace of periodic functions
on the infinite cylinder (\pr{f:inf-cylinder}).
The cylinder has radius $\radius$, the periodicity is $L_{period} = 2\pi$,
and the surface is parameterized by $\params$ along the length and $\paramt$ in angle.
This is a simple model for a closed filament ignoring center-line curvature,
also chosen for a recent slender-body BVP analysis \cite{mori21inv}.
By symmetry argument, the eigenfunctions for the Laplace single- and double-layer operators must be all complex exponentials of the form
\begin{equation}
  f_{kj}(\params,\theta) = e^{i(k\params + j\theta)},
  \qquad k,j \in \mathbb{Z}.
  \label{e:cylinder-eigenvectors}
\end{equation}
Let $\lambda^{S}_{kj}$ and $\lambda^{D}_{kj}$ be the corresponding eigenvalues of the on-surface single-layer operator $\LaplaceSL$ and the on-surface exterior double-layer operator $I/2 + \LaplaceDL$ respectively, \ie,
\begin{align*}
    \LaplaceSL f_{kj} &= \lambda^{S}_{kj} ~ f_{kj} , \\
    \bigl( I/2 + \LaplaceDL \bigr)f_{kj} &= \lambda^{D}_{kj} ~ f_{kj}.
\end{align*}
Pending an analytic study,
we computed these eigenvalues numerically by applying the quadratures described in \pr{s:numerical-algo} to a single period,
and using an approximately periodized Green's function
(sufficient accuracy of a few digits was achieved here
by naive summation of around $10^4$ periodic image sources).
The results are plotted in \pr{p:cylinder-eigenvalues}.
The eigenvalue $\lambda^{S}_{0\,0}$ of the single-layer operator is unbounded
(due to the log divergence of a line sum in 3D) and therefore not plotted.
For computing the condition number of $\LaplaceSL$, we assume that the operator acts on the space of mean-zero functions and therefore $\lambda^{S}_{0\,0}$ can be neglected.
The eigenvalue $\lambda^{D}_{0\,0}$ of the exterior double-layer operator $I/2 + \LaplaceDL$ is always zero;
therefore, we report the condition number of the latter after removing this null-space by adding a constant-kernel rank-1 update to the operator.
Since both the single- and double-layer operators have the same eigenfunctions $f_{kj}$, so does any linear combination of these operators;
the corresponding eigenvalues are the same linear combination
of $\lambda^{S}_{kj}$ and $\lambda^{D}_{kj}$.

For each
$k\in\mathbb{Z}$, $k\neq 0$, but staying in the limit
$|k \radius| \ll 1$ ($s$-wavelength much larger than the radius), we observed the eigenvalues $\lambda^{S}_{kj}$ and $\lambda^{D}_{kj}$ to be approximately given by
\be
    \lambda^{S}_{kj} \approx
    \begin{cases*}
        \radius \log |k \radius|^{-1}, \quad & if $j = 0$ \\
        \radius/2|j|,                       & otherwise,
    \end{cases*}
   \qquad
    \lambda^{D}_{kj} \approx
    \begin{cases*} 
        \frac{1}{2}(k \radius)^2 \log |k \radius|^{-1}, & if $j = 0$ \\
        0.5,                                           & otherwise.
    \end{cases*}
    \label{cyleig}
\ee
Firstly, note that $\lambda^{S}_{kj} \rightarrow 0$ as $|j| \rightarrow \infty$,
with the same form as the 2D single-layer operator on a radius-$\radius$
circle.
The combined-field operator of
\pr{e:laplace-bie} has the eigenvalues $\lambda^{D}_{kj} + \eta\lambda^{S}_{kj}$.
Since $\lambda^{D}_{kj} \approx 0.5$ for $j \neq 0$, these eigenvalues remain
safely bounded away from zero.
However, the eigenvalues corresponding to $j = 0$ are approximately $\bigO{\radius \log |k \radius|^{-1}}$ and therefore, for fixed $\eta$,
the condition number blows up as $\bigO{1/(\radius \log \radius^{-1})}$ as $\radius \rightarrow 0$.

\paragraph{Laplace slender combined field integral equation formulation}
We propose a new slender-body admixture where the single-layer operator is scaled by $\eta=1/(2 \radius \log \radius^{-1} )$.
This choice fixes the problematic $j=0$ eigenvalues so that
$\lambda^{D}_{k\, 0} + \eta\lambda^S_{k\, 0} \approx 0.5$ for all $k\neq 0$
in \pr{cyleig}, yet allows decay towards $0.5$ for each other $j\neq 0$.
In the exterior of $\Omega$ (which may comprise one or many bodies of similar slenderness),
the solution $u$ to the \pr{e:laplace-dirichlet0,e:laplace-dirichlet1,e:laplace-dirichlet2} is thus represented as
\begin{align}
  u &= \LaplaceDLR[\sigma] +
  \LaplaceSLR[\sigma] \, / \, (2 \radius \log \radius^{-1})
  & \text{in } \Real^3 \setminus \overline\Omega. \label{e:laplace-bie-rep-new}
\end{align}
Taking the exterior limit to $\partial\Omega$ and
applying the boundary conditions we get the BIE
\begin{align}
    \LaplaceCF \sigma &= u_0  & \text{on } \partial\Omega , \label{e:laplace-bie-new}
\end{align}
where $\LaplaceCF := I/2 + \LaplaceDL + \LaplaceSL / (2 \radius \log \radius^{-1} )$.
Indeed, for the cylinder in \pr{p:cylinder-eigenvalues} (last row), we show that the condition number of $\LaplaceCF$ remains small as $\radius$ shrinks.
This BIE can therefore be solved efficiently using GMRES.

\begin{figure}[htb!]
\resizebox{0.99\textwidth}{!}{\begin{tikzpicture}
  \begin{scope}[shift={(-6.9,0)}]
    \begin{axis}[width=7cm, height=8cm, xmin=0.125, xmax=131072, ymax=1e8, xmode=log,ymode=log,xlabel={$\eta$}, ylabel={$\kappa$}, ylabel near ticks, ylabel shift={-2pt}, legend style={draw=none,at={(0.4,1)},anchor=north west}]
      \addplot [thick,color=c3,mark=triangle] table [x={sl-scaling},y={cond1e-2}] {data/cond-torus-stokes-fixed-eps}; \addlegendentry{$\radius$=1e-2};
      \addplot [thick,color=c2,mark=square  ] table [x={sl-scaling},y={cond1e-4}] {data/cond-torus-stokes-fixed-eps}; \addlegendentry{$\radius$=1e-4};
      \addplot [thick,color=c1,mark=o       ] table [x={sl-scaling},y={cond1e-6}] {data/cond-torus-stokes-fixed-eps}; \addlegendentry{$\radius$=1e-6};

      \draw[thick, dashed, dash pattern=on 2pt off 2pt, color=c3] (axis cs:10.8574,1) -- (axis cs:10.8574,20);
      \draw[thick, dashed, dash pattern=on 2pt off 2pt, color=c2] (axis cs:542.868,1) -- (axis cs:542.868,38);
      \draw[thick, dashed, dash pattern=on 2pt off 2pt, color=c1] (axis cs:36192.2,1) -- (axis cs:36192.2,58);
    \end{axis}
  \end{scope}

  \begin{scope}[shift={(0,0)}]
    \begin{axis}[width=7cm, height=8cm, xmin=0.125, xmax=131072, ymax=22, xmode=log, xlabel={$\eta$}, ylabel={$\gmresiter$}, ylabel near ticks, ylabel shift={-3pt}, legend style={draw=none,at={(0.4,1)},anchor=north west}]
      \addplot [thick,color=c3,mark=triangle] table [x={sl-scaling},y={iter1e-2}] {data/cond-torus-stokes-fixed-eps}; \addlegendentry{$\radius$=1e-2};
      \addplot [thick,color=c2,mark=square  ] table [x={sl-scaling},y={iter1e-4}] {data/cond-torus-stokes-fixed-eps}; \addlegendentry{$\radius$=1e-4};
      \addplot [thick,color=c1,mark=o       ] table [x={sl-scaling},y={iter1e-6}] {data/cond-torus-stokes-fixed-eps}; \addlegendentry{$\radius$=1e-6};

      \draw[thick, dashed, dash pattern=on 2pt off 2pt, color=c3] (axis cs:10.8574,0) -- (axis cs:10.8574,7.5);
      \draw[thick, dashed, dash pattern=on 2pt off 2pt, color=c2] (axis cs:542.868,0) -- (axis cs:542.868,6.0);
      \draw[thick, dashed, dash pattern=on 2pt off 2pt, color=c1] (axis cs:36192.2,0) -- (axis cs:36192.2,5.0);
    \end{axis}
  \end{scope}
 
  \begin{scope}[shift={(7.7,3.3)}]
  \node at (0,0) {\includegraphics[angle=90,origin=c,width=3.28cm]{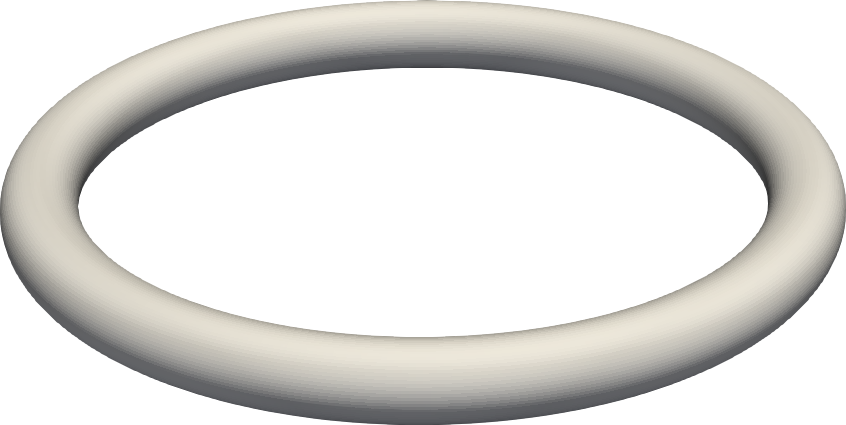}};
  \draw [red, thick, -{Latex[length=4.5pt,width=6.5pt]}] (-0.1,-2.36) -- (-0.1,-2.66);
  \draw [red, thick, -{Latex[length=4.5pt,width=6.5pt]}] (-0.1,-3.57) -- (-0.1,-3.27);
  \draw [red, thick] (-0.2,-2.66) -- (0.0,-2.66);
  \draw [red, thick] (-0.2,-3.27) -- (0.0,-3.27);
  \node at (-0.15, -2.94) {\color{red} \large $2\radius$};
  \end{scope}

  \begin{scope}[shift={(-6.9,-7.7)}]
    \begin{axis}[width=7cm, height=8cm, xmin=0.125, xmax=131072, ymax=1e8, xmode=log,ymode=log,xlabel={$\eta$}, ylabel={$\kappa$}, ylabel near ticks, ylabel shift={-2pt}, legend style={draw=none,at={(0.4,1)},anchor=north west}]
      \addplot [thick,color=c3,mark=triangle] table [x={sl-scaling},y={cond1e-2}] {data/cond-torus-stokes-var-eps}; \addlegendentry{$\radius$=1e-2};
      \addplot [thick,color=c2,mark=square  ] table [x={sl-scaling},y={cond1e-4}] {data/cond-torus-stokes-var-eps}; \addlegendentry{$\radius$=1e-4};
      \addplot [thick,color=c1,mark=o       ] table [x={sl-scaling},y={cond1e-6}] {data/cond-torus-stokes-var-eps}; \addlegendentry{$\radius$=1e-6};

      \draw[thick, dashed, dash pattern=on 2pt off 2pt, color=c3] (axis cs:10.8574,1) -- (axis cs:10.8574,26.0);
      \draw[thick, dashed, dash pattern=on 2pt off 2pt, color=c2] (axis cs:542.868,1) -- (axis cs:542.868,48.5);
      \draw[thick, dashed, dash pattern=on 2pt off 2pt, color=c1] (axis cs:36192.2,1) -- (axis cs:36192.2,74.9);
    \end{axis}
  \end{scope}

  \begin{scope}[shift={(0,-7.7)}]
    \begin{axis}[width=7cm, height=8cm, xmin=0.125, xmax=131072, ymax=200, xmode=log, xlabel={$\eta$}, ylabel={$\gmresiter$}, ylabel near ticks, ylabel shift={-5pt}, legend style={draw=none,at={(0.4,1)},anchor=north west}]
      \addplot [thick,color=c3,mark=triangle] table [x={sl-scaling},y={iter1e-2}] {data/cond-torus-stokes-var-eps}; \addlegendentry{$\radius$=1e-2};
      \addplot [thick,color=c2,mark=square  ] table [x={sl-scaling},y={iter1e-4}] {data/cond-torus-stokes-var-eps}; \addlegendentry{$\radius$=1e-4};
      \addplot [thick,color=c1,mark=o       ] table [x={sl-scaling},y={iter1e-6}] {data/cond-torus-stokes-var-eps}; \addlegendentry{$\radius$=1e-6};

      \draw[thick, dashed, dash pattern=on 2pt off 2pt, color=c3] (axis cs:10.8574,0) -- (axis cs:10.8574,29.4);
      \draw[thick, dashed, dash pattern=on 2pt off 2pt, color=c2] (axis cs:542.868,0) -- (axis cs:542.868,33.0);
      \draw[thick, dashed, dash pattern=on 2pt off 2pt, color=c1] (axis cs:36192.2,0) -- (axis cs:36192.2,28.0);
    \end{axis}
  \end{scope}
  
  \begin{scope}[shift={(7.7,-4.4)}]
  \node at (0,0) {\includegraphics[angle=90,origin=c,width=3.28cm]{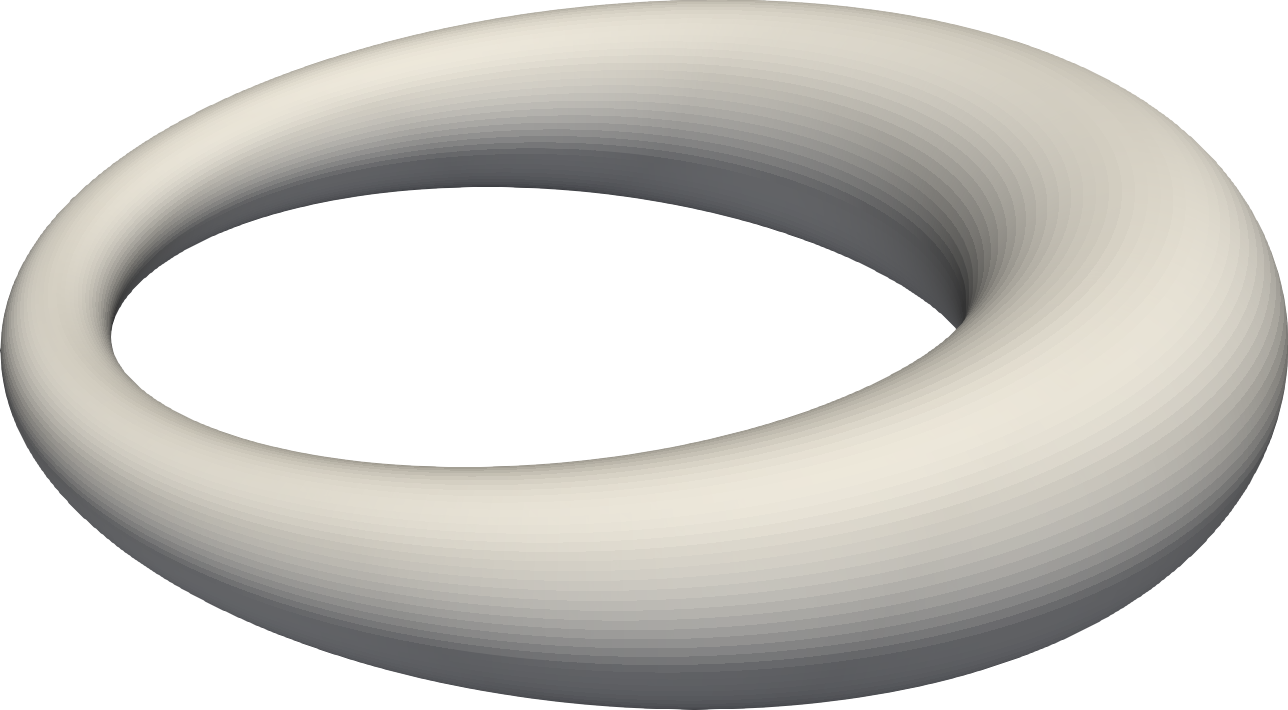}};
  \draw [red, thick, -{Latex[length=5.5pt,width=7.5pt]}] (-0.15,2.5) -- (-0.15,3.0);
  \draw [red, thick, -{Latex[length=5.5pt,width=7.5pt]}] (-0.15,2.0) -- (-0.15,1.5);
  \draw [red, thick] (-0.3,3.0) -- (0.0,3.0);
  \draw [red, thick] (-0.3,1.5) -- (0.0,1.5);
  \node at (-0.15, 2.25) {\color{red} \large $6\radius$};
  
  \draw [red, thick, -{Latex[length=4.5pt,width=6.5pt]}] (-0.1,-2.15) -- (-0.1,-2.46);
  \draw [red, thick, -{Latex[length=4.5pt,width=6.5pt]}] (-0.1,-3.28) -- (-0.1,-2.98);
  \draw [red, thick] (-0.2,-2.46) -- (0.0,-2.46);
  \draw [red, thick] (-0.2,-2.98) -- (0.0,-2.98);
  \node at (-0.15, -2.7) {\color{red} \large $2\radius$};
  \end{scope}
 
\end{tikzpicture}}
\caption{\label{p:stokes-torus-eig}Condition numbers $\kappa(\StokesCF)$
  of the Stokes combined field integral operator
  $\StokesCF := I/2 + \StokesDL + \eta \StokesSL$
  (left), and
  the corresponding number of GMRES iterations for an exterior
  Stokes Dirichlet BVP solve to $\gmrestol=10^{-8}$ (middle),
  for three slenderness choices, and different scaling parameters $\eta>0$.
  {The predicted optimal values $\eta = 1/(2 \radius \log \radius^{-1} )$ are shown with vertical dashed lines.}
In the top row, the surface geometry is a torus with major radius $R=1$ and a fixed minor radius $\radius$.
In the bottom row, the minor radius varies from $\radius$ to $3\radius$.
  Geometries are shown on the right.
  The surface is discretized using eight $10$th-order Chebyshev panels, each with Fourier order of $20$.
  Note that the varying-radius case demands $\gmresiter$ around 8 times larger
  than the constant-radius case.
}
\end{figure}
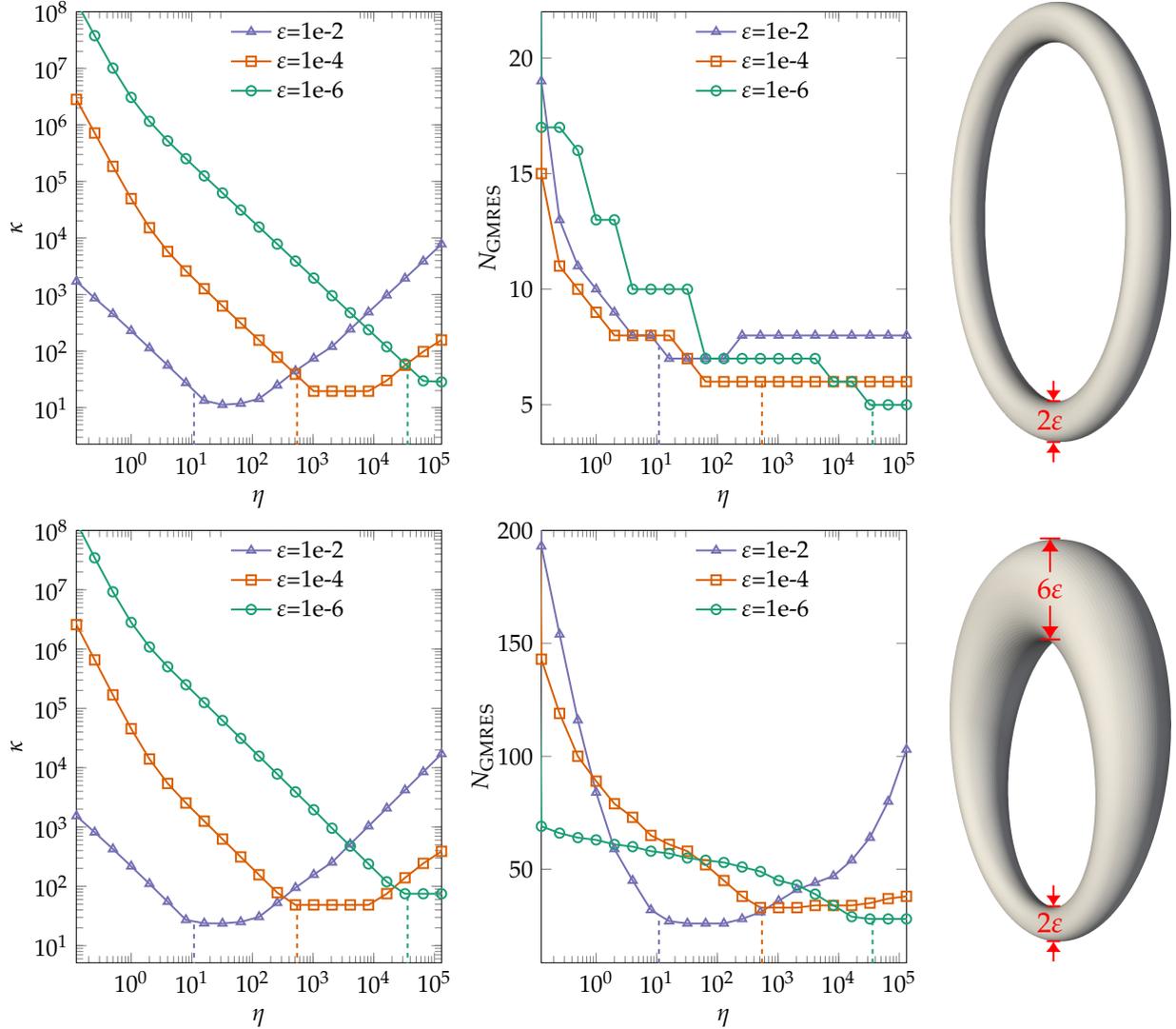

\paragraph{Stokes slender combined field integral equation formulation}
Combined field representations with mixing parameter $\eta\approx 1$
are already used for the Stokes exterior Dirichlet BVP,
\pr{e:stokes-eq0,e:stokes-eq1,e:stokes-eq2,e:stokes-dirichlet3},
for which their unique solvability has been proven
\cite{hsiao85,hebeker}.
However, inspired by the above Laplace analysis,
we propose a new admixture for the slender geometry case.
In $\Real^3 \setminus \overline\Omega$, the solution $\vct{u}$ is represented as
(recalling the velocity potentials \pr{e:stokes-sl,e:stokes-dl}),
\begin{align}
    \vct{u} &= \StokesDLR[\vct{\sigma}] + \StokesSLR[\vct{\sigma}] ~/~ (2 \radius \log \radius^{-1} ) & \text{in } \Real^3 \setminus \overline\Omega. \label{e:stokes-bie-rep-new}
\end{align}
Then, $\vct{u}$ satisfies the Stokes equations \pr{e:stokes-eq0,e:stokes-eq1,e:stokes-eq2} by construction.
Taking the limit to $\partial\Omega$ and applying the boundary conditions in \pr{e:stokes-dirichlet3}, we get the second-kind BIE
\begin{align}
    \StokesCF\vct{\sigma} &= \vct{u}_0  & \text{on } \partial\Omega, \label{e:stokes-bie-new}
\end{align}
where $\StokesCF := I/2 + \StokesDL + \StokesSL / (2 \radius \log \radius^{-1} )$, and
$\StokesSL$ and $\StokesDL$ are the (tensor-valued) boundary integral operators,
$\StokesDL$ being taken in the principal value sense.
This BIE can be solved efficiently using GMRES, and then \pr{e:stokes-bie-rep-new} can be used to evaluate $\vct{u}$ in $\Real^3 \setminus \overline\Omega$.

To test this claim, in \pr{p:stokes-torus-eig}, we return to true toroidal (closed filament)
geometries, and track the condition number $\kappa$ of the operator and
the number of GMRES iteration $\gmresiter$ required to solve the combined field integral equation
for different values of the single-layer parameter $\eta>0$.
We clearly observe minima in both the condition number and the number of GMRES iterations around $\eta = 1/(2 \radius \log \radius^{-1} )$,
matching our Laplace analysis.
We show similar results for more complicated geometries in \pr{ss:results-tangle}, for both Laplace and Stokes Dirichlet BVPs,
where GMRES converged in a few tens of iterations, instead of the hundreds of iterations for the standard $\eta=1$ admixture formulation.

\subsection{Stokes mobility problem \label{ss:mobility-formulation}} 
We turn to the Stokes mobility problem stated in \pr{s:problem-setup}.
We first present a standard 
``completed double-layer'' BIE formulation
\cite{Pozrikidis1992} with a slight variant of the completion flow.
The fluid velocity is represented in terms of the Stokes double- and single-layer potentials as,
\begin{align}
  \vct{u} &=
  \StokesDLR[\vct{\sigma}] +
  \StokesSLR[\vct{\CompletionFlowDen}]
  & \text{in } \Real^3 \setminus \Omega \label{e:mobility-dl-rep}
\end{align}
where $\vct{\sigma}$ is an unknown vector density field on $\partial\Omega$, and $\StokesSLR[\vct{\CompletionFlowDen}]$ is a ``completion flow''.
Some form of completion flow is necessary because the double-layer potential by itself can only represent flows with zero net force and zero net torque on each closed surface.
In the literature, completion flows using an interior point
stokeslet and rotlet dominate
\cite{power87,Pozrikidis1992,KK91book,afk16},
although in a slender body they would behave poorly
since they would give surface data that is only smooth on the $\radius$ scale.
The latter motivated the use of a line source completion flow
\cite{keaveny11}.
Instead, we propose a single-layer potential representation of the completion flow,
with density $\vct{\CompletionFlowDen}$ given by
the restriction to each surface of a certain rigid-body velocity field,
\begin{align}
    \vct{\nu}(\vct{x}) &= \vct\alpha_b + \vct{\beta}_b \times (\vct{x} - \vct{x}^c_b), & \vct{x} \in \partial\Omega_b , \; b=1,\dots,B.
\end{align}
Here $\vct\alpha_b\in\Real^3$ and $\vct\beta_b\in\Real^3$ are chosen such that
$\int_{\partial\Omega_b} \vct{\CompletionFlowDen}(\vct{x}) \,dS_\xx = \vct{F}_b$,
the given force,
and
$\int_{\partial\Omega_b} \vct{\CompletionFlowDen}(\vct{x}) \times (\vct{x} - \vct{x}^c_b) \,dS_\xx = \vct{T}_b$, the given torque.
This guarantees
{(see, e.g., \cite[(4.1.5--6)]{Pozrikidis1992})}
that the representation \pr{e:mobility-dl-rep} has
the required body forces and torques.
For each body independently, 
such a pair $(\vct\alpha_b,\vct\beta_b)$ is found by solving a simple
$6\times6$ linear system whose entries are approximated using the
surface quadrature used in the Nystr\"om method.
This form is convenient since it needs only existing BIE quadratures for $S$,
and the completion flow data
remains smooth with respect to $\params$ and $\paramt$ as $\radius\to 0$.

By construction, $\vct{u}$ in \pr{e:mobility-dl-rep}
satisfies the Stokes equations \pr{e:stokes-eq0,e:stokes-eq1,e:stokes-eq2}.
%
%
Taking the limit of \pr{e:mobility-dl-rep} to $\partial\Omega$, applying the boundary conditions \pr{rbm,e:mobility-bc}, and rearranging the terms, we get,
\begin{align}
  (I/2 + \StokesDL)\vct{\sigma} - \RigidBodyVel &= \vct{u}_s - \StokesSL\CompletionFlowDen & \text{on } \partial \Omega ,
  \label{e:mobility-bie-tmp}
\end{align}
where both terms on the right are known
(recall $\vct{u}_s$ is the given slip velocity on $\partial\Omega$).
On the left, $\vct{\sigma}$ is an unknown density
and $\RigidBodyVel$ is an unknown rigid body motion of the form
\pr{rbm}.
For $\Nobj$ rigid bodies, let $\RigidBodySpace$ be the $6\Nobj$ dimensional space of all rigid body motions on $\partial\Omega_b$, $b=1,\dots,\Nobj$
(\ie, it contains functions of the form $\vct{f}(\vct{x}) = \vct{v} + \vct{\omega} \times (\vct{x} - \vct{x}^c_b)$ on each $\partial\Omega_b$), and
let $\{\mathfrak{v}_1, \dots, \mathfrak{v}_{6\Nobj}\}$ be an orthonormal basis of $\RigidBodySpace$ in $L^2(\partial\Omega)$.
{To construct these orthonormal basis vectors, we build the three translation and three rotation
vectors for each rigid body and then use Gram–Schmidt orthogonalization to orthogonalize these six vectors.
The inner product in the Gram-Schmidt process is the $L^2$-inner product on the surface.
}
Then (as in \cite{corona3dmob}),
\be
L := \sum_{i=1}^{6\Nobj} \mathfrak{v}_i \mathfrak{v}_i^T
\label{L}
\ee
is the orthogonal projector
onto the subspace of rigid body motions in $L^2(\partial\Omega)$.
We will need the following well known fact that
$\RigidBodySpace$ is in
the null space of the Stokes exterior double-layer operator
\cite[Ch.~16, Thm.~4]{KK91book}; we include a concise proof in \pr{a:pf}.
\begin{pro}   
  \label{p:null}
  $(I/2 + \StokesDL)\vct{v} = \vct{0}$ for all $\vct{v} \in \RigidBodySpace$.
\end{pro}
Therefore, the solution $\vct{\sigma}$ to \pr{e:mobility-bie-tmp} is not unique and is determined only up to a vector in $\RigidBodySpace$.
Yet, we can use the extra degrees of freedom in $\vct{\sigma}$ and the fact that $\RigidBodyVel \in \RigidBodySpace$ to represent $\RigidBodyVel$ in terms of $\vct{\sigma}$ as,
\begin{align}
    \RigidBodyVel =  -L \vct{\sigma}. \label{e:recover-rigidbody-motion}
\end{align} 
Substituting in \pr{e:mobility-bie-tmp} gives us the following BIE formulation,
\begin{align}
  (I/2 + \StokesDL + L) \vct{\sigma} &= \vct{u}_s - \StokesSL\CompletionFlowDen & \text{on } \partial \Omega .
  \label{e:mobility-bie}
\end{align}
The boundary integral operator on the left side of \pr{e:mobility-bie} is
Fredholm and invertible;
this follows by Riesz--Fredholm theory \cite[Ch.~3-4]{LIE}
since it is the adjoint of the injective Fredholm operator arising in the interior
traction mobility formulation \cite[Lem.~6]{corona3dmob}.
Finally, after solving \pr{e:mobility-bie} for $\vct{\sigma}$, we can recover the rigid body motion $\RigidBodyVel$ from \pr{e:recover-rigidbody-motion}.

An issue with this formulation is that, as shown in \pr{p:cylinder-eigenvalues} for the Laplace case where $\kappa(I/2+D^L+{\bf 11}^T) = \bigO{\radius^{-2} / \log \radius^{-1}}$,
the Stokes exterior double-layer operator
$I/2+D$ is extremely ill-conditioned for slender-body geometries.
The cause---eigenvalue clustering around zero---is empirically the same.
Therefore, solving the BIE formulation in \pr{e:mobility-bie} with GMRES requires unreasonably large numbers of iterations.

\paragraph{Stokes slender mobility combined field integral equation formulation}
We now present a new BIE formulation for the Stokes mobility problem which gives well-conditioned discretizations for slender-body geometries.
The representation is
\begin{align}
  \vct{u} &= 
  \left( \StokesDLR + \StokesSLR / (2 \radius \log \radius^{-1} ) \right)
  \bigl[(I-L)\vct{\sigma}\bigr]
+ \StokesSLR[\vct{\CompletionFlowDen}]
  & \text{in } \Real^3 \setminus \overline\Omega . \label{e:new-mobility-rep}
\end{align}
Notice that this combines two features:
i) it replaces the double-layer in the above completed formulation
by the same carefully scaled admixture of the single- and double-layer potentials as in \pr{ss:dirichlet-bvp}, and
ii) the density $\vct{\sigma}$ is first projected to the space orthogonal to $\RigidBodySpace$ to remove any additional net force or net torque on any rigid body $\Omega_k$ due to the new single-layer source.
Taking the exterior limit of \pr{e:new-mobility-rep} to $\partial\Omega$ and applying the boundary condition in \pr{e:mobility-bc}, we get,
\begin{align}
  \StokesCF(I-L)\vct{\sigma}
  - \RigidBodyVel
  &= \vct{u}_s - \StokesSL\CompletionFlowDen & \text{on } \partial \Omega ,
  \label{e:new-mobility-bie-tmp}
\end{align}
where, as before, $\StokesCF := I/2 + \StokesDL + \StokesSL / (2 \radius \log \radius^{-1} )$ is the injective combined-field operator,
and $\RigidBodyVel$ is the unknown rigid body velocity.
As before, since $\vct{\sigma}$ is determined only up to a vector in $\RigidBodySpace$, and $\RigidBodyVel \in \RigidBodySpace$, we can choose the representation for $\RigidBodyVel$ in \pr{e:recover-rigidbody-motion} to get the new BIE,
\begin{align}
  \bigl( \StokesCF(I-L) + L \bigr) \vct{\sigma} 
  &= \vct{u}_s - \StokesSL\CompletionFlowDen   & \text{on } \partial \Omega.
  \label{e:new-mobility-bie}
\end{align}
The following shows that this proposed BIE is uniquely solvable.
Its proof is in \pr{a:pf}; note that the sign of the admixture parameter $\eta$ is crucial.
\begin{thm}  
  \label{t:BIE}
  Let $\Omega\subset\Real^3$ be the union of one or more smooth bounded bodies.
  Let $L$ be the orthogonal projector onto the space of rigid-body motions
  on $\partial\Omega$, as in \pr{L}.
  Let $\StokesCF := I/2 + \StokesDL + \eta\StokesSL$ be the Stokes
  exterior combined-field integral operator with parameter $\eta>0$ on $\partial\Omega$.
  Then
  \be
  \bigl( \StokesCF(I-L) + L \bigr) \vct{\sigma} = \vct{f}
  \label{BIE}
  \ee
  has a unique solution $\vct{\sigma}$
  for any right-hand side $\vct{f} \in C(\partial\Omega)^3$.
\end{thm}
After discretizing \eqref{e:new-mobility-bie}
and solving the linear system for $\vct{\sigma}$,
we can recover the rigid body motion using \pr{e:recover-rigidbody-motion}.
We expect the left-hand side operator in \pr{e:new-mobility-bie} to have a small condition number, comparable to the condition number of operator $\StokesCF$.
Therefore, it should converge rapidly when solved using iterative solvers like GMRES, uniformly as $\radius\to 0$.
This is confirmed through the numerical experiments presented next.

\section{Numerical results \label{s:results}}
We present numerical results to validate our method and demonstrate its efficiency.  We first compare the accuracy
and computational cost of our method with a generic boundary integral code in \cref{ss:results-biest}.  In
\cref{ss:results-tangle}, we solve Laplace and Stokes boundary value problems on a slender-body geometry and in
\cref{ss:results-nearly-touching}, we show similar results for close-to-touching geometries.
In \cref{ss:results-sbt}, we compare with slender-body theory (SBT) to show when SBT is a good approximation to the true solution and when it fails to give accurate results.
We solve the Stokes mobility problem in \cref{ss:results-mobility} and present parallel scalability results in \cref{ss:results-scalability}.
{All experiments use a discretization order of $\Ns=10$. All errors reported in this section are relative errors.}

\paragraph{Hardware} All numerical experiments were performed on the Skylake nodes of the Iron cluster at the Flatiron
Institute.  Each node has two 20-core Intel Xeon Gold 6148 CPUs running at 2.4GHz and 768GB of RAM.
All experiments, except those in \pr{ss:results-scalability}, were run on a single node (using up to $40$ cores).

\paragraph{Software} Our software was run on a Linux operating system and compiled using GCC-11.3.0, with the OpenMPI library version 4.0.7, and
optimization flags ``\texttt{-O3 -march=native}''.  We linked against the Intel MKL library version 2023.0.0 and the FFTW library version 3.3.10.
Most code is in C++, makes use of the first author's SCTL library, and
is publicly available at
~\url{https://github.com/dmalhotra/CSBQ}. See \cite{CSBQ-zenodo}.

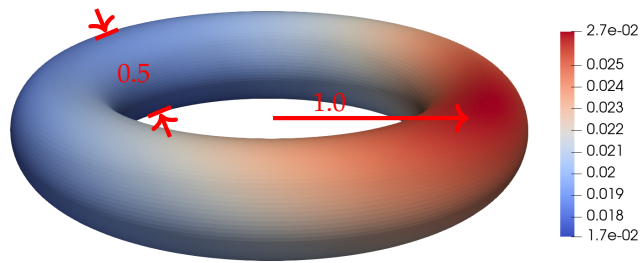
\begin{figure}[htb!]
  \begin{tikzpicture}
    \node[anchor=south west,inner sep=0] at (0,0) {\includegraphics[width=8.341cm]{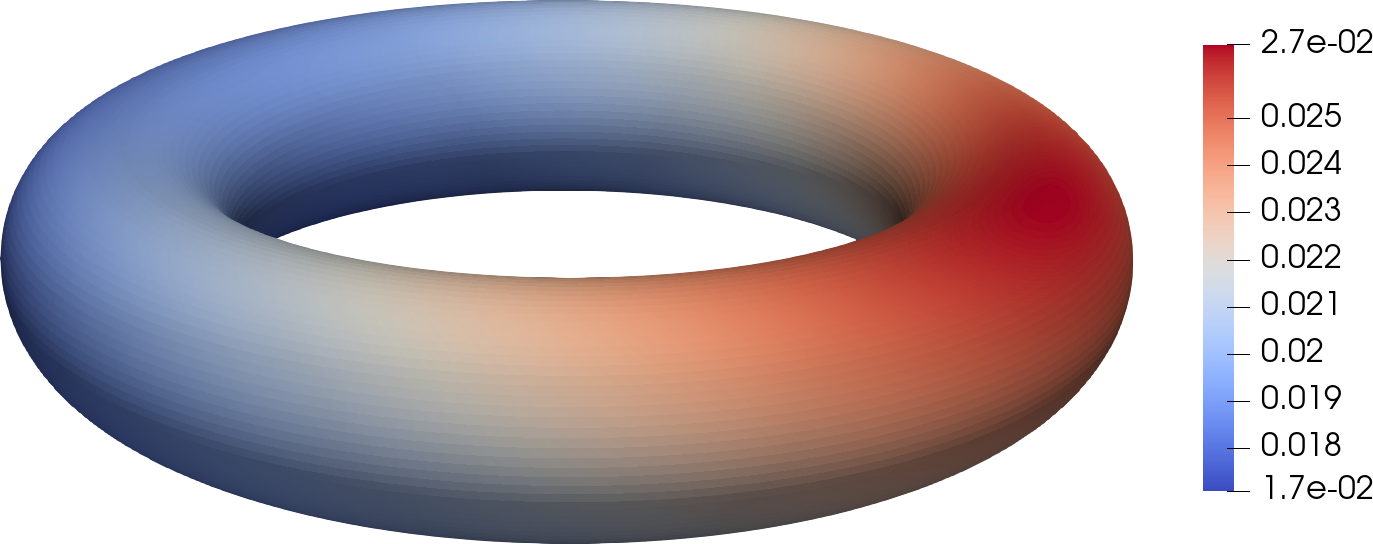}}; 
    \draw [red, ultra thick, ->|](1.15,3.35) -- (1.30,2.99);
    \node at (1.65, 2.50) {\color{red} $0.5$};
    \draw [red, ultra thick, ->|](2.13,1.65) -- (1.98,2.01);

    \draw [red, ultra thick, ->](3.5,1.9) -- (6.1,1.9);
    \node at (4.25, 2.12) {\color{red} $1.0$};
  \end{tikzpicture}
  \caption{\label{f:conv-biest}
    A torus with major radius $1.0$ and minor radius $\radius=0.25$ {(centered at the origin with axis in the Z-direction)} used for the convergence comparison in \cref{t:conv-biest}.
    The colors represent the reference harmonic potential $u$. It is generated by a somewhat distant point charge in the exterior {at coordinates (1,2,3)}.
  }
\end{figure}
\begin{table}[htb!] 
  \begin{tabular}{r r r r r r | c r r r r}
    \hline
    \multicolumn{6}{c|}{Our proposal (CSBQ)} & \multicolumn{5}{c}{Existing toroidal solver (BIEST)} \\
    $\Nunknown$ & {$\Nelem$} & {$\Nt$} & $\left\|e\right\|_{\infty}$ & \Tsetup & \Teval    & $~$ &   $\Nunknown$ &    $\left\|e\right\|_{\infty}$ & \Tsetup & \Teval \\
    \hline
    320  &  4 &  8 & 3.3e-03 &       0.016 &        0.0007 & $~$ &   $39\times13$ &    2.0e-03 &      0.1319 &     0.0017 \\
    720  &  6 & 12 & 7.6e-07 &       0.055 &        0.0015 & $~$ &   $63\times21$ &    4.0e-06 &      1.4884 &     0.0042 \\
    1280 &  8 & 16 & 1.1e-09 &       0.166 &        0.0037 & $~$ &   $87\times29$ &    4.3e-09 &      6.6825 &     0.0313 \\
    2000 & 10 & 20 & 2.1e-10 &       0.357 &        0.0083 & $~$ &  $111\times37$ &    3.5e-10 &     15.4711 &     0.0862 \\
    \hline
  \end{tabular}
  \caption{\label{t:conv-biest}
    Results for the Laplace Green's representation test for the geometry of \cref{f:conv-biest}. With increasing mesh refinement ($\Nunknown$ unknowns),
    we report the $L^\infty$-norm of the error in \pr{e:greens-rep} for the slender-body quadrature and for the BIEST code.
    For the same accuracy, the quadrature setup time $\Tsetup$ and evaluation time $\Teval$ are compared for both methods.
  }
\end{table}

\subsection{Comparison with a general boundary integral method \label{ss:results-biest}}
We first compare the accuracy and efficiency of our method
for Laplace BVPs against BIEST
\cite{Malhotra-jcp-2019,Malhotra-ppcf-2019}, a general boundary integral code for smooth surfaces of genus one.  BIEST uses a uniform
doubly-periodic grid to discretize the toroidal parameterization, with a smooth blending to a polar coordinate transform to compute the singular integral. Unlike our
method, BIEST does not require the cross-section to be circular; however, it cannot handle 
slender geometries
efficiently. To do a fair comparison of the two methods, we choose the geometry shown in \cref{f:conv-biest}, which can
be handled efficiently by both methods. We use the interior
Green's representation theorem
\cite[Thm.~6.5]{LIE} to validate the quadrature accuracy. It states that, for a
harmonic function $u$ on a bounded domain $\Omega \subset \Real^3$,
\[
u(\vct{x}) = \int_{\partial\Omega} {G}(\vct{x}-\vct{y}) \frac{\partial u(\vct{y})}{\partial\normal} dS_\yy - \int_{\partial\Omega} u(\vct{y}) \frac{\partial G(\vct{x}-\vct{y})}{\partial \normal_\yy} \, dS_\yy,
\qquad \text{for } \vct{x} \in \Omega,
\]
recalling ${G}(\vct{r}) = (4 \pi |\vct{r}|)^{-1}$, the free-space Green's function for Laplace's equation.
Taking the limit as $\vct{x} \to \partial\Omega$ with $\vct{x} \in \Omega$ and applying jump relations, we have
\begin{align}
  u(\vct{x}) &= \frac{u(\vct{x})}{2}
                  + \biggl(\LaplaceSL \frac{\partial u}{\partial\normal}\biggr)(\vct{x})
                  - (\LaplaceDL u)(\vct{x}),
  & \text{ for } \vct{x} \in \partial\Omega. 
    \label{e:greens-rep}
\end{align}
We evaluate the RHS in
\pr{e:greens-rep} using the two boundary integral methods and compare it to the reference potential~$u$ on~$\partial\Omega$.  The
reference potential $u$ (visualized in \cref{f:conv-biest}) is generated using a unit point charge at $\vct{x}_0$ outside of but near to the domain $\Omega$; ~ \ie,
~$u(\vct{x}) = {G}(\vct{x} - \vct{x}_0)$.

In \pr{t:conv-biest}, we show convergence in $L^\infty$-norm with mesh refinement for both
methods.
{For CSBQ, the geometry is discretized into $\Nelem$ uniform elements, each with the same discretization order $\Ns=10$ and $\Nt$.
The parameters for BIEST were chosen to match the errors of CSBQ and require about $2\times$ as many unknowns.}
We also report the quadrature setup time $\Tsetup$ and evaluation time $\Teval$ for each method. For the same accuracy,
the setup for slender-body quadrature is up to $40\times$ faster that for BIEST.
The evaluation time for both methods is dominated by the far-field computation, which has $\bigO{\Nunknown^2}$ cost.
The constant for the slender-body code is about half that of BIEST due to better vectorization.
Also note that the slender-body method supports FMM acceleration and should scale as \bigO{\Nunknown} for larger \Nunknown;
however, the FMM is not advantageous until $\Nunknown$ is in the tens-of-thousands.

\begin{figure}[htb!]
  \includegraphics[width=0.49\textwidth]{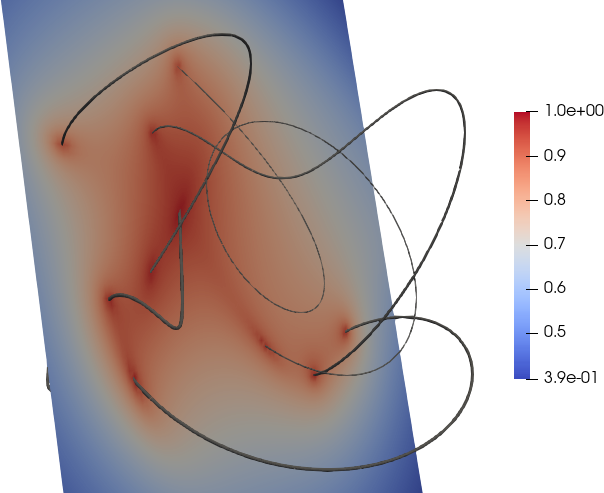}
  \includegraphics[width=0.49\textwidth]{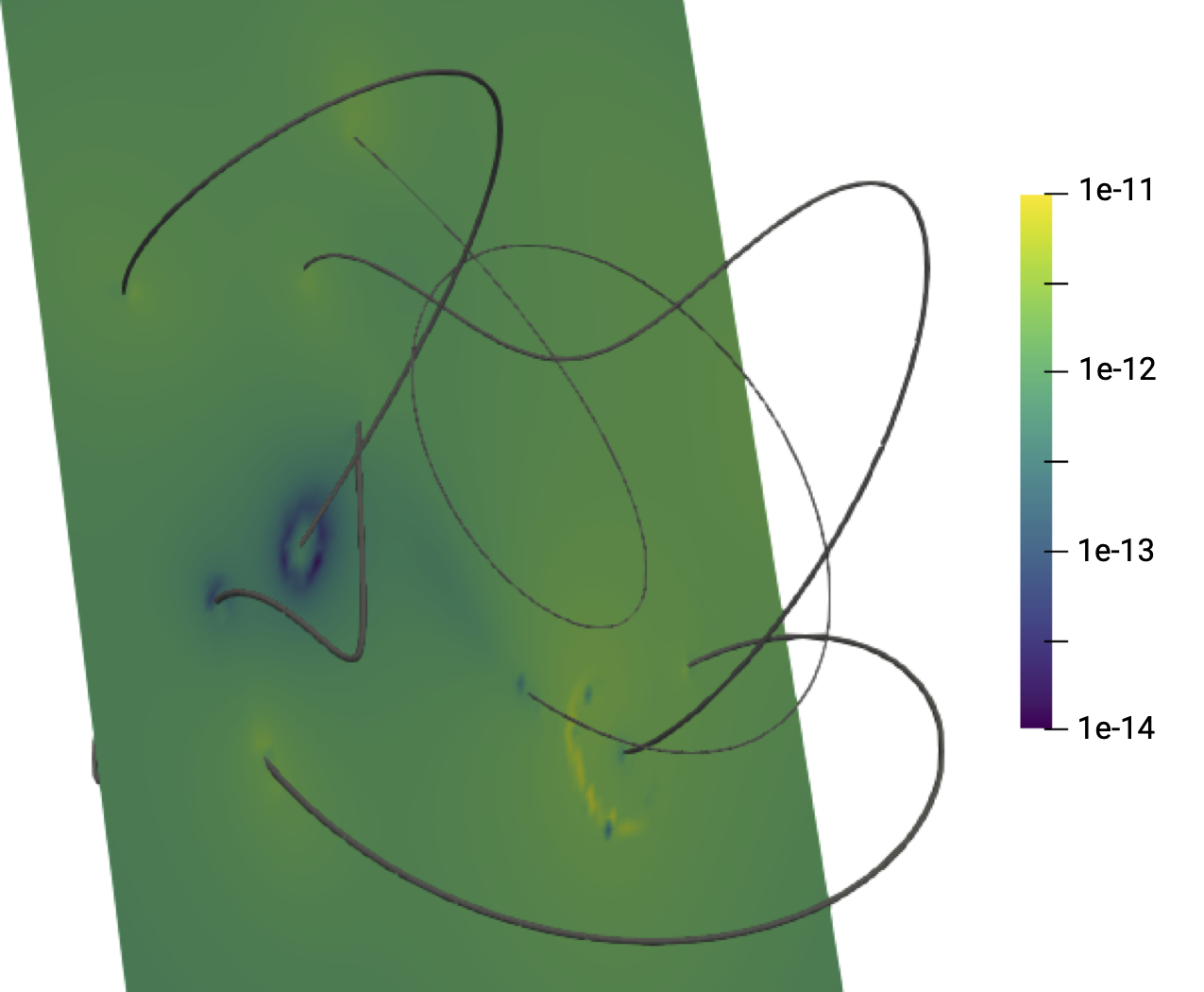}
  \caption{\label{f:conv-tangle-laplace}
    A closed fiber of length about 58, circular cross-sectional radius $\radius$ varying in the range $[0.005,0.015]$, hence aspect ratio $\sim10^3$,
    used for the convergence study in \cref{t:conv-tangle-laplace}.
    Left: the solution (visualized on a planar cross-section) corresponds to the potential in the exterior of a thin conducting wire maintained at a constant unit potential.
    Right: pointwise error magnitudes compared to a reference solution computed to much higher accuracy.
  }
\end{figure}
\begin{table}[htb!]
  \resizebox{0.99\textwidth}{!}{\begin{tabular}{c | r r r r r | r r | r r | r r}
    \hline
    &               &         &                 &           &          &            &                              &                \multicolumn{2}{c |}{1-core} & \multicolumn{2}{c }{40-cores} \\
    & $\Nunknown$   &  \Nelem & $\Nt^{\!\!max}$ & \gmrestol & $\eta$ & \gmresiter &  $\left\|e\right\|_{\infty}$ &   $\Tsetup~~(\Nunknown/\Tsetup)$ &  \Tsolve &   \Tsetup & \Tsolve \\
    \hline \parbox[t]{2mm}{\multirow{5}{*}{\rotatebox[origin=c]{90}{varying ~$\eta$}}}
    & 2.3e4         &   204   &              24 &     1e-08 &   3.1e-2 &    $>200$ &                      5.6e-06 &        1.41         ~~~~~~(1.6e4) &    104.6 &      0.12 &    5.96 \\
    & 2.3e4         &   204   &              24 &     1e-08 &   5.0e-1 &       63  &                      3.5e-07 &        1.41         ~~~~~~(1.6e4) &     32.8 &      0.12 &    1.82 \\
    & 2.3e4         &   204   &              24 &     1e-08 &    3.2e1 &       22  &                      1.2e-08 &        1.40         ~~~~~~(1.6e4) &     11.4 &      0.12 &    0.63 \\
    & 2.3e4         &   204   &              24 &     1e-08 &    2.0e3 &       46  &                      1.3e-08 &        1.41         ~~~~~~(1.6e4) &     23.9 &      0.12 &    1.33 \\
    & 2.3e4         &   204   &              24 &     1e-08 &    5.2e5 &       56  &                      1.0e-08 &        1.42         ~~~~~~(1.6e4) &     29.1 &      0.12 &    1.63 \\
    \hline \parbox[t]{2mm}{\multirow{6}{*}{\rotatebox[origin=c]{90}{no-FMM}}} 
    & 3.5e3         &    49   &               8 &     1e-02 &    3.2e1 &         4 &                      1.9e-02 &      0.093          ~~~~~~(3.7e4) &    0.029 &     0.025 &   0.009 \\
    & 8.6e3         &   103   &              12 &     1e-05 &    3.2e1 &        13 &                      1.9e-05 &      0.333          ~~~~~~(2.6e4) &    0.768 &     0.029 &   0.056 \\
    & 1.6e4         &   157   &              20 &     1e-07 &    3.2e1 &        19 &                      1.5e-07 &      0.806          ~~~~~~(1.9e4) &    4.423 &     0.071 &   0.275 \\
    & 2.8e4         &   227   &              24 &     1e-09 &    3.2e1 &        24 &                      5.3e-09 &      1.604          ~~~~~~(1.7e4) &   17.478 &     0.129 &   0.814 \\
    & 7.5e4         &   457   &              40 &     1e-11 &    3.2e1 &        30 &                      4.2e-11 &     11.190          ~~~~~~(6.7e3) &  160.174 &     1.742 &  10.026 \\
    & 1.6e5         &   893   &              48 &     1e-12 &    3.2e1 &        33 &                      1.9e-12 &     35.770          ~~~~~~(4.5e3) &  813.750 &     5.389 &  47.462 \\
    \hline \parbox[t]{2mm}{\multirow{6}{*}{\rotatebox[origin=c]{90}{FMM}}}  
    & 3.5e3         &    49   &               8 &     1e-02 &    3.2e1 &         4 &                      1.9e-02 &      0.181          ~~~~~~(1.9e4) &    0.140 &     0.051 &  0.194 \\
    & 8.6e3         &   103   &              12 &     1e-05 &    3.2e1 &        13 &                      1.9e-05 &      0.673          ~~~~~~(1.3e4) &    1.857 &     0.060 &  0.803 \\
    & 1.6e4         &   157   &              20 &     1e-07 &    3.2e1 &        19 &                      1.4e-07 &      1.614          ~~~~~~(9.6e3) &    8.649 &     0.148 &  2.199 \\
    & 2.8e4         &   227   &              24 &     1e-09 &    3.2e1 &        24 &                      7.7e-09 &      3.244          ~~~~~~(8.5e3) &   16.545 &     0.262 &  9.096 \\
    & 7.5e4         &   457   &              40 &     1e-11 &    3.2e1 &        30 &                      3.4e-11 &     22.113          ~~~~~~(3.4e3) &  120.745 &     3.445 & 31.406 \\
    & 1.6e5         &   893   &              48 &     1e-12 &    3.2e1 &        33 &                      6.8e-11 &     70.905          ~~~~~~(2.3e3) &  353.662 &    10.538 & 96.814 \\
    \hline
  \end{tabular}}
  \caption{\label{t:conv-tangle-laplace}
    Results for a Laplace Dirichlet boundary value problem shown in \cref{f:conv-tangle-laplace}.
    Here, \Nelem~ is the number of slender body elements, $\Nt^{\!\!max}$ is the maximum discretization order of the elements in the angular direction, $\Nunknown$ is the number of unknowns in the boundary integral equation, and $\gmrestol$ is the tolerance for the GMRES solve.
    The solution is computed using the Laplace combined field integral operator ($\LaplaceCF = I/2 + \LaplaceDL + \eta \LaplaceSL$), and we show results for different values of the mixing parameter $\eta$.
    We report the number of GMRES iterations $\gmresiter$, the quadrature setup time $\Tsetup$, the setup rate $\Nunknown/\Tsetup$ and the BIE solve time $\Tsolve$ for the serial implementation on 1 core and in parallel on 40 cores.
  }
\end{table}

\subsection{Dirichlet boundary value problems on slender geometry \label{ss:results-tangle}}
We demonstrate the performance of our method for the slender closed fiber in \pr{f:conv-tangle-laplace}. Its centerline $\xx_c(s)$ coordinates
are Fourier series with iid normal random sine and cosine coefficients
over the frequency index range $k=1,\dots,10$,
with decaying standard deviations $1/(1+k/3)$.
We choose a varying $\radius$ with typical aspect ratio
(circumference to length) of order $10^3$,
making it extremely challenging for prior boundary integral methods.
We solve a Laplace Dirichlet boundary value problem in the exterior of the domain.
The solution $u$ satisfies Laplace's equation $\Delta u = 0$ for all points $\vct{x} \in \Real^3 \setminus \overline\Omega$, has boundary condition
$u|_{\partial\Omega} \equiv 1$, and decay $u(\vct{x}) \rightarrow 0$ as $|\vct{x}| \rightarrow \infty$.
This corresponds to an electrostatic problem where a thin conducting wire loop is maintained at a constant unit potential.
We compute the numerical solution by solving the BIE formulation in \pr{e:laplace-bie-new}, for an unknown density $\sigma$ on the surface $\partial\Omega$.
Then, using \pr{e:laplace-bie-rep-new}, we evaluate the solution $u$ on a uniform 3D grid of dimensions $100\times100\times100$ enclosing $\Omega$.
We estimate the error by comparing it to a reference solution computed to much higher precision on a fine boundary mesh.
\pr{f:conv-tangle-laplace} shows the potential $u$ and the error in the computed solution on a planar cross-section.

In \pr{t:conv-tangle-laplace}, we present results for different boundary mesh resolutions (number of unknowns $\Nunknown$),
GMRES tolerance ($\gmrestol$), and mixing parameter $\eta$ (in $\LaplaceCF = I/2 + \LaplaceDL + \eta \LaplaceSL$).
In each case, we report the number of GMRES iterations in the linear solve (\gmresiter), and the maximum solution error on the grid ($\left\|e\right\|_{\infty}$) compared to the reference solution.
We also report the timings for the quadrature setup stage (\Tsetup) and the solve time (\Tsolve, excluding \Tsetup) on 1 core and 40 cores (using MPI parallelization)
when evaluated directly (``no-FMM'' case) and using FMM acceleration.

By varying the mixing factor $\eta$
we determined the optimal value to be around $\eta=32$;
this is similar to the expected $1/(2\radius \log \radius^{-1}) \approx 19$
using the lower end $0.005$ of the radius range for the geometry.
As discussed in \pr{s:bie}, $\eta$ affects the condition number of the discretized linear system and therefore the number of GMRES iterations (\gmresiter) required. 
The top rows of the table shows that
a poor choice (such as $\eta=\bigO{1}$) not only takes three
times longer to solve due to the iteration count, but also
has noticeably worse accuracy.

By refining the surface mesh and reducing the GMRES tolerance, we observe convergence to about 12 digits in the $L^\infty$ norm.
For 8-digit accuracy, we get quadrature setup rates ($\Nunknown/\Tsetup$) of 17K unknowns/s on 1 core, and 214K unknowns/s on 40 cores.
The solve time ($\Tsolve$) is dominated by the cost of the far-field evaluation.
With $\Nunknown$, this cost scales as $\bigO{\Nunknown^2}$ when evaluated directly and as $\bigO{\Nunknown}$ with FMM acceleration;
however, due to the larger constants in the FMM cost, this benefit is not apparent until $\Nunknown$ is large.
With our current implementation, when using the FMM, the single-layer and the double-layer operators have to be evaluated separately.
With some modifications, it would be possible to use a single combined-field kernel function (as we already do in the ``no-FMM'' case).
This would reduce $\Tsetup$ and $\Tsolve$ by a factor of two for the FMM accelerated case.


\begin{figure}[htb!]
  \includegraphics[width=0.49\textwidth]{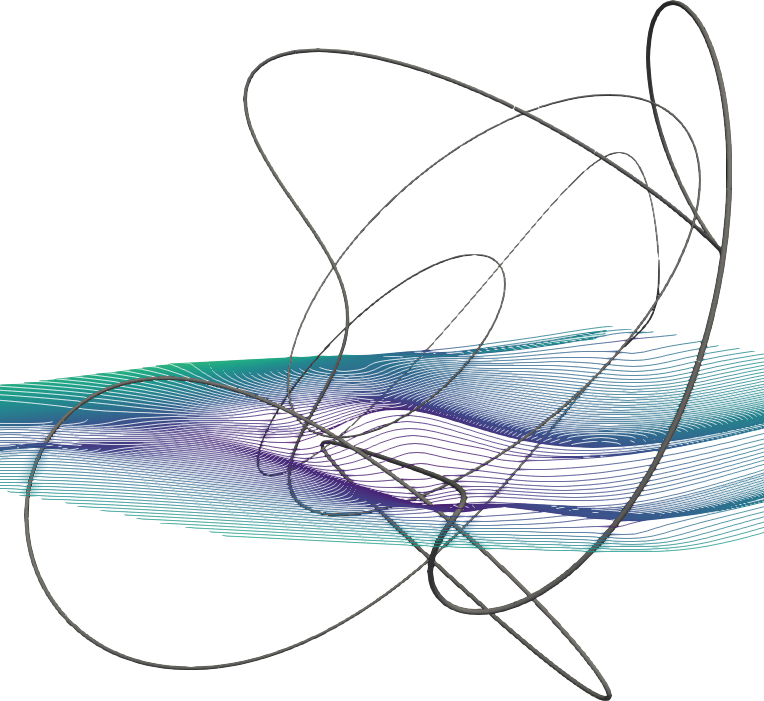}
  \includegraphics[width=0.49\textwidth]{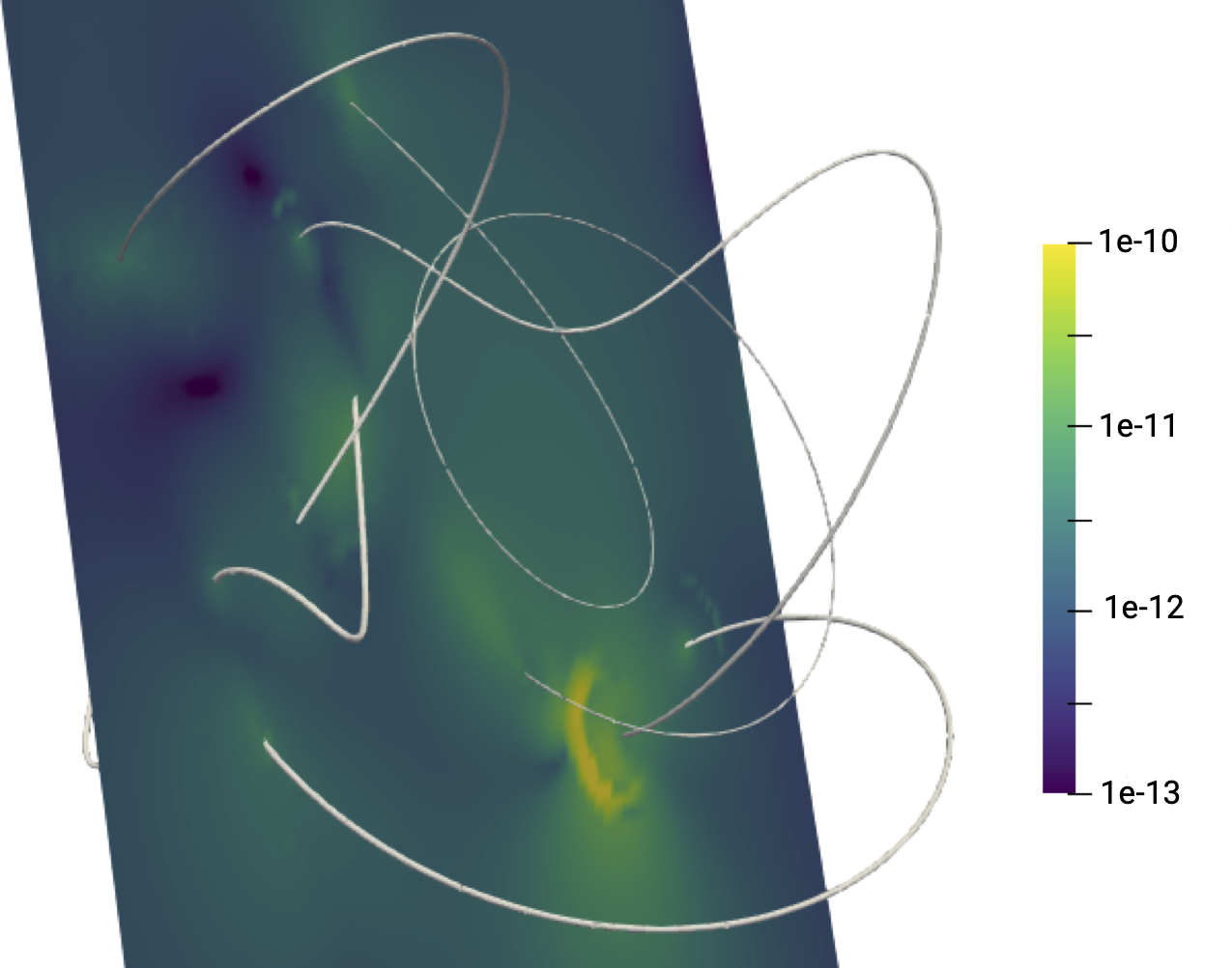}
  \caption{\label{f:conv-tangle-stokes}
Stokes BVP tests for the slender fiber geometry of
\pr{f:conv-tangle-laplace}, with convergence results in \cref{t:conv-tangle-stokes}.
    Left: the streamlines for the flow around the geometry. The flow corresponds to a thin rigid wire dragged through a stationary viscous fluid.
    Right: pointwise error magnitude compared to a reference solution computed to much higher accuracy.
  }
\end{figure}
\begin{table}[htb!]
  \resizebox{0.99\textwidth}{!}{\begin{tabular}{c | r r r r r | r r | r r | r r}
    \hline
    &               &         &                 &           &          &            &                              &                \multicolumn{2}{c |}{1-core} & \multicolumn{2}{c }{40-cores} \\
    & $\Nunknown$   &  \Nelem & $\Nt^{\!\!max}$ & \gmrestol & $\eta$ & \gmresiter &  $\left\|e\right\|_{\infty}$ &   $\Tsetup~~(\Nunknown/\Tsetup)$ &  \Tsolve &   \Tsetup & \Tsolve \\
    \hline \parbox[t]{2mm}{\multirow{5}{*}{\rotatebox[origin=c]{90}{varying ~$\eta$}}}
    & 8.3e4         &   227   &              24 &     1e-08 &   1.3e-1 &    $>200$ &                      1.6e-05 &        3.6          ~~~~~~(2.3e4) &    418.9 &      0.32 &   20.21 \\
    & 8.3e4         &   227   &              24 &     1e-08 &  2.0~~~~ &       95  &                      2.2e-07 &        3.6          ~~~~~~(2.3e4) &    198.1 &      0.35 &    9.80 \\
    & 8.3e4         &   227   &              24 &     1e-08 &    6.4e1 &       38  &                      4.5e-08 &        3.6          ~~~~~~(2.3e4) &     78.9 &      0.32 &    3.76 \\
    & 8.3e4         &   227   &              24 &     1e-08 &    2.0e3 &       63  &                      4.0e-08 &        3.6          ~~~~~~(2.3e4) &    130.9 &      0.32 &    6.12 \\
    & 8.3e4         &   227   &              24 &     1e-08 &    1.3e5 &    $>200$ &                      1.8e-07 &        3.6          ~~~~~~(2.3e4) &    418.8 &      0.32 &   20.13 \\
    \hline \parbox[t]{2mm}{\multirow{6}{*}{\rotatebox[origin=c]{90}{no-FMM}}}
    & 1.0e4         &    49   &               8 &     1e-02 &    6.4e1 &         5 &                      3.5e-02 &    0.193          ~~~~~~(5.4e4) &     0.130 &    0.042 &   0.017 \\
    & 2.6e4         &   103   &              12 &     1e-05 &    6.4e1 &        22 &                      5.5e-05 &    0.572          ~~~~~~(4.5e4) &     4.039 &    0.045 &   0.215 \\
    & 4.7e4         &   157   &              20 &     1e-07 &    6.4e1 &        33 &                      6.6e-07 &    1.416          ~~~~~~(3.3e4) &    19.518 &    0.134 &   1.162 \\
    & 8.3e4         &   227   &              24 &     1e-08 &    6.4e1 &        38 &                      4.5e-08 &    3.623          ~~~~~~(2.3e4) &    78.907 &    0.324 &   3.689 \\
    & 2.2e5         &   457   &              40 &     1e-10 &    6.4e1 &        49 &                      2.9e-10 &   21.949          ~~~~~~(1.0e4) &   746.966 &    4.458 &  48.494 \\
    & 4.8e5         &   893   &              48 &     1e-11 &    6.4e1 &        54 &                      2.4e-11 &   84.363          ~~~~~~(5.7e3) &  3788.948 &   15.177 & 227.747 \\
    \hline \parbox[t]{2mm}{\multirow{6}{*}{\rotatebox[origin=c]{90}{FMM}}}
    & 1.0e4         &    49   &               8 &     1e-02 &    6.4e1 &         5 &                      3.5e-02 &    0.373          ~~~~~~(2.8e4) &     0.335 &    0.090 &    0.603 \\
    & 2.6e4         &   103   &              12 &     1e-05 &    6.4e1 &        22 &                      5.6e-05 &    1.186          ~~~~~~(2.2e4) &     7.059 &    0.096 &    8.412 \\
    & 4.7e4         &   157   &              20 &     1e-07 &    6.4e1 &        33 &                      4.0e-06 &    2.833          ~~~~~~(1.6e4) &    31.694 &    0.287 &   22.120 \\
    & 8.3e4         &   227   &              24 &     1e-08 &    6.4e1 &        38 &                      7.2e-08 &    7.140          ~~~~~~(1.2e4) &    88.176 &    0.648 &   57.439 \\
    & 2.2e5         &   457   &              40 &     1e-10 &    6.4e1 &        49 &                      1.2e-09 &   44.177          ~~~~~~(5.1e3) &   488.440 &    8.747 &  216.993 \\
    & 4.8e5         &   893   &              48 &     1e-11 &    6.4e1 &        54 &                      1.4e-10 &  167.555          ~~~~~~(2.9e3) &  1818.084 &   30.235 &  387.249 \\
    \hline
  \end{tabular}}
  \caption{\label{t:conv-tangle-stokes}
    Results for a Stokes Dirichlet boundary value problem shown in \cref{f:conv-tangle-stokes}. For different surface discretizations (\Nelem, $\Nt$), number of unknowns $\Nunknown$, and GMRES tolerance $\gmrestol$,
    we report the number of iterations $\gmresiter$, the quadrature setup time $\Tsetup$, the setup rate $\Nunknown/\Tsetup$ and the total BIE solve time $\Tsolve$ for the serial implementation on 1 core and in parallel on 40 cores.
  }
\end{table}

In \pr{f:conv-tangle-stokes} and \pr{t:conv-tangle-stokes}, we present corresponding results for a Stokes Dirichlet (resistance) boundary value problem for this same fiber.
The solution $\vct{u}$ satisfies the Stokes equations in the exterior of the wire loop: ~$\Delta \vct{u} - \nabla p = \vct{0}$,~ ~$\nabla \cdot \vct{u} = 0$~ for all points $\vct{x} \in \Real^3 \setminus \overline\Omega$,
with boundary conditions ~$\vct{u}|_{\partial\Omega} = \vct{u}_0$
where $\vct{u}_0 = (1,1,1)$ is constant,
and~ $|\vct{u}(\vct{x})| \rightarrow 0$ as $|\vct{x}| \rightarrow \infty$. This corresponds to a rigid body dragged through a stationary viscous fluid.
We solve the boundary integral equation and evaluate the velocity field $\vct{u}$ on a $100\times100\times100$ grid as before.
The $L^\infty$-norm errors on the grid (compared to a reference solution, computed to much higher accuracy) are reported.
The quadrature setup time, setup rate and the BIE solve times are again reported for 1 core and 40 cores, with and without using FMM acceleration.
For single precision accuracies (about 7 digits), we find setup rate of about 23K unknowns/s on one core and over 255K unknowns/s on 40 cores.

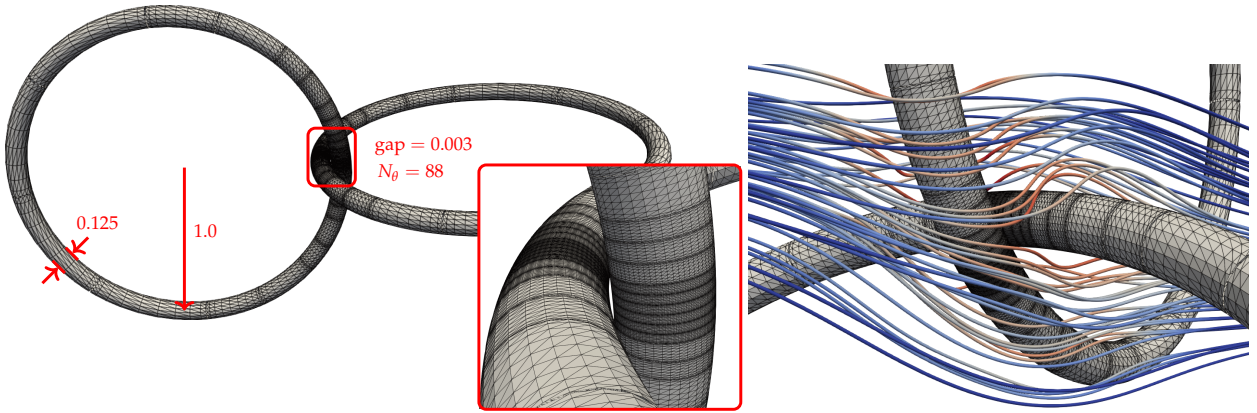
\begin{figure}[htb!]
  \resizebox{0.59\textwidth}{!}{\begin{tikzpicture}
      \node[anchor=south west,inner sep=0] at (0,0) {\includegraphics[width=12.65cm]{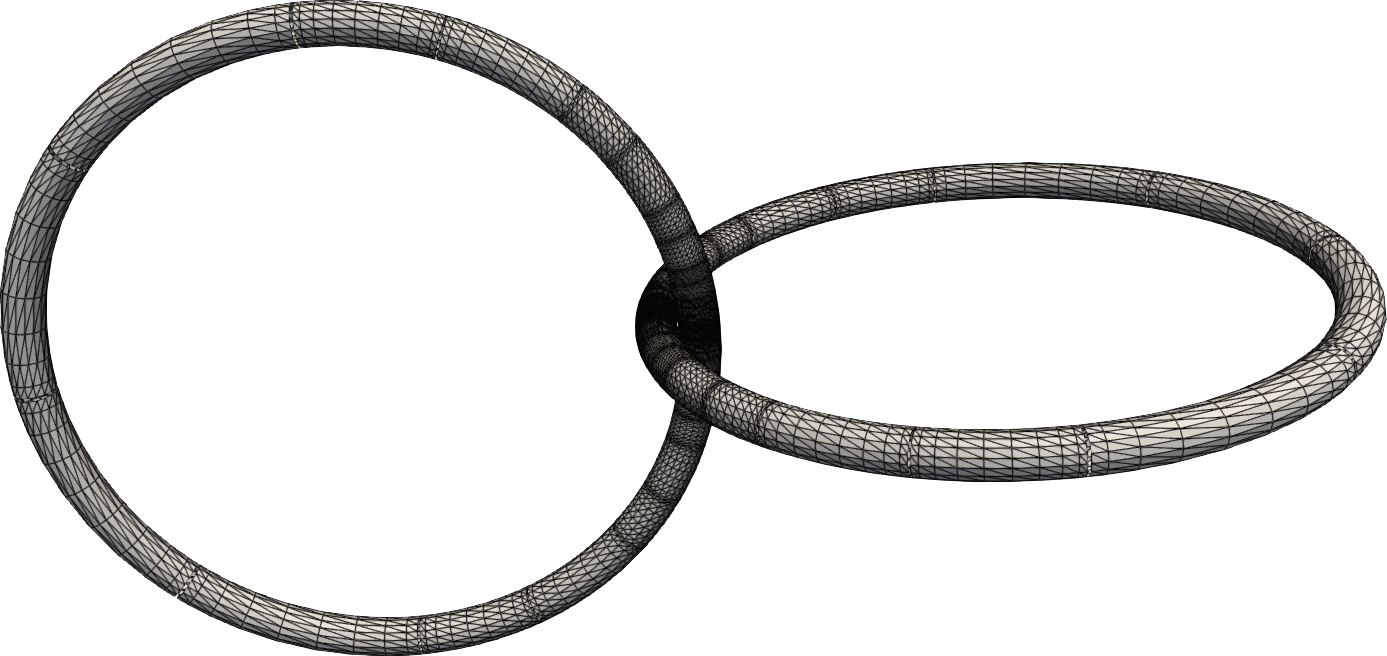}};
      \node[anchor=south west,inner sep=0] at (9,-1.7) {\adjincludegraphics[height=4.62cm,trim={0 0 {0.11\width} 0},clip]{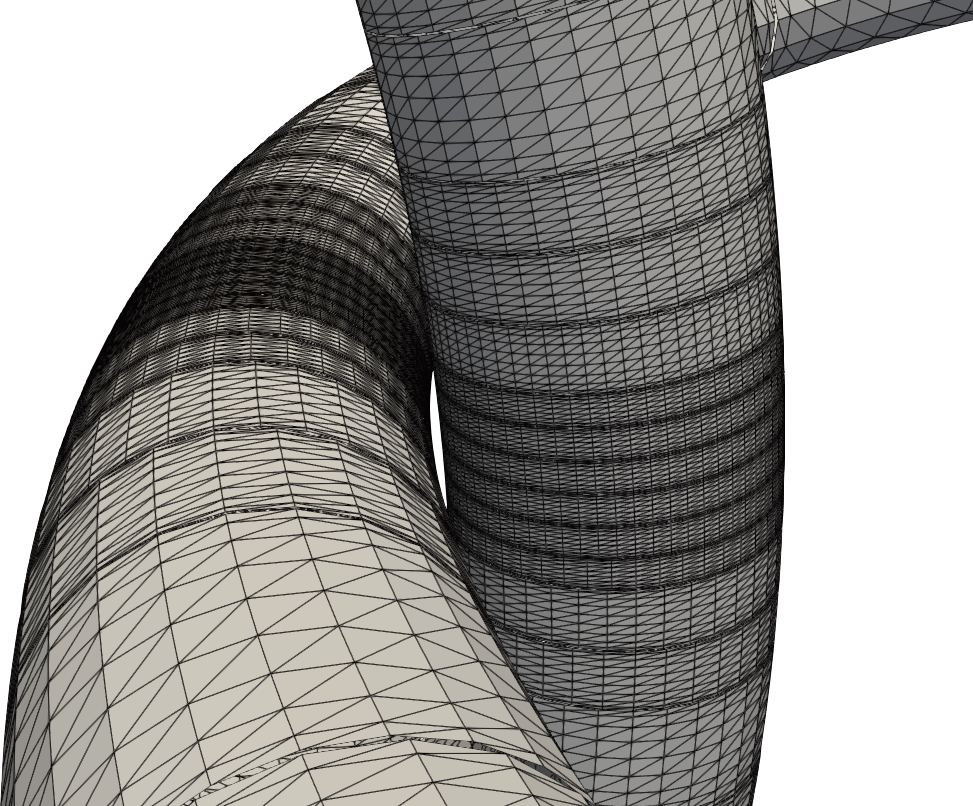}};
      \draw[red,ultra thick,rounded corners] (5.75,2.55) rectangle (6.65,3.65);

      \draw[red,ultra thick,rounded corners] (9,-1.7) rectangle (13.98,2.95);

      \draw [red, ultra thick, ->|](0.7,0.7) -- (1.03,1.03);
      \draw [red, ultra thick, ->|](1.57,1.57) -- (1.24,1.24);
      \node at (1.75, 1.85) {\color{red} $0.125$};

      \draw [red, ultra thick, ->](3.4,2.9) -- (3.4,0.18);
      \node at (3.8, 1.7) {\color{red}  $1.0$};

      \node at (7.95, 3.3) {\color{red}  gap $= 0.003$};
      \node at (7.7, 2.8) {\color{red}  $N_\theta = 88$};

  \end{tikzpicture}}
  \includegraphics[width=0.4\textwidth]{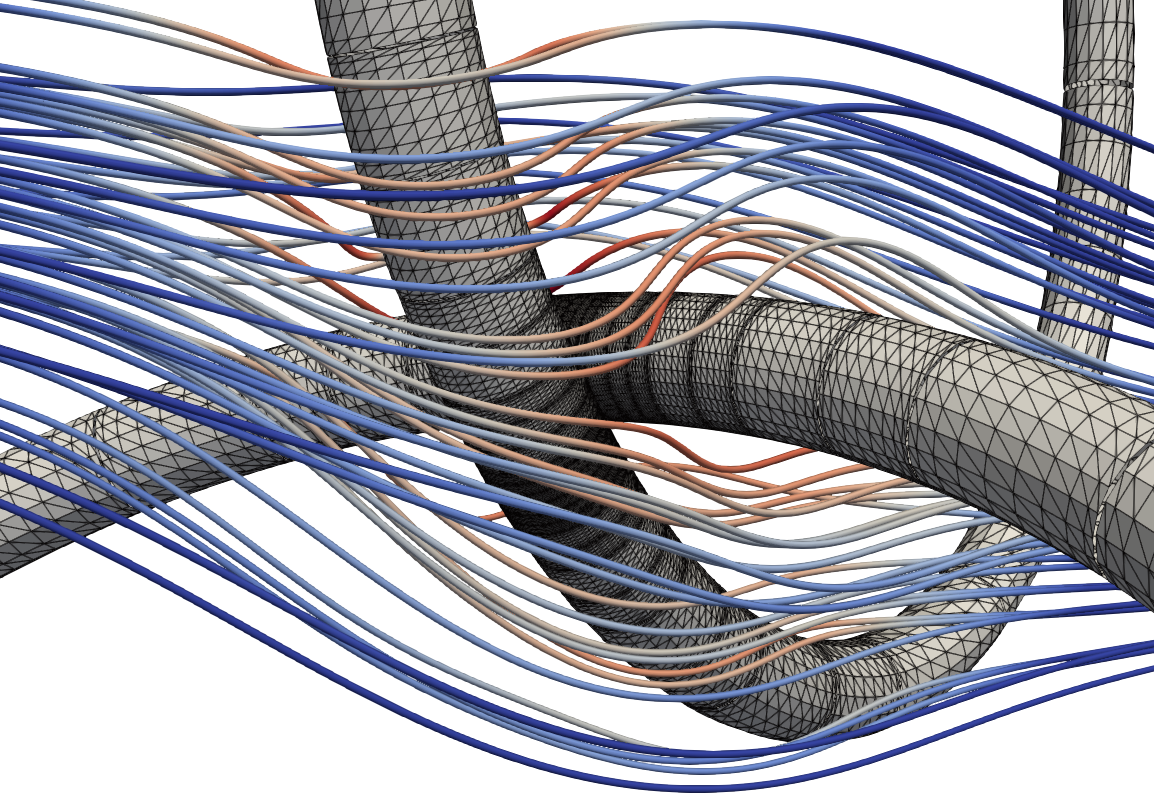}
  \caption{\label{f:conv-touching-stokes}
    Stokes flow around two tori with a narrow separation between them
    (see plot for parameters).
  Left: the surface discretization mesh used to solve the BIE for the Stokes Dirichlet BVP. The finest elements (see zoomed view)
  require $\Nt = 88$ points in the angular direction to resolve the interactions in the gap.
  Right: streamlines showing the flow around the geometry.
  }
\end{figure}
\begin{table}[htb!]
  \begin{tabular}{r r r | r r | r r | r r}
    \hline
                &            &             &            &                              &        \multicolumn{2}{c |}{1-core} & \multicolumn{2}{c }{40-cores} \\
    $\Nunknown$ & $\quadtol$ & $\gmrestol$ & \gmresiter &  $\left\|e\right\|_{\infty}$ &   $\Tsetup~~(N/\Tsetup)$ &  \Tsolve &   \Tsetup & \Tsolve \\
    \hline
    6.5e4       &     1e-03  &       1e-02 &          4 &        2.0e-02 &          6.8   ~~~~~~(9.5e+3) &          4.1 &          1.10 &         0.5 \\
    6.5e4       &     1e-06  &       1e-05 &         21 &        1.5e-05 &         12.3   ~~~~~~(5.3e+3) &         28.0 &          1.80 &         3.6 \\
    6.5e4       &     1e-08  &       1e-07 &         31 &        3.1e-07 &         16.5   ~~~~~~(3.9e+3) &         48.5 &          2.24 &         6.2 \\
    6.5e4       &     1e-12  &       1e-10 &         45 &        3.5e-10 &         32.4   ~~~~~~(2.0e+3) &         82.3 &          3.39 &        10.0 \\
    6.5e4       &     1e-14  &       1e-12 &         52 &        6.7e-12 &         47.5   ~~~~~~(1.4e+3) &        108.1 &          4.07 &        12.2 \\
    \hline
  \end{tabular}
  \caption{\label{t:conv-touching-stokes}
    Convergence for a Stokes Dirichlet boundary value problem for close-to-touching tori as in \cref{f:conv-touching-stokes}.
    The error $\left\|e\right\|_{\infty}$ converges to about 10 digits as we reduce the GMRES tolerance $\gmrestol$ and the quadrature accuracy tolerance $\quadtol$.
    We also report the quadrature setup time $\Tsetup$ and the BIE solve time $\Tsolve$ for the serial case on 1 core and in parallel on 40 cores.
  }
\end{table}

\subsection{Close-to-touching interactions \label{ss:results-nearly-touching}}

In order to demonstrate that our slender BIE quadrature remains accurate
even in the presence of close-to-touching surfaces,
we test the geometry of \pr{f:conv-touching-stokes} with
two rings (tori) of unit major radius, minor radius
$\radius=0.0625$, and separation less than $\radius/20$.
Note that slender body theory would give meaningless answers in this case,
since it already breaks down at separations of order $\radius$
(see next section)---here we are an order of magnitude closer yet.
A true PDE solve is required, and it is difficult to imagine
a method other than a boundary integral equation achieving the
accuracy that we now exhibit.
The mesh is refined adaptively with $\Nunknown=65K$ unknowns {($\Nelem=62$ elements)} and a maximum Fourier discretization order $\Nt^{\!\!max}=88$.
We solve the Stokes Dirichlet boundary value in the exterior
of the rings with boundary conditions $\vct{u}|_{\partial\Omega} = \vct{u}_0$
where $\vct{u}_0 = (1,1,1)$.
In \pr{t:conv-touching-stokes}, we report maximum error
compared to a reference
solution computed to much higher accuracy evaluated on a $100\times100\times100$ grid as before.
We observe convergence to about 11 digits
in $L^\infty$-norm as we reduce the quadrature accuracy tolerance $\quadtol$
and the GMRES tolerance $\gmrestol$.

Most of the unknowns are in the close-touching region, due to panel adaptivity
in $\params$ along the centerline, and in the
angular discretization order $\Nt$ for each resulting slender element.
The high resulting $\Nt$ values for many of the elements makes the quadrature setup relatively expensive, since it scales quadratically with $\Nt$.
For 6-digits of accuracy we get a setup rate of 3.9K unknowns/s on 1 core and 29K unknowns/s in parallel on 40 cores. These are about 10 times slower than
for the previous (non-close) geometries tested.

\begin{rmk}  
  The tori geometries and separation in \pr{f:conv-touching-stokes}
  are similar to those in the recent BIE tests of
  \cite[Sec.~6.3]{hedgehog} using
``hedgehog'' quadrature (not designed for slender bodies, nor parallelized),
which reported 6 CPU hours at 5-digit accuracy.
While a strict comparison is not meaningful,
we note that our CSBQ 1-core solve time of 40 seconds at this accuracy
is roughly three orders of magnitude faster.
\end{rmk}

\begin{figure}[htb!]  
  \begin{minipage}[t]{.6\linewidth}
    \centering
    \begin{tabular}{|r r | r r |}
      \hline
      $\radius$ &     $U_{exact}$ & Error-BIE & Error-SBT \\
      \hline
      1e-1      & 6.14921383598558e-2 &   0.1e-12 &   0.5e-02 \\
      1e-2      & 9.09845223245838e-2 &   0.9e-12 &   0.1e-03 \\
      1e-3      & 1.20156558899037e-1 &   0.6e-14 &   0.2e-05 \\
      1e-4      & 1.49319329075867e-1 &   0.2e-13 &   0.2e-07 \\
      1e-5      & 1.78481913130970e-1 &   0.1e-13 &   0.3e-09 \\
      \hline
    \end{tabular}
    \captionof{table}{\label{t:conv-sbt}
      Sedimentation speed of a ring of major radius 1 and cross-sectional (minor) radius $\radius$ subjected to an axial force $(0,0,-1)$, in a unit-viscosity Stokes fluid.
      We report the relative errors in the downward velocity computed using our boundary integral method (Error-BIE) and using slender body theory (Error-SBT). $U_\tbox{exact}$
      is taken from the reciprocals of the 3rd column of \cite[Tbl.~2]{mitchell2022},
      divided by $6\pi(1+\radius)$, and is accurate to at least 13 digits.
      On the right is shown the geometry for $\eps=10^{-1}$,
      with a sheet of solution streamlines.
    }
  \end{minipage}
  \begin{minipage}[t]{.39\linewidth}
    \centering
    $\vcenter{\includegraphics[width=0.99\textwidth]{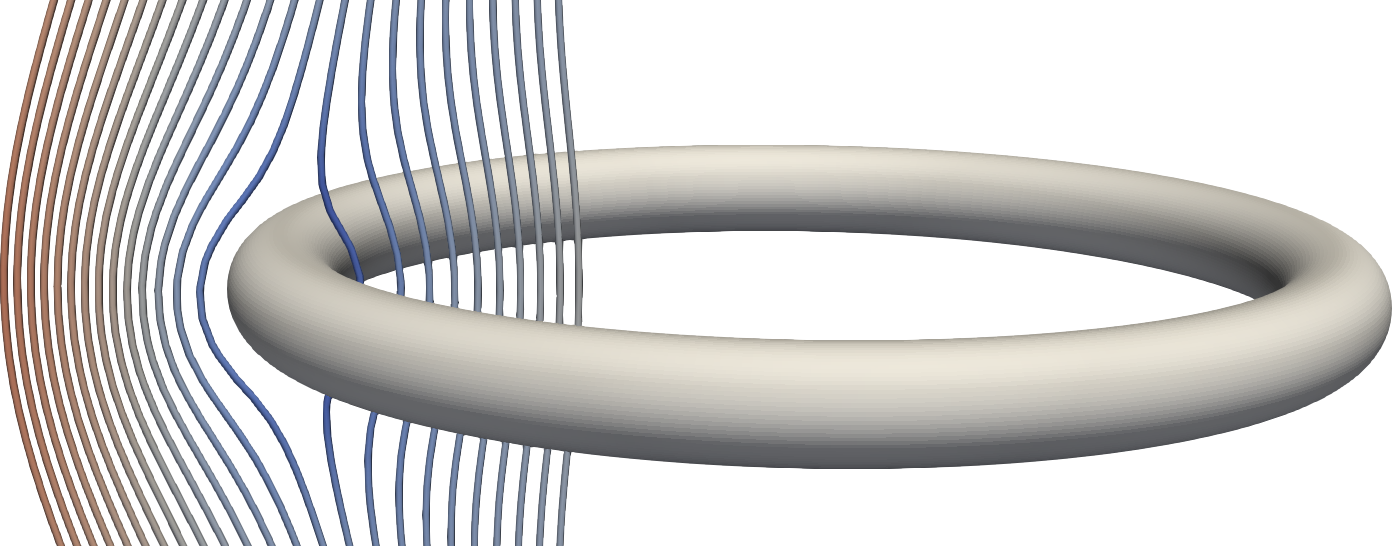}}$
  \end{minipage}
\end{figure}

\subsection{Comparisons with slender body theory \label{ss:results-sbt}}

Here we compare numerical solutions based on an SBT asymptotic approximation
to converged solutions of the true Stokes Dirichlet BVP.
Indeed, the present CSBQ method enables one to quantify
errors in using SBT as a numerical tool for \textit{rigid} body
hydrodynamics, which has not yet received numerical study in the close-touching
case. (See \cite{keaveny11} for a comparison for a helical fiber
without close-touching interactions, and \cite{mitchell2022} for
a close-touching flexible loop).
The only relevant rigorous error bound is the recent special case 
of a rigid straight periodic fiber, with regularized SBT kernel,
for which a $\bigO{\eps^2}$ error bound in the $L^2$-norm was given
for $H^2$-regular velocity data \cite{mori21inv}.

Following \cite{mitchell2022},
we start with a simple case without close-to-touching issues:
the axial mobility (sedimentation)
of tori of varying $\radius$, shown in \pr{t:conv-sbt}.
Due to symmetry, all functions of $\params$ are constant
and the SBT formulation (see \pr{a:sbt}; here $K\vct{f}$ and $\shat^T\vct{f}$ vanish)
degenerates into a single scalar ratio between force and velocity.
In \pr{t:conv-sbt} we show numerical results for rings (tori)
of unit radius and varying cross-sectional radius $\radius$, when a unit downward force $(0,0,-1)$ is applied to them.
The true downward velocity $(0,0,-U_{exact})$ of the rings is obtained from the drag coefficients (computed to at least 13 digits using a semi-analytic approach) reported in Table 2 of \cite{mitchell2022}.
We compare this with the velocity computed using SBT, and that using our boundary integral (CSBQ) method.
The latter exhibits a relative accuracy of at least 12 digits for all values of $\radius$, which serves as an independent validation of CSBQ.
The error of SBT agrees extremely well with expectations,
namely its leading asymptotic omitted term of $\bigO{\radius^2 \log \radius^{-1}}$ \cite{johnson80}.


\begin{figure}[htb!]
  $\vcenter{\hbox{\resizebox{0.53\textwidth}{!}{\begin{tikzpicture}
    \node[anchor=south west,inner sep=0] at (0,0) {\includegraphics[width=10cm]{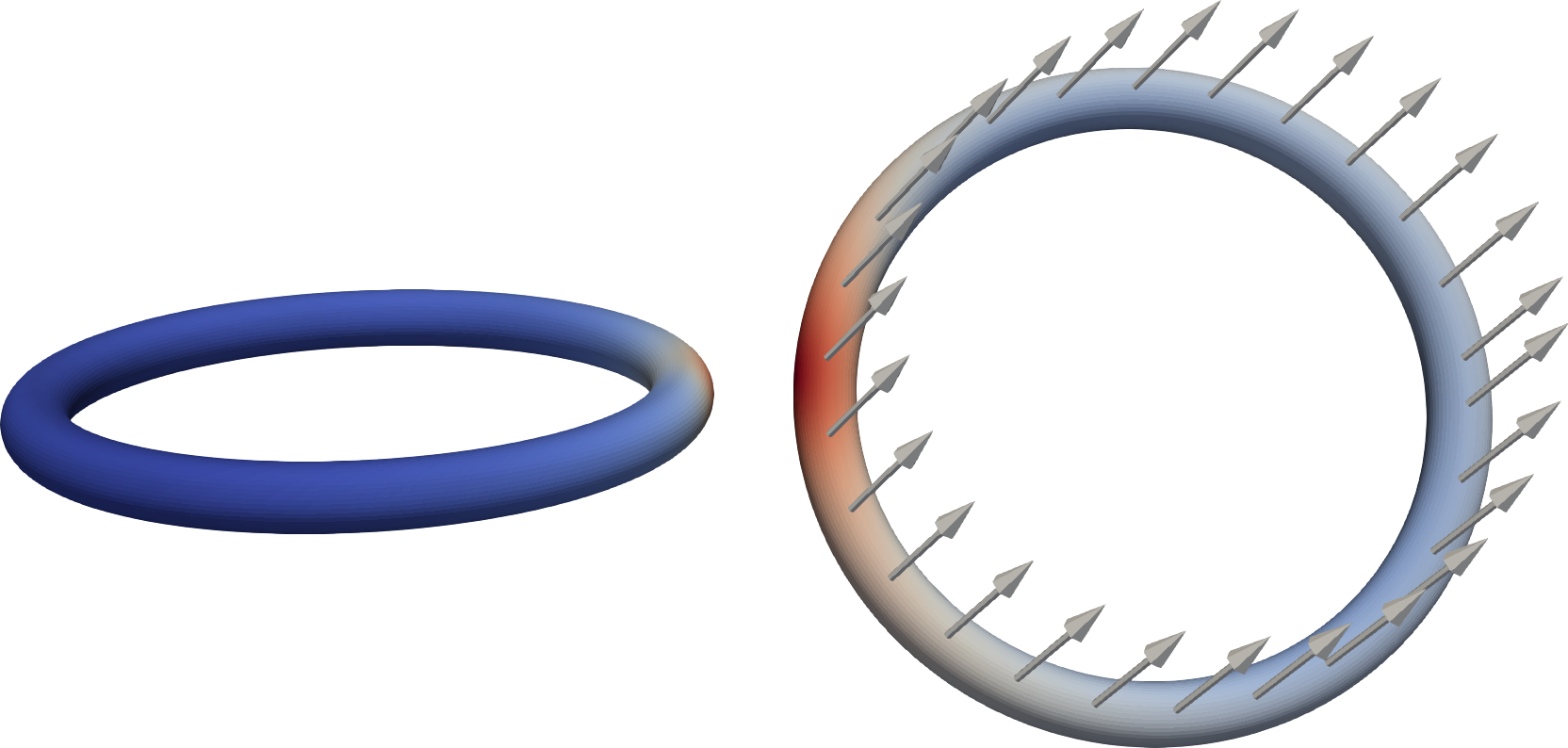}};
    
%
%
    
    \draw [red, thick, -{Latex[length=4.5pt,width=6.5pt]}] (2.3,2.21) -- (2.3,2.56);
    \draw [red, thick, -{Latex[length=4.5pt,width=6.5pt]}] (2.3,3.28) -- (2.3,2.93);
    \draw [red, thick] (2.4,2.56) -- (2.2,2.56);
    \draw [red, thick] (2.4,2.93) -- (2.2,2.93);
    \node at (2.65, 3.2) {\color{red} \large $2\radius$};
    
    \draw [red, thick, {Latex[length=4.5pt,width=6.5pt]}-{Latex[length=4.5pt,width=6.5pt]}] (4.57,2.3) -- (5.07,2.3);
    \draw [red, thick] (4.57,2.2) -- (4.57,2.4);
    \draw [red, thick] (5.07,2.2) -- (5.07,2.4);
    \node at (4.82, 2.6) {\color{red} \large $\delta$};
  \end{tikzpicture}}}}$%
  $\vcenter{\hbox{\resizebox{0.47\textwidth}{!}{\begin{tikzpicture}
    \node[anchor=south west,inner sep=0] at (0,0) {\includegraphics[width=8cm,trim = 1.9cm 1.1cm 1.9cm 1.1cm, clip]{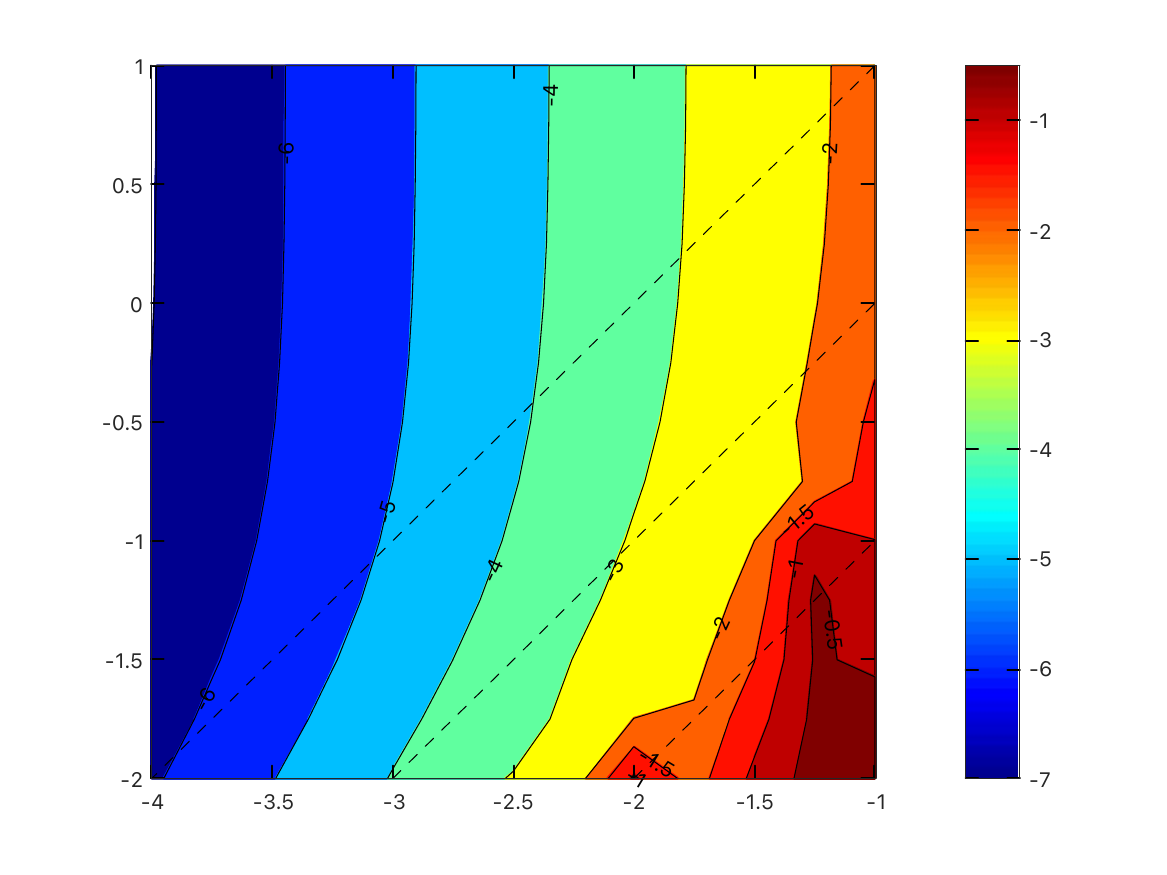}};
    \node at (3.4, 6.6) {$\log_{10} \| \vct{F}_{SBT} - \vct{F}_{BIE} \|_2 / \| \vct{F}_{BIE} \|_2$};
    
    \node at (3.30, 4.0) {$\delta = 100\radius$};
    \node at (3.30, 1.9) {$\delta = 10\radius$};
    \node at (4.30, 0.7) {$\delta = \radius$};
    
    \node at (-0.3, 2.7) {$\log_{10} \delta$};
    \node at (3.5, -0.18) {$\log_{10} \radius$};
  \end{tikzpicture}}}}$
  \caption{\label{f:twisted-tori}
    Left: two tori, each with major radius 1 and minor radius $\radius$, with surface separation $\delta$.
    The left ring is centered at the origin and has axis $(0,0,1)$, while the right one is centered at $(2+2\radius+\delta,0,0)$ and has a twisted axis $(0,-\sin \pi/3,\cos \pi/3)$.
  We compute the total drag force $\vct{F}$ on the left ring (held stationary) due to the right ring being translated with the velocity $(0.8,-0.5,1)$ shown by arrows.
  The colors denote the magnitude of the density $\vct{\sigma}$ in the boundary integral solution.
  Right: contour plot of the relative error in the mutual drag force computed using the slender-body theory inverse problem ($\vct{F}_{SBT}$), compared to the fully resolved solution ($\vct{F}_{BIE}$) computed using the boundary integral formulation,
  for different values of $\radius$ and $\delta$.
  Note the $\bigO{1}$ errors for $\delta \le \eps$.
  }
\end{figure}

We turn to a more interesting case in \pr{f:twisted-tori}:
the accuracy of numerical SBT for close-touching rigid bodies
with relative motion.
There are two rings of unit radius, each with cross-sectional radius $\radius$, separated by a minimum surface-to-surface distance $\delta$.
The first ring is held stationary while the second ring is translating with a given velocity, and we compute the total drag force on the first ring.
Recall that SBT expresses velocity in terms of centerline forces,
thus for such a resistance (\eg, rigid body drag) problem
one must invert this to solve for centerline force given centerline velocities.
This amounts to solving a 1D linear integral equation.
\pr{a:sbt} outlines our
high-order Nystr\"om method for this,
which uses the classical (un-regularized)
SBT kernel for self-interactions, and the
correction of \cite{johnson80} for interactions between bodies.
With CSBQ we resolve the solution to about 10 digits, and use it as the reference solution $\vct{F}_{BIE}$.
We plot the relative error in the SBT-approximated force for different values of $\delta$ and $\radius$.
For separations $\delta \gg \radius$ (upper part of the plot),
the error is well explained by $\bigO{\eps^2 \log \eps^{-1}}$ as seen for the single ring in \pr{t:conv-sbt}.
For $\delta=10\radius$ (see dashed line)
we see 2--3 digit accuracy, similar to
that found by Mitchell \etal~in a flexible case \cite{mitchell2022}.
However, the lower-left of the plot indicates about
$10\%$ error when $\delta=\radius$.
Finally, the
net force magnitude using SBT is wrong by greater than a factor of two when $\delta=\radius/10$.

\begin{figure}[htb!]
    \centering
    \includegraphics[width=0.36\textwidth]{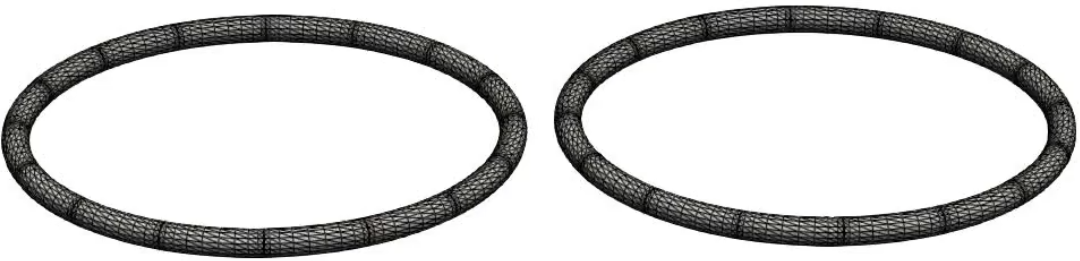}
    \includegraphics[width=0.2107\textwidth]{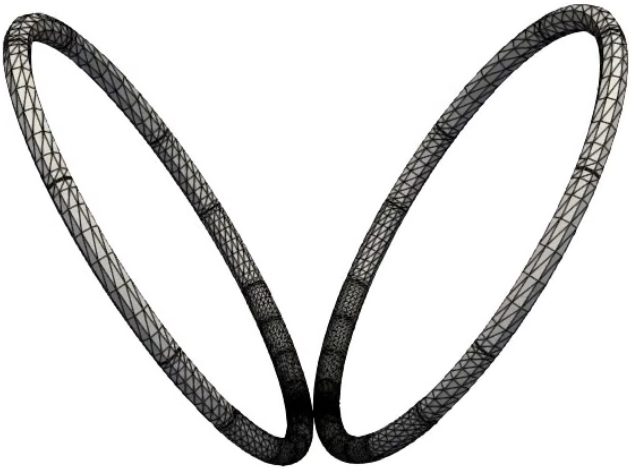}
    \hspace{0.05\textwidth}
    \includegraphics[width=0.059\textwidth]{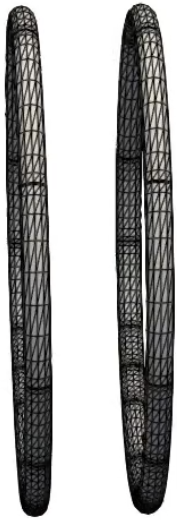}
    \hspace{0.05\textwidth}
    \includegraphics[width=0.1886\textwidth]{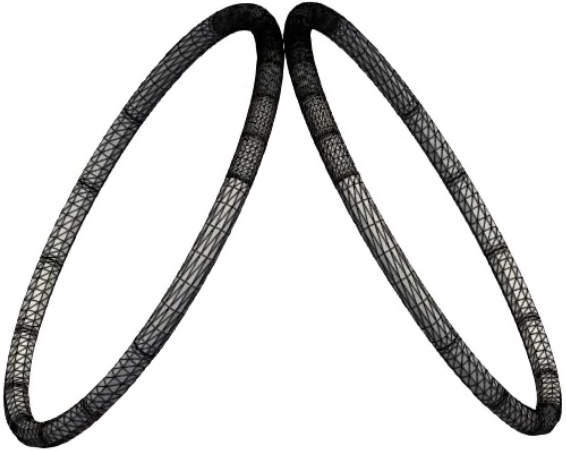}
    \caption{Two rings in a time-periodic sedimentation flow visualized at different points in time.
    At each instant, the velocity and angular velocity of each ring is computed by solving the Stokes mobility problem.
    The solution is resolved to 8-digit accuracy in spatial discretization, quadrature accuracy, and GMRES residual;
    we use a 5\textsuperscript{th}-order adaptive spectral deferred correction (SDC) scheme for time-stepping with a relative error tolerance of $10^{-7}$ per unit time.
    }
    \label{f:sed2}
\end{figure}

\begin{figure}[htb!]
  \resizebox{0.8\textwidth}{!}{\begin{tikzpicture}
    \pgfplotsset{
      xmin=1.0, xmax=305,
      width=0.8\textwidth, height=0.4\textwidth,
      xlabel={$T$}, xtick distance=50,
    }
    \begin{axis}[axis y line*=left, ymin=0, ymax=100, ylabel={$\gmresiter$}, legend style={draw=none,at={(0,1)},anchor=north west}]
      \addplot [thick,color=blue] table [x={t},y={noprecond}] {data/sed2}; \addlegendentry{no-preconditioner};
      \addplot [thick,color=red] table [x={t},y={precond}] {data/sed2}; \addlegendentry{block-preconditioner};
    \end{axis}

    \begin{axis}[axis y line*=right, ymin=0, ymax=62000, ylabel={$\Nunknown$}, legend style={draw=none,at={(0.97,0.97)},anchor=north east}]
      \addplot [thick,dashed,color=black] table [x={t},y={N}] {data/sed2}; \addlegendentry{$\Nunknown$};
    \end{axis}


  \end{tikzpicture}}
  \caption{Number of GMRES iterations $\gmresiter$, and number of unknowns $\Nunknown$ (dotted line and right-side axis), versus time, for sedimentation of two close-to-touching rings as shown in \pr{f:sed2}.
  GMRES iteration counts are shown for the case with (red) and without (blue) a block-diagonal preconditioner.
  The time interval shown (just over 300 units) represents one period of the motion (one full revolution of each ring). This required 149 time steps, with a smallest step size of 1.08 and a largest of 3.89.
  }
  \label{f:sed2-plot}
\end{figure}
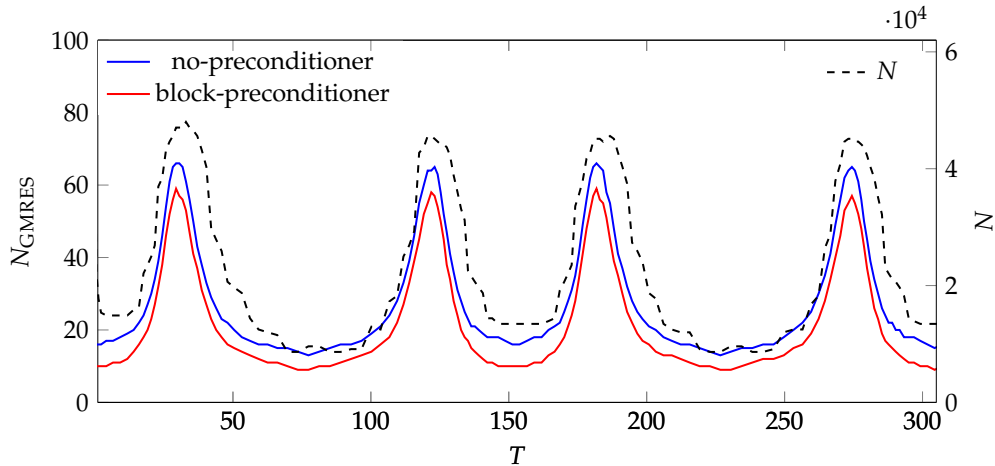

\subsection{Stokes mobility problem\label{ss:results-mobility}}
In \pr{f:sed2}, we show the solution of a mobility problem with two sedimenting rings using the BIE formulation in \pr{e:new-mobility-bie}.
Each ring has a major radius of $0.45$, and $\radius = 0.025$.
The rings are suspended in a Stokesian fluid and a unit downward force is applied to each ring at their center of mass.
We use 5\textsuperscript{th}-order spectral deferred correction (SDC)
\cite{SDC} to evolve in time.
We use adaptive time-stepping with an error tolerance of 1e-7 and use the method of \cite{Quaife2016} to determine the step size.
The smallest and largest step sizes are 1.08 and 3.89 respectively.
The minimum separation between the rings is $0.0014$.
In each time step, we adaptively refine (and coarsen) the geometry to a tolerance of 1e-8.
In \pr{f:sed2-plot}, we show the number of unknowns in the discretization at different points in time.
We also show the number of GMRES iterations required to solve the mobility boundary integral equation with and without using a preconditioner, as follows.
We precompute the exact (dense direct) inverse of the discretized boundary integral operator for one ring at the finest discretization in a reference orientation.
Then, we use this precomputed matrix to construct a block diagonal preconditioner with two diagonal blocks by appropriately rotating and refining the mesh to the reference geometry and back.
The reference geometry used to construct this
preconditioner contains 44K unknowns.
From the plots in \pr{f:sed2-plot}, the preconditioner reduces the number of GMRES iterations by fewer than 10 iterations; it is not very effective.

\begin{figure}[htb!]
    \centering
    \includegraphics[width=0.30\textwidth]{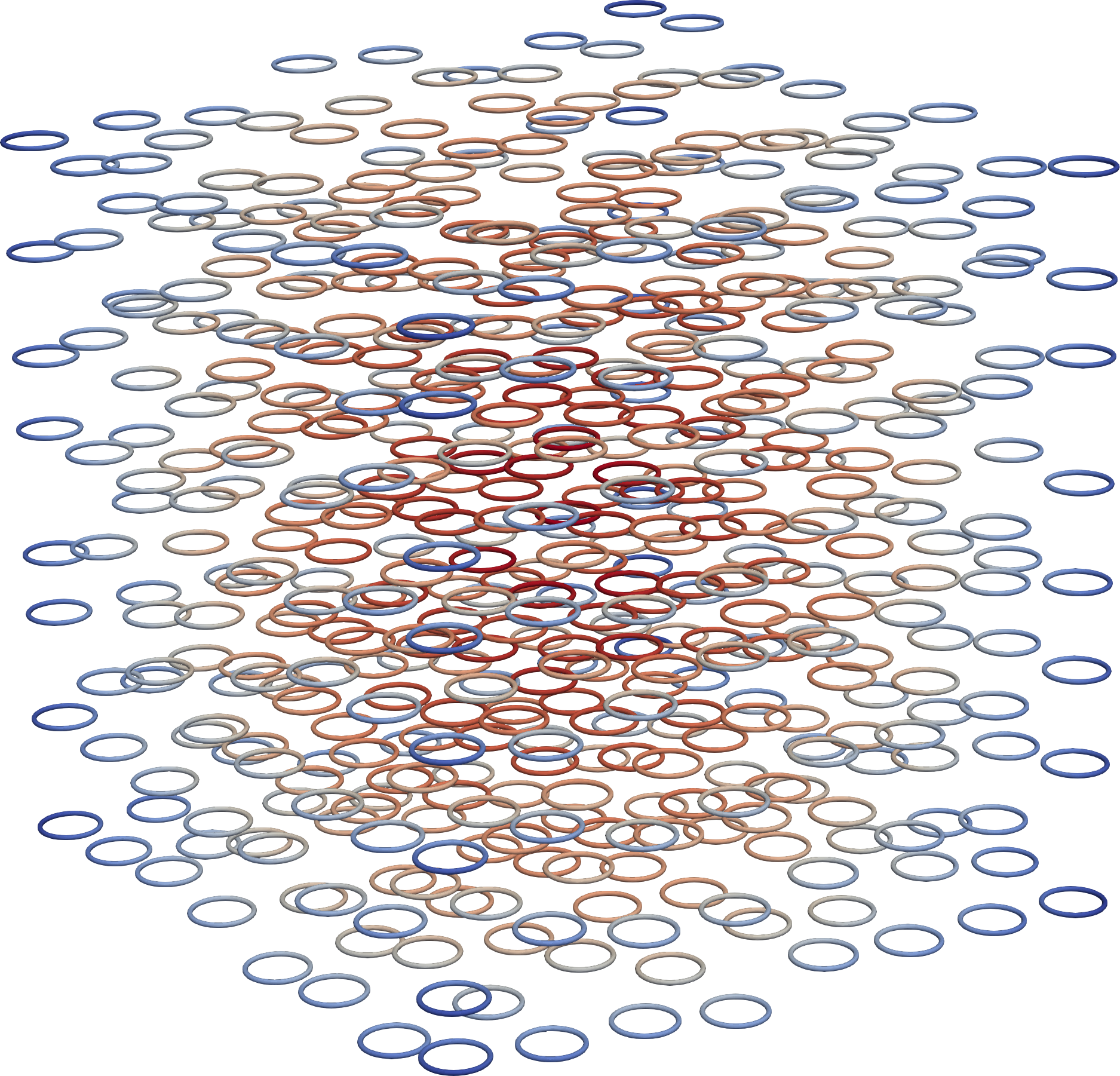}
    \hfill
    \includegraphics[width=0.30\textwidth]{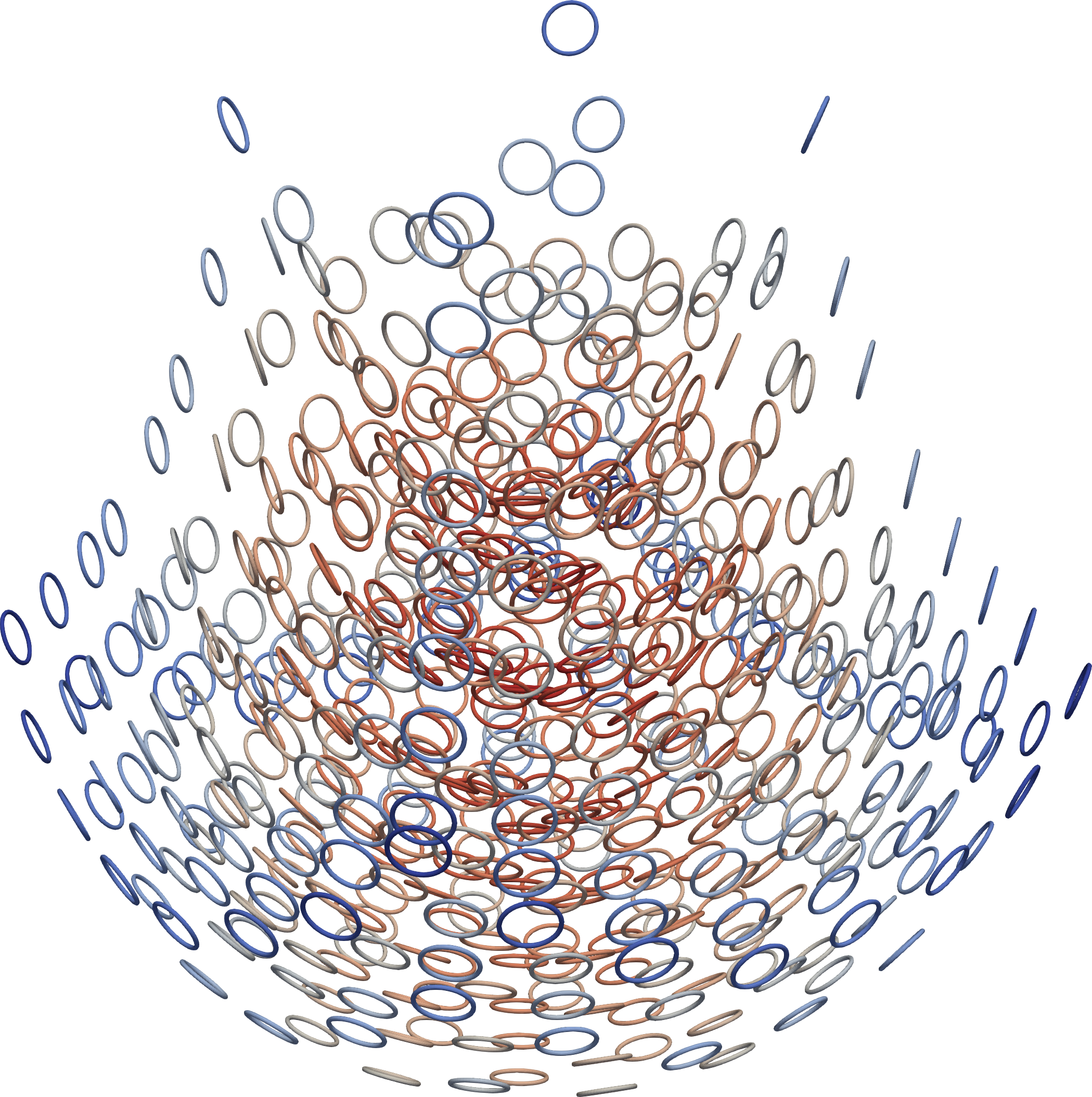}
    \hfill
    \includegraphics[trim={0 0 0 20cm},clip,width=0.25\textwidth]{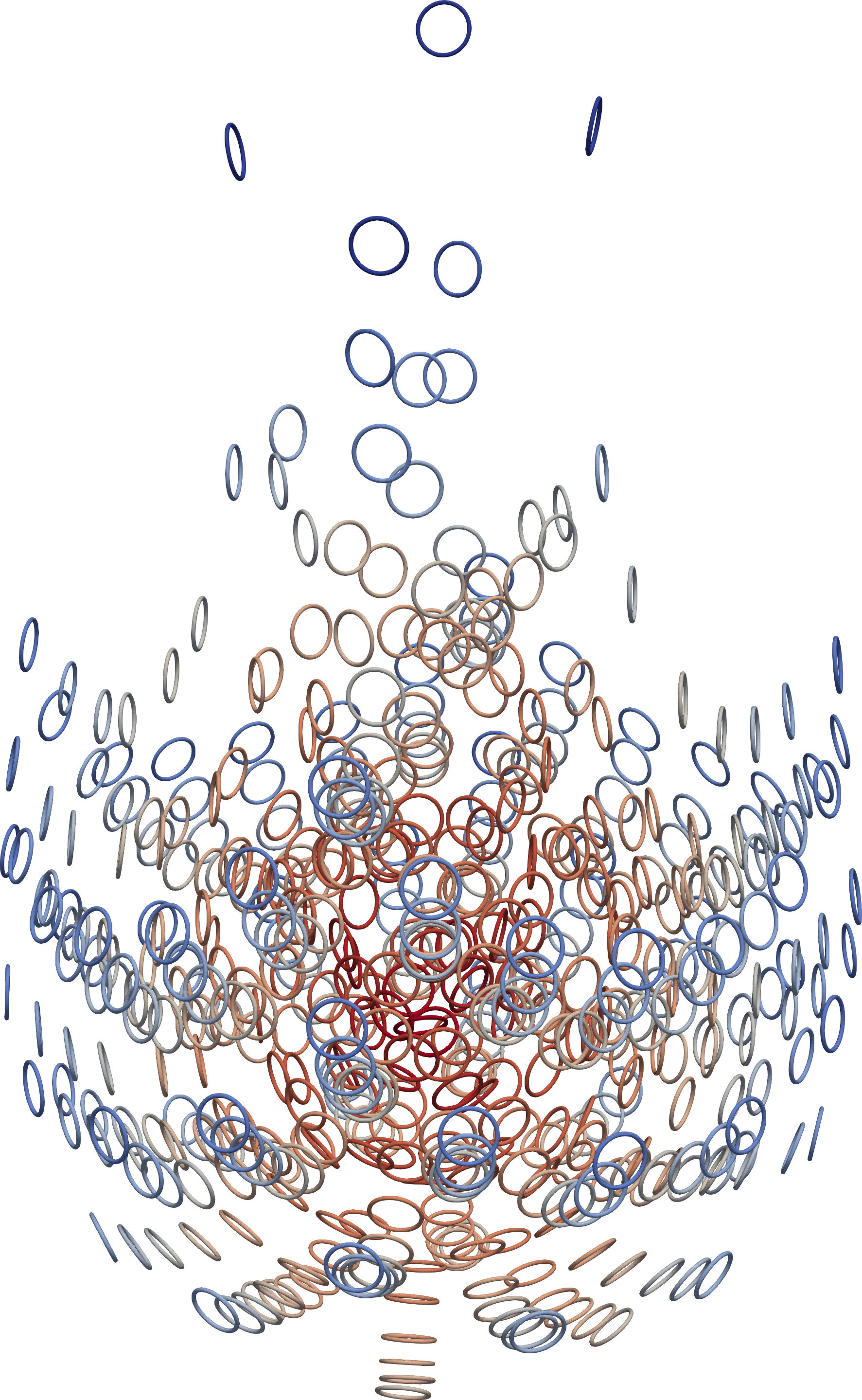}
    \caption{Sedimentation flow with 512 rings, each with major radius $0.45$ and $\radius=0.025$.
      The rings are visualized at times $T = 0$ (left), $6.6$ (center), and $13.2$ (right,
      just before two rings become exponentially close and cause breakdown).
      The solution is computed using 5\textsuperscript{th}-order spectral deferred correction (SDC) time-stepping,
      solved to 7 digits of accuracy in time and 8 digits of accuracy in spatial discretization.
      The color denotes the magnitude of the velocity on $\partial\Omega$.
      In \pr{f:sscal}, we show strong scaling results on $160$ CPU cores for one time step of this flow.}
    \label{f:sed512}
\end{figure}

\begin{figure}[htb!]
  \resizebox{0.75\textwidth}{!}{\begin{tikzpicture}
    \begin{semilogxaxis}
      [xtick=data,xticklabels={1,2,5,10,20,40,80,160,320,640}, xmin=0.707, xmax=230, xlabel=cores,
      ymin=0,ybar stacked, ymajorgrids,ylabel=wall-time $\times$ cores (s),
      legend style={draw=none},legend columns=4, legend entries={$\Tsetup~~$, $\Time_{near}~~$, $\Time_{FMM}~~$, $\Time_{other}$}, legend pos=north west,
      bar width=0.06\textwidth, width=0.8\textwidth, height=0.35\textwidth]

      \addplot[color=black, fill=c1]
      table[x=cores,y=Tsetup] {data/sscal};

      \addplot[color=black, fill=c2]
      table[x=cores,y=Tnear] {data/sscal};

      \addplot[color=black, fill=c3]
      table[x=cores,y=Tfar] {data/sscal};

      \addplot[color=black, fill=gray]
      table[x=cores,y=Tother] {data/sscal};
    \end{semilogxaxis}
\end{tikzpicture}}
  \caption{Strong scaling results on $4$ nodes ($160$ cores total) for one time step (at $T=6.6$) of sedimentation flow in \pr{f:sed512} with $1.65$ million unknowns.
    We report the cumulative CPU time for different stages of the computation.
    $\Tsetup$ is the quadrature setup time for building the local quadrature corrections.
    The quadrature evaluation time $\Teval$ is composed of $\Time_{near}$ for applying the near corrections, and
    $\Time_{FMM}$ for the FMM computation.
    $\Time_{other}$ is the time for the remaining steps in the mobility solve.
  }
  \label{f:sscal}
\end{figure}
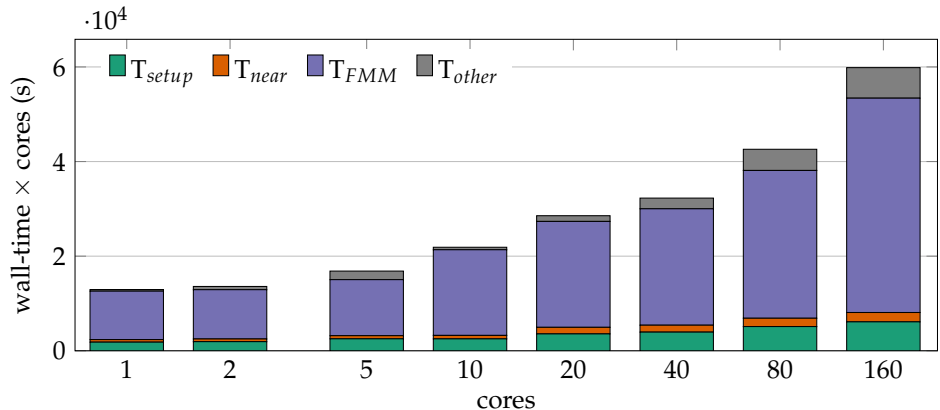

\subsection{Parallel scalability\label{ss:results-scalability}}
We consider the sedimentation flow in \pr{f:sed512} with $512$ slender rings,
each with aspect ratio about 20.
The full simulation took about 2 days of wall-clock time on four nodes
(160 cores); most of this time is devoted to the last 10\% of simulated time
when rings approach each other, $N$ grows, and the time step $\Delta T$
shrinks.
We present a strong scaling study in \pr{f:sscal},
showing a breakdown of the total CPU time for one time step at $T=6.6$ with step size $\Delta T = 0.16$.
This one time step of the 5\textsuperscript{th}-order spectral deferred correction required 18 solves of the Stokes mobility problem,
each of which required on average 16.5 GMRES iterations.
This corresponds to $18$ quadrature setups and $297$ quadrature evaluations each for the single- and double-layer operators.
The computation is overwhelmingly dominated by FMM, comprising $70$--$83\%$ of the total time.
Scaling from $1$ core to $160$ cores (4 nodes), we get a $34\times$ speedup, thus a parallel efficiency of $21.6\%$.
Here, results on up to $40$ cores are for a single node, and the rest are multi-node results using all $40$ cores per node.
We tuned the number of OpenMP threads and MPI processes per node for best performance.
It is preferable to avoid threads accessing data on different NUMA nodes, so it is best to have at least one MPI process per NUMA node (same as a CPU socket in this case);
however, too many MPI processes (\eg, pure MPI) makes the load imbalance worse.

\section{Conclusions}
\label{s:concl}

We have presented an efficient quadrature scheme (CSBQ)
for second-kind boundary integral formulations of rigid 3D slender-body
Stokes Dirichlet and mobility BVPs,
whose cost is independent of the body aspect ratio $\eps^{-1}$.
We have shown that adaptivity along the centerline and in the angular discretization can achieve close to machine accuracy,
even down to lubrication-dominated separations ($<\eps/10$).
This convergent scheme contrasts the commonly used asymptotic
slender body theory (SBT), which is non-convergent at any $\eps$,
has uncontrolled errors, and has $\bigO{1}$
errors for distances of order $\eps$ or less.
We include the first known study (enabled by CSBQ)
of such SBT errors in the setting of close-to-touching rigid bodies.

We have strived for efficiency in implementation,
including using precomputed generalized Chebyshev quadratures
to accelerate near-singular toroidal and singular centerline quadratures.
Our quadrature setup rate is thus about 20,000 unknowns/sec on a single core
at 7-digit accuracy.
We present and test newly-scaled combined-field formulations
for Dirichlet (Laplace and Stokes) and Stokes mobility problems,
with condition number bounded independent of $\eps\to 0$.
We prove that our new (projected) mobility formulation has
a unique solution.
We expect these combined tools to give an efficient and trustworthy alternative to SBT in simulations of viscous rigid-fiber hydrodynamics.
In an HPC distributed-memory parallel FMM-accelerated implementation
we show its use when coupled with an iterative solver and high-order time-stepping for
rigid multi-body sedimentation problems.

\begin{rmk}[Collision handling]
  \label{r:collisions}
  In this work we have chosen to present a high-order convergent scheme for
  what could be called the ``mathematical''
  mobility problem, i.e., integrating velocities and angular velocities given by an
  accurate Stokes BVP solve, no matter how close the surfaces come.
  Yet, in 3D, smooth rigid bodies approach exponentially fast in the lubrication
  limit under constant forces \cite[\S7.1]{KK91book} \cite{guazzellibook}.
  The decay time scales as the typical radius of curvature, which is $\bigO{\eps}$,
  so is extremely short. Thus, once an approach starts,
  miniscule distances of, say, $10^{-100}$ are rapidly
  reached in the mathematical solution. Of course,
  due to surface roughness, friction, and molecular effects \cite{guazzellibook},
  this model no longer matches any physical experiment.
  Numerical breakdown (as in \pr{ss:results-mobility})
  occurs well before this, certainly when
  machine precision fails to distinguish surface coordinates correctly.
  Thus, in practical solvers, explicit short-range force pairs are usually
  added to prevent rigid body collisions \cite{yan20} (although not always for
  flexible fibers \cite{maxian21}).
  Choice of such ad hoc forces is application dependent and
  beyond the scope of this work.
  However, we believe that our work will allow correct hydrodynamics
  to be used with shorter-range collision-avoiding forces than for SBT-based schemes.
\end{rmk}

There are several other fruitful avenues to extend the presented techniques, including:
\bi
\item Open fibers. Here, the centerline panelization would need to be adjusted
  to respect possible endpoint singularities in $s$, when the
  radius does not have parabolic behavior at endpoints.
  In informal studies we have found no difficulties with
  $\eps$-radius hemispherically rounded endpoints,
  and will report results at a later date.
\item
  Non-circular cross-sections, as in \cite{keaveny11,laugabook}.
  The same quadrature rules should work as long as the cross-section curve is not too irregular.
  In evaluating the $\paramt$-integral (\pr{sss:modal-greens-fn}), we would need to also evaluate the surface coordinates,
  normals and Jacobians at each quadrature node from a Fourier representation of the cross-section (as we do now for the surface density).
\item
  Flexible fibers.
  We expect that the ideas presented could accelerate the completed
  single-layer formulation of the slender-body BVP of \cite{mitchell2022},
  or a second-kind formulation of this (non-classical, angle-averaged) BVP.
\item
  Nonuniform adaptivity in the angular direction. This
  could increase efficiency somewhat in the lubrication case of separations less than $\eps$. Likewise, the use of generalized Gaussian
\cite{Bremer2010}
(as opposed to Chebyshev) quadratures could increase quadrature throughput.
\item
  Extension to the alternative single-layer mobility formulation
  involving the interior traction BVP \cite{karrilakim,rachhgreengard}.
  We believe that the question is open:
  Does there exist a single-layer mobility formulation that remains well-conditioned
  for slender bodies as $\eps\to 0$ ?
\ei




\section{Acknowledgments}
We are very grateful to Manas Rachh
for suggesting the projection in the slender Stokes mobility formulation
\pr{e:new-mobility-rep} and help in proving \pr{t:BIE}.
We also thank Laurel Ohm and Ondrej Maxian for their expert help in
understanding SBT, and Mike Shelley for suggesting the problem.
The Flatiron Institute is a division of the Simons Foundation.

\appendix

\gdef\thesection{\Alph{section}} 
\makeatletter
\renewcommand\@seccntformat[1]{\appendixname\ \csname the#1\endcsname.\hspace{0.5em}}
\makeatother

\section{Proofs}
\label{a:pf}

\begin{proof}[Proof of \pr{p:null}]
  The interior Green's representation theorem \cite[(19),~Ch.~3]{Ladyzhenskaya}
  states that if any velocity field $\vct{u}$ with pressure field $p$ satisfies
  the Stokes equations in $\Omega$, then
  $\vct{u} = \StokesSLR[\vct{T}^-] - \StokesDLR[\vct{u}^-]$,
  where the traction vector is
  $\vct{T}(\vct{u},p) := -p\vct{n} + (\nabla\vct{u}+ \nabla\vct{u}^T)\vct{n}$,
  and the ``$-$'' superscript indicates the interior limit on $\partial\Omega$.
  Taking the interior limit, applying jump relations,
  $(I/2+\StokesDL)\vct{u}^- = \StokesSL\vct{T}^-$.
  Choosing $\vct{u}$ a rigid body motion with $p\equiv 0$,
  which are a Stokes solution with $\vct{u}^-\in\RigidBodySpace$,
  the stress and hence traction vanishes,
  showing that $\vct{u}^- \in \mbox{Nul}(I/2+\StokesDL)$.
\end{proof}

\begin{proof}[Proof of \pr{t:BIE}]
  Since $\StokesCF(I-L) + L$ is a compact perturbation of $I/2$,
  Riesz--Fredholm theory
  applies, and it is enough to show that its adjoint is injective. To this end,
  let $\bp$ solve the homogeneous adjoint equation,
  \be
  (I-L)\StokesCF^T\bp + L\bp = \vct{0}.
  \label{orth}
  \ee
  Since the two terms lie in orthogonal spaces, they are both zero.
  In particular, $(I-L)\StokesCF^T\bp = \vct{0}$, and, expanding this
  using $L(I/2+\StokesDL^T) = 0$ which is the adjoint of
  \pr{p:null}, gives $(I/2 + \StokesDL^T)\bp + \eta (I-L)\StokesSL\bp = \vct{0}$.
  Defining the layer potential
  $\vct{u}=\StokesSLR[\bp]$ with $p$ the corresponding pressure,
  the jump relations imply $(I/2 + \StokesDL^T)\bp = \vct{T}^-$ and
  $\StokesSL\bp = \vct{u}^-$, using notation from the proof of \pr{p:null}.
  Thus $(\vct{u},p)$ is a Stokes solution in $\Omega$ with a generalized Robin
  boundary condition
  \be
  \vct{T}^- + \eta (I-L)\vct{u}^- = \vct{0}.
  \label{genBC}
  \ee
  Inserting this into the interior Green's 2nd identity for Stokes solutions
  \cite[p.~53]{Ladyzhenskaya} gives
  \[
  \frac{1}{2}\int_\Omega \sum_{i,j=1}^3
  \biggl( \frac{\partial u_i}{\partial x_j} + \frac{\partial u_j}{\partial x_i} \biggr)^2
  d\xx
  = \int_{\partial\Omega} \vct{u}^- \cdot \vct{T}^- dS
  = -\eta \int_{\partial\Omega} \vct{u}^- \cdot (I-L)\vct{u}^- dS.
  \]
  The left side is nonnegative, but the right side
  nonpositive since $\eta>0$ and $I-L$ is positive-semidefinite.
  Thus both sides vanish, so the stress tensor $\nabla\vct{u}+\nabla\vct{u}^T$
  is identically zero, so $\vct{u}$ is a rigid body motion in $\Omega$,
  so $\vct{u}^- \in \RigidBodySpace$. Again using \eqref{genBC}
  gives $\vct{T}^-=\vct{0}$.
  Thus by the jump relations
  $\vct{u}^+=\vct{u}^-$ and
  $\vct{T}^+=-\bp$,
  the latter obeying $L\bp=\vct{0}$ as shown above from \pr{orth}.
  Thus $\vct{T}^+\in \RigidBodySpace^\perp$,
  making $\int_{\partial\Omega} \vct{u}^+ \cdot \vct{T}^+ dS = 0$.
  Now applying the exterior Green's 2nd identity gives that stress of $\vct{u}$ is identically
  zero in the exterior, so that the corresponding pressure field $p$ must be constant.
  Yet this constant must be zero since $p=\bigO{1/|\xx|^2}$ for a single-layer potential.
  Thus $\vct{T}^+=\vct{0}$, so $\bp=\vct{0}$, proving the desired injectivity.
\end{proof}

\section{Generalized Chebyshev Quadratures \label{s:cheb-quad}}
\newcommand{\Nfun}{{\ensuremath{k}}} 
\newcommand{\Ngl}{{\ensuremath{m}}} 
\newcommand{\Nq}{{\ensuremath{n}}} 

We briefly describe our implementation of the algorithm of
\cite[Sec.~4]{Bremer2010}
for constructing a Chebyshev quadrature rule for a given set of $\Nfun$ integrand functions,
over a given interval $[a,b]\subset\Real$, to a desired tolerance $\epsilon$.
The algorithm first compresses the $\Nfun$ integrands to $\Nq \le \Nfun$ orthogonal basis functions, then constructs an $\Nq$-point Chebyshev quadrature rule for these basis functions.
The result is nodes $\{y_i\}_{i=1}^n$ and weights $\{\omega_i\}_{i=1}^n$.
These steps are described below.

\paragraph{Compression of integrands to an orthonormal basis}
Given an accuracy tolerance $\epsilon$, and $\Nfun$ piecewise continuous square-integrable integrand functions $\{ \varphi_1, \varphi_2, \dots, \varphi_\Nfun \}$ on $[a,b]$,
we first need a quadrature rule on $[a,b]$ with nodes $\{x_1, \dots, x_\Ngl\}$ and weights $\{w_1, \dots, w_\Ngl\}$
which integrates products of any two integrands to the required tolerance.
A simple way to do this is to i)
adaptivity construct a piecewise degree-$(p-1)$ polynomial interpolant for the vector valued function $\vct{\varphi} = (\varphi_1, \varphi_2, \dots, \varphi_\Nfun)$ to a relative accuracy $\epsilon$, using the usual algorithm
of sampling at $p$ Chebychev points per panel and splitting panels
until their Chebyshev coefficient magnitudes decay to $\epsilon$.
A predetermined dyadically refined panelization may instead be used for
a known singularity location.
Then, ii) build a composite quadrature rule with the $p$-node Gauss--Legendre rule in each of the panels. Since the degree that GL integrates exactly
is $2p-1 > 2(p-1)$, this is accurate for all products of function pairs.

Given the above $m$-node rule for $[a,b]$, we construct the following
$m$-by-$k$ matrix
to which compression will be applied:
\[
A =
\begin{pmatrix}
\varphi_1(x_1)\sqrt{w_1} & \varphi_2(x_1)\sqrt{w_1} & \cdots & \varphi_\Nfun(x_1)\sqrt{w_1} \\
\varphi_1(x_2)\sqrt{w_2} & \varphi_2(x_2)\sqrt{w_2} & \cdots & \varphi_\Nfun(x_2)\sqrt{w_2} \\
\vdots & \vdots & & \vdots \\
\varphi_1(x_\Ngl)\sqrt{w_\Ngl} & \varphi_2(x_\Ngl)\sqrt{w_\Ngl} & \cdots & \varphi_\Nfun(x_\Ngl)\sqrt{w_\Ngl} \\
\end{pmatrix} .
\]
The scaling by the square-root of the quadrature weights means that the dot-product of any two columns of A is equal to the $L^2$-inner product of the corresponding functions since
$\int_a^b \varphi_i \varphi_j dx = \sum_{k=1}^{\Ngl} \varphi_i(x_k) \varphi_j(x_k) w_k$.
Given a tolerance $\epsilon$, we construct an orthonormal basis for the columns of $A$,
which can approximate any column of $A$ to an accuracy of $\epsilon$ in the $l^2$-norm.
This can be done by computing a singular value decomposition ($A = U \Sigma V^T$)
then
truncating the $U$ matrix to the first $\Nq$ columns such that the singular values corresponding to the discarded columns are smaller than $\epsilon$.
(A rank-revealing QR decomposition $A = U \widetilde{R}$ may similarly be truncated.)
This is equivalent to approximating an $L^2([a,b])$-orthonormal basis
$\{u_1, \dots, u_\Nq\}$ for Span $\{\varphi_1, \dots, \varphi_\Nfun\}$,
to an accuracy of $\epsilon$,
and then discretizing these basis functions using the above
adaptive quadrature rule.
Therefore the truncated $m$-by-$n$ matrix is
\[
U =
\begin{pmatrix}
u_1(x_1)\sqrt{w_1} & u_2(x_1)\sqrt{w_1} & \cdots & u_\Nq(x_1)\sqrt{w_1} \\
u_1(x_2)\sqrt{w_2} & u_2(x_2)\sqrt{w_2} & \cdots & u_\Nq(x_2)\sqrt{w_2} \\
\vdots & \vdots & & \vdots \\
u_1(x_\Ngl)\sqrt{w_\Ngl} & u_2(x_\Ngl)\sqrt{w_\Ngl} & \cdots & u_\Nq(x_\Ngl)\sqrt{w_\Ngl} \\
\end{pmatrix} .
\]

\paragraph{Finding stable quadrature nodes and weights}
We now build an $\Nq$-point quadrature rule that integrates each of $\{u_1, \dots, u_\Nq\}$.
We first compute a column pivoted QR decomposition of $U^T$.
Each pivot column corresponds to a node in $\{x_1, \dots, x_\Ngl\}$ and
we take these $\Nq$ nodes to be our quadrature nodes
$\{y_1, \dots, y_\Nq\} \subset \{x_1, \dots, x_\Ngl\}$.
The quadrature weights $\omega_j$
are obtained by solving the square ``Vandermonde transpose'' linear system,
\begin{equation*}
\sum\limits_{j=1}^{\Nq} u_i(y_j) \omega_j = \int_a^b u_i dx, \qquad \text{for all}~~ 1 \leq i \leq \Nq.
\end{equation*}
The use of column pivoting insures
system has a small condition number and
that the quadrature weights can be stably computed
\cite[Thm.~3.2]{Bremer2010}.
In all of the above precomputations we use quad (128-bit real) precision,
although double precision is often adequate for
$\epsilon > 10^{-10}$.
Note that further optimization to generalized Gaussian
rules with less nodes is possible \cite{Bremer2010}; however,
we favored the above since it is completely automatic and
performed well enough.

\section{Numerical solution of the slender-body theory inverse problem}
\label{a:sbt}

Here we outline a high-order accurate discretization of 
the periodic version \cite{ueda} of
classical nonlocal SBT \cite{johnson80,gotz00},
and a solution of the resulting 1D integral equation needed for
\pr{ss:results-sbt}.
We consider multiple
fibers with constant circular cross-sectional radius $\radius$,
and no background flow.
For notational simplicity,
each closed centerline curve $\gamma_b$, $b=1,\dots,\Nobj$ has the same length
$L$.
SBT is most conveniently expressed in terms of arc-length,
so we assume here that $\xx_b(s)$ is an arc-length parameterization
of the $b$th centerline curve
(note that the body index replaces the subscript $c$ used in \pr{ss:disc}).
Numerically, any smooth parameterization may be converted to arc-length
using panel-wise antiderivatives of its speed function.
Then $\gamma_b = \xx_b([0,L))$.
Classical SBT expresses velocity as a line integral over a given
force density $\vct{f}$ on all centerlines.
Precisely, letting $\shat = d\xx_b(s)/ds$ be the unit tangent at $s$,
the velocity $\uu_b$ and force density $\vct{f}_b$
on the centerline of the $b$th body are related by
\be
\uu_b(s) =
\frac{1}{8\pi}\left[(I - 3\shat\shat^T) - 2(I+\shat\shat^T)
  \log \frac{\pi\radius}{4L}      
  \right] \vct{f}_b(s)
+ \sum_{b'=1}^\Nobj K_{bb'}[\vct{f}_{b'}](s).
\label{SBT}
\ee
Here the first term is a local drag term, whereas the second
$\radius$-independent term is nonlocal and defined for the self-interaction of the
$b$th fiber (diagonal blocks) by
\be
K_{bb}[\vct{f}_{b}](s) :=
\int_{0}^{L_{b}}
\left(
  S(\vct{r}) \vct{f}_b(s') - 
  \frac{I + \shat\shat^T}{8\pi d(s,s')} \vct{f}_b(s)
  \right)
\,ds',
\qquad \mbox{ with }\;
d(s,s') = \frac{L}{\pi}\sin\biggl(\frac{\pi}{L}|s-s'|\biggr),
\label{SBTself}
\ee
where $\vct{r}:= \xx_b(s)-\xx_b(s')$ is the displacement vector,
and $S$ is the Stokeslet velocity kernel defined beneath \pr{e:stokes-sl}.
The between-fiber interactions (off-diagonal blocks $b'\neq b$) are
\be
K_{bb'}[\vct{f}_{b'}](s) :=
\int_{0}^{L_{b'}}
\left[S(\vct{r}) + \frac{\radius^2}{2}\Delta S(\vct{r})\right]
\vct{f}_{b'}(s') \,ds' ,
\label{SBTfar}
\ee
where the between-fiber displacement is $\vct{r}:= \xx_b(s)-\xx_{b'}(s')$,
and the 
so-called ``doublet''
\footnote{The naming is somewhat confusing in the literature:
it is also known variously as a
``potential dipole'' \cite[(7.2.3)]{Pozrikidis1992},
``degenerate quadrupole'' \cite[(10.62)]{KK91book},
or ``doublet flow'' since
$\Delta S(\xx-\yy)\vct{f}$ is the flow generated by the limit of
a source and sink separated infinitesimally in the $\vct{f}$ direction.
It should not be confused with the force or ``Stokeslet'' doublet
\cite[(7.2.15)]{Pozrikidis1992}.}
kernel
\cite[(3.3.8)]{Pozrikidis1992}
is
$\frac{1}{2}\Delta S(\vct{r}) = -\frac{1}{8\pi}\nabla\nabla |\vct{r}|^{-1} =
\frac{1}{8\pi}\left( \frac{\vct{I}}{|\vct{r}|^3} - 3\frac{\vct{r} \vct{r}^{T}}{|\vct{r}|^5} \right)$.

To devise a high-order Nystr\"om quadrature, the diagonal singularity of
the self-interaction kernel must be understood.
The two terms in \pr{SBTself}
are each divergent like $|s-s'|^{-1}$,
thus neither is integrable, but the
periodized arc-distance function $d(s,s')$ is such that their
singularities cancel.
It can be shown
that for a smooth curve and force density
the integrand is smooth
apart from a discontinuity at $s'=s$, whose size is related
to curvature and to $d\vct{f}_b/ds$ \cite[Sec.~3.1]{tornbergprod}.
For example, for the unit circle in the
$xy$-plane with force density $\vct{f}=(0,0,f)$, the
$z$-component of \pr{SBTself} is
$\int_0^{2\pi} (f(s')-f(s))/|e^{is'}-e^{is}|\,ds'$, whose
integrand jumps from $-f'(s)$ to $f'(s)$ at $s'=s$.

Force and velocity are discretized on the centerline using panel-based
Gauss--Legendre quadrature.
The centerline parameter is split into
intervals $\{I_1,\dots, I_\Nelem\}$, whose union is $[0,L)$.
Let $w_{i}^{(k)}$ be the GL weights and
$s_i^{(k)}$ the nodes, $i=1,\dots,\Ns$, for GL quadrature on $I_k$.
It was sufficient for our experiments to fix
$\Nelem$ equal-length panels for each body, and fix an order $\Ns=12$.
The nodes map to points $\yy_{b,i}^{(k)} = \xx_b(s_i^{(k)})$
with force density samples $\vct{f}_{b,i}^{(k)}$ and velocity samples
$\vct{u}_{b,i}^{(k)}$.
Between-fiber blocks are discretized using plain Nystr\"om quadrature, so
\[
K_{bb'}[\vct{f}_{b'}](s) \approx
\sum_{k=1}^\Nelem \sum_{i=1}^{\Ns}
\left[ S(\xx_{b}(s) - \yy_{b',i}^{(k)})
  + \frac{\radius^2}{2} \Delta S(\xx_{b}(s) - \yy_{b',i}^{(k)})
  \right] w_i^{(k)}
\vct{f}_{b',i}^{(k)}
.
\]
By setting $s$ to each of the $\Nelem\Ns$ nodes for centerline $b$,
this defines elements of each offdiagonal
block $A_{bb'}$ of the overall Nystr\"om matrix
$A$.
Unlike in \pr{ss:quad}, we do not use near-singular corrections, so
must push $\Nelem$ high enough so that all panels are in each
other's far fields.
Blocks of self-interaction matrices $A_{bb}$ between different
panels are filled similarly using plain Nystr\"om quadrature for
\pr{SBTself}; note here that the 2nd term involving $\vct{f}_b(s)$
subtracts only diagonal entries given by the row-sums of the Nystr\"om
matrix for $(I+\shat\shat^T)/8\pi d(s,s')$.

A special rule is needed to handle the self-interation of each
panel in \pr{SBTself}.
Two options are product quadratures \cite{tornbergprod} or auxiliary
nodes; for simplicity we choose the latter.
(Note that, unlike with algebraic singularities,
no special rule is needed for neighboring panels.)
Let the panel parameter interval $I_k$ be $[a,b)$, and consider
a target parameter $s\in I_k$.
Define auxiliary nodes as the union of $\Ns$ GL nodes for
$(a,s)$ with $\Ns$ nodes for $(s,b)$, and define corresponding weights.
This split at the target point handles the diagonal discontinuity.
We interpolate from the original GL nodes on $I_k$ onto these auxiliary nodes
using a barycentric Lagrange matrix $P$ as in \pr{sss:near-sing-quad}.
Finally, the local drag term (first term in \pr{SBT}) simply adds to the
diagonal of the Nystr\"om matrix.

We fill the $3N\times 3N$ Nystr\"om matrix $A$ as above,
where $N$ is the total number of nodes on all fibers,
thus get the discretization of \pr{SBT} as
\be
A \vct{F} = \vct{U},
\label{SBTlinsys}
\ee
where $\vct{U}, \vct{F} \in \Real^{3N}$
are the vectors of velocities and force densities at all nodes.
The resistance problem specifies $\vct{U}$, thus we
solve \eqref{SBTlinsys} for $\vct{F}$
(the SBT inverse problem \cite{mori21inv}), then
sum $\vct{F}$ using the quadrature weights to get the total force (drag)
on each body.
For convergence to 7 digit accuracy in the tests of \pr{ss:results-sbt},
we found $\Nelem \ge \mbox{max}\bigl(50,3.8/(\delta+2\radius)\bigr)$
panels sufficient, requiring up to $\Nelem \approx 400$.
A dense direct solution of \eqref{SBTlinsys} was then adequate for this task.
This completes the high-order accurate SBT numerical solution.
Clearly, many efficiency gains would be possible, but are beyond the needs
of this paper.

\begin{rmk}
  The classical SBT self-interaction operator $K_{bb}:C(\gamma_b)\to C(\gamma_b)$
  in \pr{SBTself} is unbounded,
  with negative eigenvalues growing logarithmically
  in magnitude \cite{gotz00,ueda}, causing the spectrum of the total
  SBT operator on the right side of \pr{SBT} to pass close to zero.
  This may cause resonance, or even lack of invertibility,
  but is only numerically relevant when $N \gtrsim \radius^{-1}$,
  \ie, for larger $\radius$ values.
  Our high-order Nystr\"om discretization allowed convergence to
  around 7 digits even though resonances (oscillatory $\vct{f}(s)$)
  were sometimes visible at the largest $\radius$.
  Regularized versions (either by mollifying the kernel
  \cite{Tornberg2004} or by deleting a small interval in $s'$
  about the target point \cite{maxian21})
  have been proposed which avoid the resonance problem; yet, since
  they cause a change in the right side of \pr{SBT} of $\bigO{\radius^2}$,
  we do not expect them to give different conclusions in
  \pr{ss:results-sbt}.
  A comparison of different approaches to regularizing SBT is
  beyond the scope of this paper.
\end{rmk}


The documented MATLAB implementation used is found in the {\tt SBT}
directory of the {\tt CSBQ} repository discussed in \pr{s:results}.

\bibliographystyle{elsarticle-num-names}
\bibliography{ref}

\end{document}